\documentclass[11pt]{amsart}
\usepackage{amssymb}
\usepackage{graphicx}
\usepackage{xcolor} 
\usepackage{tensor}
\usepackage{fullpage} 
\usepackage{amsmath}
\usepackage{amsthm}
\usepackage{mathrsfs} 
\usepackage{verbatim}

\usepackage[colorlinks=true,citecolor=blue,linkcolor=blue]{hyperref}
\usepackage{esint}
\usepackage{enumitem}

\usepackage{array}
\usepackage{multirow}
\usepackage{makecell}
\usepackage{pdflscape}
\usepackage{caption}

\usepackage{tabularx}
\newcolumntype{Y}{>{\centering\arraybackslash}X}

\definecolor{purple}{rgb}{0.5, 0, 1}
\definecolor{orange}{rgb}{1,.5,0}

\providecommand{\MR}{\relax\ifhmode\unskip\space\fi MR }

\usepackage{seqsplit} 
\definecolor{green}{rgb}{0,0.8,0} 

\newcommand{\redregion}{%
  \begingroup
  \setlength{\unitlength}{1pt}%
  \begin{picture}(14,9)
    \put(1,1){\color{red!12}\rule{12pt}{7pt}}%
    \put(1,1){\color{red!75!black}\dashbox{1}(12,7){}}%
  \end{picture}%
  \endgroup
}

\newtheorem{theorem}{Theorem}[section]
\newtheorem{corollary}[theorem]{Corollary}
\newtheorem{lemma}[theorem]{Lemma}
\newtheorem{proposition}[theorem]{Proposition}

\theoremstyle{definition}

\theoremstyle{definition}
\newtheorem{remark}[theorem]{Remark}

\newtheorem{qu}[theorem]{\textbf{Question}}

\numberwithin{equation}{section}

\newcommand{\nnrm}[1]{{\vert\kern-0.25ex\vert\kern-0.25ex\vert #1 
		\vert\kern-0.25ex\vert\kern-0.25ex\vert}}

\newcommand{\dist}{\mathrm{dist}\,}

\newcommand{\supp}{{\mathrm{supp}}\,}

\newcommand{\bbR}{\mathbb R}

\newcommand{\bfone}{\mathbf1}

\begin{document}
	

\title{The Kelvin--Hicks Formula from the Energy--Casimir Principle}
    
	\author{Ken Abe}
	\address{Department of Mathematics, Graduate School of Science, Osaka Metropolitan University, 3-3-138 Sugimoto, Sumiyoshi-ku, Osaka, 558-8585, Japan}
	\email{kabe@omu.ac.jp}
	
	\author{Kwan Woo}
	\address{Departement Mathematik und Informatik, Universit\"at Basel, CH-4051 Basel, Switzerland}
	\email{kwan.woo@unibas.ch}
	
	\date{\today}

\renewcommand{\thefootnote}{\fnsymbol{footnote}}
\footnotetext[1]{%
}


\begin{abstract}
The classical Kelvin--Hicks formula expresses the propagation speed of a thin vortex ring in an incompressible ideal fluid in terms of its circulation, distance from the symmetry axis, and core radius.
In this paper, we rigorously derive this formula from a variational principle based on Arnold's energy--Casimir method, without prescribing any of these geometric or physical quantities a priori.
We introduce a scale-invariant variational principle involving the energy, impulse, and $L^p$ Casimir invariants. For $1 <  p  < \infty$, we construct orbitally stable energy-maximizing vortex rings and investigate how their geometry depends on the Casimir exponent $p$.
For sufficiently large $p$, the vortex rings have cores homeomorphic to a ball and converge to Hill's spherical vortex as $p \to \infty$.
In contrast, as $p\to1^+$, the energy-maximizing vortex rings concentrate onto a circular filament of radius $e^{1/3}$, uniquely selected by the variational principle.
After suitable translations and rescalings, their vortex-core profiles converge to an optimizer of the logarithmic Hardy--Littlewood--Sobolev inequality.

Introducing an effective core radius $\epsilon_p$ determined by the entropy of the vorticity distribution, we establish the Kelvin--Hicks formula
\begin{align*}
W_p =  \frac{\Gamma_p}{4\pi R_p} \left(\log\frac{8R_p}{\epsilon_p}-1+o(1)\right) \quad \textrm{as}\ p\to1^+,
\end{align*}
where $W_p$, $\Gamma_p$, and $R_p$ denote the propagation speed, circulation, and effective ring radius, respectively.
In particular, the constant correction $-1$ is obtained for nonuniform vortex-core profiles without imposing any prescribed core structure.
\end{abstract}

\maketitle

\tableofcontents

    \clearpage
    
    \section{Introduction}\label{sec:int}

\subsection{Kelvin--Hicks formula}
\begin{figure}[htbp]
\centering
\includegraphics[width=0.9\textwidth]{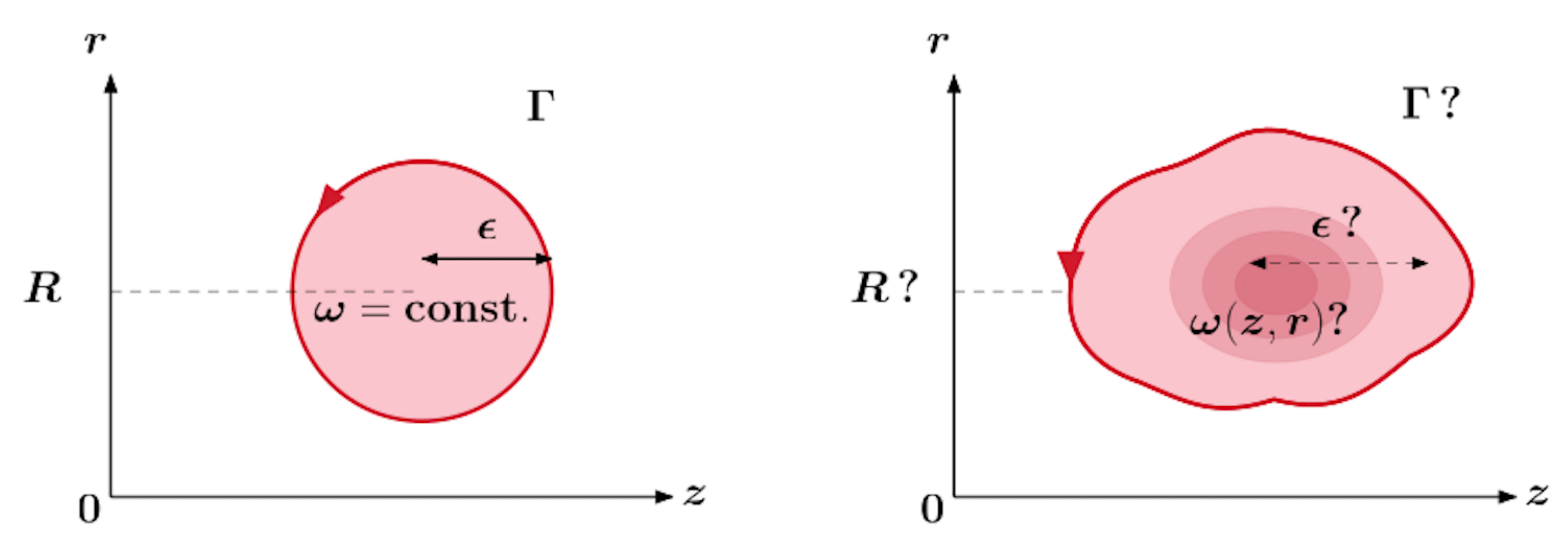}
\caption{
Comparison of idealized and general vortex rings.
\textbf{Left:} An idealized vortex ring to which the Kelvin--Hicks formula applies, with prescribed circulation $\Gamma$, ring radius $R$, core radius $\epsilon$, and vortex-core profile $\omega$.
\textbf{Right:} A general vortex ring for which these quantities are not prescribed a priori and whose core need not be circular.
}
\label{fig:vortex_ring0}
\end{figure}

The propagation speed $W$ of a thin axisymmetric vortex ring translating through an incompressible ideal fluid is described by the Kelvin--Hicks formula \cite[art. 163 (7), p. 241]{Lamb}, \cite[\S 10.2 (1)]{Saff}, \cite[(1.3)]{FM2000}: 
\begin{align}\label{eq:Kelvin-Hicks}
W = \frac{\Gamma}{4\pi R} \left( \log\frac{8R}{\epsilon} - \frac{1}{4} + o(1) \right) \quad \text{as } \frac{\epsilon}{R}\to0,
\end{align}
where $\Gamma$ is the circulation, $R$ is the major radius of the vortex ring, and $\epsilon$ is the radius of its cross-sectional core; see Figure \ref{fig:vortex_ring0}.
This formula was derived by Kelvin \cite{Kelvin1867} and Hicks \cite{Hicks1885} through a formal asymptotic expansion of the Biot--Savart law.
\footnote[2]{The leading-order logarithmic dependence in the Kelvin--Hicks formula can already be traced back to Helmholtz's original work on circular vortex filaments \cite{Hel58}; see also \cite[pp.~781--782]{MGK12} for a historical account.}
The Kelvin--Hicks formula provides a quantitative description of the singular self-induced velocity of a vortex filament in terms of the parameters $\Gamma$, $R$, and $\epsilon$.
It also depends on the vortex-core profile $\omega$, as the constant correction $-1/4$ reflects the assumption of a uniform vorticity distribution $\omega \equiv \textrm{const.}$ \cite{Saf70}.

Fraenkel \cite{Fra70,Fra72} investigated the connection between the Kelvin--Hicks formula \eqref{eq:Kelvin-Hicks} and steadily translating solutions of the Euler equations through a variational construction.
More recently, de Valeriola and Van Schaftingen \cite{dVVS13} constructed steadily translating vortex rings whose propagation speeds exhibit the leading-order logarithmic behavior of the Kelvin--Hicks formula \eqref{eq:Kelvin-Hicks} for prescribed parameters $\Gamma>0$, $R>0$, and sufficiently small $\epsilon>0$.
Meyer and Seis \cite{MeyerSeis2026} constructed steadily translating vortex rings with interfacial vortex sheets in two-fluid flows with surface tension and derived a generalized Kelvin--Hicks formula for their propagation speeds.
Beyond steadily translating solutions, the propagation of thin-core vortex rings has been studied by Benedetto, Caglioti, and Marchioro \cite{BCM2000} and Guo, Jeong, and Zhao \cite{GJZ2026} for inviscid flows, and by Gallay and \v{S}ver\'ak \cite{GS24b} for viscous flows.

In this paper, we derive the Kelvin--Hicks formula with a constant correction corresponding to a nonuniform vortex-core profile for steadily translating vortex rings of the Euler equations, without prescribing their circulation $\Gamma$, ring radius $R$, core size $\epsilon$, or vortex-core profile $\omega$.

We address this question through an energy--Casimir variational principle. In experiments, vortex rings are typically generated by a piston ejecting water from a cylindrical tube.
As the stroke ratio increases, the leading vortex ring grows until a critical value, called the formation number \cite{GRS98}, beyond which a trailing jet develops; see also \cite{Dabiri09,LLXG25}. These experiments highlight the importance of physical quantities such as kinetic energy $E$, impulse $P$, and circulation $\Gamma$ in understanding vortex-ring formation. In particular, experiments \cite[\S\S5.5 and 6.2]{LN21} suggest that the dimensionless energy
\begin{equation}\label{eq:nondim}
\mathscr E_{\mathrm{phys}} = \frac{E}{P^{1/2}\Gamma^{3/2}}
\end{equation}
takes similar values for isolated vortex rings.

On the other hand, the axisymmetric Euler equations without swirl possess an infinite family of conserved quantities known as Casimir invariants, in addition to energy, impulse, and circulation.
These invariants are fundamental conserved quantities in ideal fluid dynamics and play an important role in the study of steady states and their stability; see, e.g., \cite[II.4]{AK21}, \cite[\S\S 3--4]{Elg26}.
Since Casimir invariants encode information about the vorticity distribution, they provide natural constraints for determining the internal structure of vortex rings.

Motivated by the Kelvin--Hicks formula \eqref{eq:Kelvin-Hicks} and by vortex ring formation experiments, we investigate the role of Casimir invariants in selecting vortex-ring geometry for three-dimensional axisymmetric Euler flows without swirl.
In this setting, the relative vorticity admits a family of conserved $L^p$ Casimirs $Q_p$, with $Q_1=2\pi\Gamma$ for nonnegative vorticity.
Replacing the circulation $\Gamma$ in \eqref{eq:nondim} by the Casimir $Q_p$, we introduce the scale-invariant energy functional
\begin{equation}\label{eq:nondim2}
\mathscr{E}_p 
=
\frac{E}{P^{\alpha_p} Q_p^{\beta_p}},
\end{equation}
where
\begin{equation}\label{eq:exponent}
\alpha_p = \frac{7p-6}{5p-3},
\qquad 
\beta_p = \frac{3p}{5p-3}.
\end{equation}
The precise definitions of $E$, $P$, $\Gamma$, and $Q_p$ can be found in \S\ref{s:1.2}.  

\begin{qu}\label{q:q1}
Do maximizers of the scale-invariant energy functional \eqref{eq:nondim2} corresponds to steadily translating vortex rings?
If so, how does the Casimir exponent $p$ determine their geometry, core structure, and propagation speed?
\end{qu}

\subsection{Energy--Casimir principle}\label{s:1.2}
To address Question \ref{q:q1}, we formulate an energy maximization problem for axisymmetric Euler flows without swirl.
The three-dimensional incompressible Euler equations take the form
\begin{equation}\begin{aligned}
u_t+u\cdot \nabla u+\nabla q&=0,\\
\nabla \cdot u &=0,
\end{aligned}
\label{eq:Euler0}
\end{equation}
where $u$ is the velocity field and $q$ is the pressure.
For axisymmetric flows without swirl, we use cylindrical coordinates $(r,\theta,z)$ and express the velocity and vorticity fields in the Clebsch representation
\[
u= \nabla \times (\psi\nabla \theta),\quad
\nabla \times u=r^{2}\zeta \nabla \theta,
\]
where $\psi$ is the axisymmetric stream function and $\zeta$ is the relative vorticity.
Taking the curl of the Clebsch representation, we obtain
\[
\nabla \times u= -L\psi \nabla \theta,
\]
where $L= \Delta_{z,r} - r^{-1}\partial_r$ is the Grad--Shafranov operator.
Thus, the vorticity and stream function are related through the elliptic boundary value problem in the half-plane $\mathbb{R}^{2}_{+} = \{ (z,r) \in \bbR^2: r > 0\}$:
\begin{equation}
- L\psi = r^{2}\zeta \quad \textrm{in } \mathbb{R}^{2}_{+},
\qquad
\psi = 0
\quad \textrm{on } \partial\mathbb{R}^{2}_{+}. \label{eq:elliptichalf}
\end{equation}
Let $G/r'=G(z,r,z',r')/r'$ denote the Green kernel of $(-L)^{-1}$.
The axisymmetric Biot--Savart law then takes the form
\begin{align}\label{eq:streamfunction}
\psi(z,r) = \int_{\mathbb{R}^{2}_{+}} G(z,r,z',r') \zeta(z',r') r'dz'dr' =:\mathcal{G}\zeta,
\end{align}
where
\begin{equation} \label{eq:Green}
G(z,r,z',r') = \frac{\sqrt{rr'}}{2\pi} F(s) = \frac{\sqrt{rr'}}{2\pi} \int_0^{\pi} \frac{\cos \theta} {\sqrt{2(1 - \cos \theta) + s}} d\theta,
\quad
s = \frac{|z - z'|^2+|r-r'|^2}{rr'}.
\end{equation}
The vorticity equation associated with \eqref{eq:Euler0} can then be written as the active scalar equation
\begin{equation}\label{eq:Euler1}
\zeta_t +  u\cdot\nabla\zeta =  0,
\qquad 
u = \nabla\times (\mathcal{G}\zeta\nabla\theta).
\end{equation}
The kinetic energy, impulse\footnote[3]{The impulse functional $P$ is twice the axial component of the physical hydrodynamic impulse.}, and $L^p$ Casimir invariants are expressed as functionals of $\zeta$:
\begin{align}
E[\zeta]
&=\pi \iint_{\mathbb R^2_+\times\mathbb R^2_+}
G(z,r,z',r')\zeta(z,r)\zeta(z',r')rr'
\,dzdr\,dz'dr',
\label{eq:kineticenergy}\\
P[\zeta]
&=2\pi \int_{\mathbb R^2_+}
r^3\zeta(z,r)\,dzdr,
\label{eq:impulse}\\
Q_p^p[\zeta]
&=2\pi \int_{\mathbb R^2_+}
r|\zeta(z,r)|^p\,dzdr,
\quad 1\le p<\infty.
\label{eq:casimir}
\end{align}
For $p=\infty$, we set
$Q_\infty[\zeta]=\|\zeta\|_{L^\infty(\mathbb R^2_+)}$.
The circulation is given by
\begin{align}
\Gamma[\zeta]
=\int_{\mathbb R^2_+}r\zeta(z,r)\,dzdr
=\frac{Q_1[\zeta]}{2\pi}
\quad (\zeta\ge0).  \label{eq:circulation}
\end{align}
To define the impulse, we introduce the weighted $L^1$ space
$L^1_*(\mathbb R^3)$, equipped with the norm
\[
\|\zeta\|_{L^1_*(\mathbb R^3)}
=\|r^2\zeta\|_{L^1(\mathbb R^3)}.
\]

We denote by $L^p_r=L^p(\mathbb R^2_+;r)$
the weighted Lebesgue space with respect to the measure
$2\pi r\,dzdr$, equipped with the norm
\[
\|\zeta\|_{L^p_r}
=
\left(
2\pi\int_{\mathbb R^2_+}
|\zeta|^p r\,dzdr
\right)^{1/p},
\quad 1\le p<\infty.
\]
We also define the weighted space
$L^1_{*,r}=L^1_*(\mathbb R^2_+;r)$,
equipped with the norm
\[
\|\zeta\|_{L^1_{*,r}}=P[|\zeta|].
\]
For axisymmetric functions, the $L^p_r$- and $L^1_{*,r}$-norms coincide with the $L^p$- and $L^1_{*}$-norms on $\mathbb{R}^3$, respectively.

For $1\le p\le\infty$, we consider the
scale-invariant variational inequality
\begin{equation}\label{eq:variational-inequality}
E[\zeta]
\le I_{p} P[\zeta]^{\alpha_p}
Q_p[\zeta]^{\beta_p},\quad \zeta\ge 0,
\end{equation}
where the exponents are given by \eqref{eq:exponent}.
To determine the optimal constant $I_{p}$,
we introduce the variational problem
\begin{align}
I_{p,\mu,\nu}
&= \sup_{\zeta\in K_{p,\mu,\nu}} E[\zeta],
\label{eq:variational}\\
K_{p,\mu,\nu}
&= \left\{
\zeta\in (L^1_*\cap L^p)(\mathbb{R}^2_+;r):
\zeta\ge 0,\,
P[\zeta]\le\mu,\,
Q_p[\zeta]\le\nu
\right\},
\label{eq:variational2}
\end{align}
for $\mu,\nu>0$. By the scaling properties of the energy and constraints,
we have
\begin{equation}\label{eq:scaling}
I_{p,\mu,\nu}
=
\mu^{\alpha_p}\nu^{\beta_p}I_{p,1,1}.
\end{equation}
Thus, the optimal constant in
\eqref{eq:variational-inequality} is
$I_p=I_{p,1,1}$, possibly infinite. We can now reformulate Question \ref{q:q1} more precisely as follows.

\begin{qu}\label{q:q2}
For which $1\le p\le\infty$ is the optimal constant $I_p$ finite and attained?
Do the maximizers generate steadily translating vortex rings, and how do their geometry, vortex-core profiles, and propagation speeds depend on $p$?
\end{qu}

\subsection{Main results}

We now present our main results, which establish a connection between the energy--Casimir variational principle \eqref{eq:variational} and the Kelvin--Hicks formula \eqref{eq:Kelvin-Hicks}.
In particular, the limit $p\to1^+$ yields a Kelvin--Hicks formula with a constant correction corresponding to a nonuniform vortex-core profile, without prescribing the circulation, ring radius, or vortex-core profile.

We consider three regimes: $1<p<\infty$, $p=\infty$, and $p\to1^+$. By scaling, all the results below extend to the variational problem \eqref{eq:variational} for arbitrary $\mu,\nu>0$. For simplicity of presentation, we state them under the normalized constraints $\mu=\nu=1$.

\subsubsection{Finite Casimir exponents: $1<p<\infty$}

We first show that, for every $1<p<\infty$, the optimal constant $I_p$ is finite and energy maximizers exist.
These maximizers generate traveling-wave solutions of the Euler equations describing steadily translating vortex rings.
Moreover, the set of energy maximizers is orbitally stable under the axisymmetric Euler equations without swirl.

\begin{theorem}[Maximizers]\label{t:MT1}
For every $1<p<\infty$, the optimal constant $I_p$ is finite and positive, and the supremum defining $I_p$ is attained.
Let $\zeta$ be any maximizer.
Then
\[
P[\zeta]=Q_p[\zeta]=1.
\]
Moreover, $\zeta$ has compact support and, up to a translation in the $z$-direction, is Steiner symmetric with respect to $z=0$.
It also satisfies
\begin{equation}\label{eq:EL_intro0}
\zeta =\kappa_p^{-\frac1{p-1}}
\left( \mathcal{G}\zeta-\frac{W_p}{2}r^2 \right)_+^{\frac1{p-1}},
\end{equation}
where $s_+ = \max\{s, 0\}$ and
\begin{equation}\label{eq:speed}
\frac{W_p}{2}=\alpha_p I_p>0,
\qquad
\kappa_p=\beta_p I_p>0.
\end{equation}
Furthermore, $\zeta$ generates a traveling-wave solution to the Euler equations \eqref{eq:Euler1} of the form
\begin{equation}\label{eq:TW}
\begin{aligned}
u_{\mathrm{tw}}(x,t)&=u(x+u_\infty t)-u_\infty,\\
\zeta_{\mathrm{tw}}(x,t)&=\zeta(x+u_\infty t),
\end{aligned}
\end{equation}
where $u_\infty = -W_p\nabla z$.
\end{theorem}

The conservation of energy, impulse, and Casimir invariants, together with the compactness of maximizing sequences, yields orbital stability of the set of energy maximizers.
We establish this stability for the unique global-in-time solutions of the Euler equations constructed by Danchin \cite[Theorem 1]{Danchin} for axisymmetric initial data without swirl satisfying the following critical $L^{3,1}$ conditions
\begin{align}
\zeta_0\in L^{3,1}(\mathbb{R}^3),\quad
r\zeta_0\in L^{\infty}\cap L^{3,1}(\mathbb{R}^{3}).
\label{eq:Danchin}
\end{align}

\clearpage

\clearpage

\begin{theorem}[Orbital stability]\label{t:MT2}
Let $1<p<\infty$, and let $S_p$ denote the set of maximizers.
Then $S_p$ is orbitally stable under the axisymmetric Euler equations without swirl \eqref{eq:Euler1}.

More precisely, for every $\varepsilon>0$, there exists $\delta>0$ such that, for any axisymmetric and nonnegative initial datum $\zeta_0\in L^{p}\cap L^{1}_{*} (\mathbb{R}^{3})$ satisfying \eqref{eq:Danchin} and
\begin{equation}\label{eq:stability}
\inf_{\eta\in S_p}
\|\zeta_0-\eta\|_{L^p(\mathbb{R}^3)}
+\left|P[\zeta_0]-1\right|
\le\delta,
\end{equation}
the corresponding unique global-in-time solution $\zeta(t)$ satisfies
\begin{equation}\label{eq:stability2}
\inf_{\eta\in S_p}
\|\zeta(t)-\eta\|_{ L^p\cap L^{1}_{*} (\mathbb{R}^3)}
\le\varepsilon,
\quad t\in\mathbb{R}.
\end{equation}
\end{theorem}

\subsubsection{Infinite Casimir exponent: $p=\infty$}

We now investigate how the Casimir exponent $p$ determines the shape and topology of energy-maximizing vortex rings.
The endpoint $p=\infty$ provides a natural starting point, since Hill's spherical vortex is characterized by the energy maximization problem \eqref{eq:variational}.

\begin{figure}[htbp]
\centering
\includegraphics[width=1.0\textwidth]{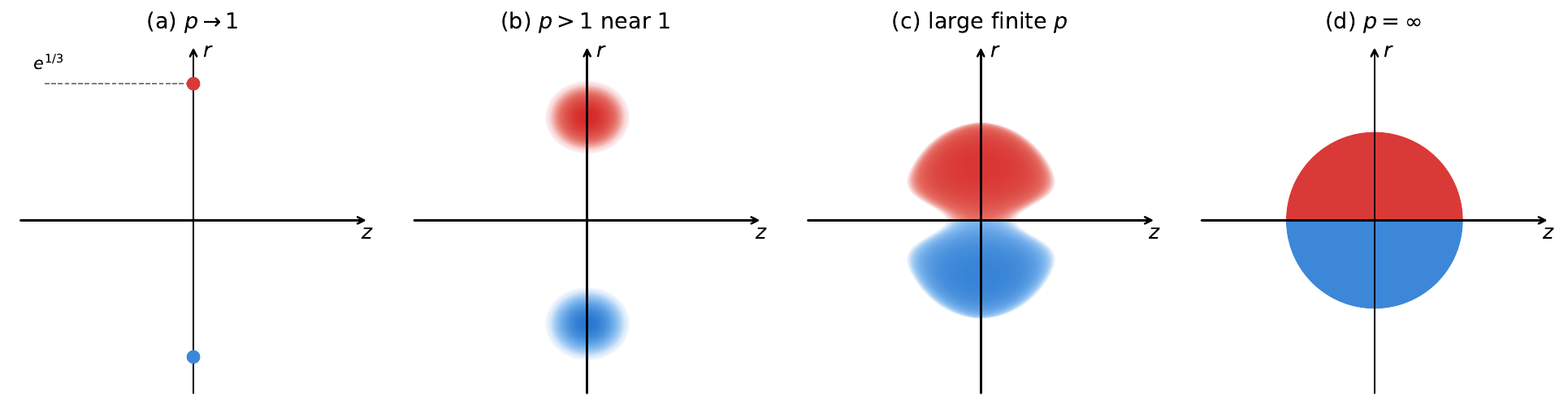}
\caption{
The relationship between the topology of vortex rings and the Casimir exponent $p$, as established in Theorems \ref{t:MT3} and \ref{t:MT4}. The vortex cores touch the symmetry axis already for sufficiently large finite $p$.
}
\label{fig:vortex_ring}
\end{figure}

Hill's vortex is the only known explicit example of a steadily translating vortex ring for the three-dimensional axisymmetric Euler equations without swirl.
Its relative vorticity and traveling speed are given by 
\begin{align}
\zeta_H=\mathbf{1}_{B_{\varrho_H}(0)_+},\qquad
B_{\varrho_H}(0)_+&=\{(z,r)\in\mathbb R^2_+:
z^2+r^2<\varrho_H^2\},\qquad
\varrho_H=\left(\frac{15}{8\pi}\right)^{1/5},
\label{eq:normalizedHill0}\\
W_{\infty}&=\frac{2}{15}\varrho_H^2\approx 0.1084, \label{eq:speedHill} 
\end{align}
where the normalization is chosen so that $P[\zeta_H]=Q_\infty[\zeta_H]=1$; see Figure \ref{fig:vortex_ring}(d).
At the endpoint $p=\infty$, Burton \cite[Theorem 2(i)]{Bur03} proved that Hill's vortex is the unique energy maximizer of \eqref{eq:variational} up to axial translations.
Consequently, the optimal constant is
\begin{align}  \label{eq:Iinfty}
I_\infty= \frac{1}{21}\varrho_H^{2}
\approx 0.03874.
\end{align}

For sufficiently large $p$, the vortex rings generated by energy maximizers have vortex cores homeomorphic to a ball, and their limiting profile as $p\to\infty$ is Hill's vortex; see Figure \ref{fig:vortex_ring}(c).
At the endpoint $p=\infty$, the uniqueness of the energy maximizer yields the orbital stability of Hill's vortex.

\begin{theorem}[Spherical vortex rings]\label{t:MT3}
(i) (Geometric characterization). There exists $1<p_L<\infty$ such that, for every $p_L\le p\le \infty$ and every maximizer $\zeta\in S_p$, there exists $\varrho_p\in C^1[0,\pi]$ such that, after a translation in the $z$-direction,
\begin{align}
\operatorname{supp}\zeta
=
\left\{
(\varrho\cos\vartheta,\varrho\sin\vartheta)\in \overline{\mathbb{R}^{2}_{+}}:
0\le \varrho\le \varrho_p(\vartheta),\
0\le\vartheta\le\pi
\right\}.
\label{eq:suppgraph}
\end{align}
In particular, the support of $\zeta$ in $\mathbb R^3$ is homeomorphic to a ball.

\noindent
(ii) (Convergence to Hill's vortex). 
We have
\begin{align}
\lim_{p\to \infty}I_p=I_\infty,\quad 
\lim_{p\to\infty}W_p=W_{\infty}. \label{eq:Iinftyconv}
\end{align}
As $p\to\infty$, after suitable translations
in the $z$-direction, the maximizers $\zeta$ converge to $\zeta_H$ locally uniformly in $\mathbb{R}^{3}\setminus (\partial B_{\varrho_H}(0)\cup (B_{\varrho_H}(0)\cap \{r=0\}))$. Furthermore, we have
\begin{align}
\lim_{p\to\infty}
\|\varrho_p-\varrho_H\|_{C^1[0,\pi]}
=0.
\label{eq:lar_p_C1}
\end{align}
(iii) (Uniqueness). For $p=\infty$, the maximizers of the variational problem
\eqref{eq:variational} coincide with Hill's vortex up to translations in the $z$-direction. Namely, 
\begin{align}
S_\infty
=
\left\{
\zeta_H(z+c,r):c\in\mathbb{R}
\right\}.
\end{align}

\noindent
(iv) (Orbital stability). For every $q\in(1,\infty)$ and every $\varepsilon>0$, there exists $\delta>0$ such that, for any axisymmetric and nonnegative initial datum $\zeta_0\in L^{q}\cap L^{1}_{*}(\mathbb{R}^{3})$ satisfying \eqref{eq:Danchin} and
\begin{equation}\label{eq:stab_in_hill}
\|\zeta_0-\zeta_H\|_{L^q(\mathbb{R}^3)}
+
\left|P[\zeta_0]-1\right|
\le\delta,
\end{equation}
the corresponding unique global solution $\zeta(t)$ of \eqref{eq:Euler1} satisfies
\begin{equation}\label{eq:stab_hill}
\inf_{c\in\mathbb R}
\|\zeta(t)-\zeta_H(\cdot+c,\cdot)\|_{L^q\cap L^{1}_{*}(\mathbb{R}^3)}
\le\varepsilon,\quad t\in\mathbb R.
\end{equation}
\end{theorem}

\begin{remark}[Comparison with Choi's stability result] 
The orbital stability of Hill's vortex was established by Choi \cite{Choi2024}.
Theorem \ref{t:MT3}(iv) strengthens this result by removing the circulation constraint on perturbations and establishing orbital stability in $L^q(\mathbb{R}^{3})$ for every $1<q<\infty$.
The stability statement remains valid if $L^q(\mathbb{R}^{3})$ and $L^q\cap L^{1}_{*}(\mathbb{R}^{3})$ are replaced by $L^q\cap L^{1}(\mathbb{R}^{3})$ and $L^q\cap L^{1}\cap L^{1}_{*}(\mathbb{R}^{3})$, respectively.
A more detailed comparison is given in Remark \ref{r:stabilitycomparison}.
\end{remark}

\subsubsection{Small Casimir exponents: $p\to1^+$}

In contrast to the regime $p\to\infty$, the opposite limit $p\to1^+$ leads to toroidal vortex rings concentrating on a circular vortex filament.
At the endpoint $p=1$, the scale-invariant energy functional \eqref{eq:nondim2} is proportional to the dimensionless energy defined in terms of circulation
in \eqref{eq:nondim}, and its optimal constant is infinite.
Our main result in this regime determines the sharp asymptotic behavior of $I_p$, the limiting filament radius, and the circulation as $p\to1^+$. We denote by $\mathcal{P}(\mathbb{R}^{2}_{+})$ the space of probability measures on $\mathbb{R}^{2}_{+}$. 

\begin{theorem}[Toroidal vortex rings]\label{t:MT4}
(i) (Geometric characterization). There exists $1<p_{S}<2$ such that, for every $1 <  p \le p_{S}$, every maximizer $\zeta\in S_p$ forms a toroidal vortex ring whose vorticity support is separated from the symmetry axis. 

\noindent
(ii) (Sharp energy asymptotics). We have
\begin{equation}
\label{eq:p1_exact_en}
\lim_{p\to1^+}(p-1)I_p
=\frac{1}{24\pi^2e},\quad \lim_{p\to1^+}(p-1)W_p
=\frac{1}{24\pi^2e}.
\end{equation}

\noindent
(iii) (Convergence to a vortex filament).
As $p\to1^+$, the maximizers concentrate on a circular vortex filament of radius $R=e^{1/3}$.
Namely, for any sequence $p_k\to1^+$ and any sequence of maximizers $\zeta_k\in S_{p_k}$, there exist translations in the $z$-direction such that
\begin{equation}\label{eq:p1_mu_concen}
\begin{aligned}
2\pi r^3\zeta_k(z,r)\,dz\,dr
&\rightharpoonup \delta_{(0,e^{1/3})},\\
2\pi r\zeta_k(z,r)^{p_k}\,dz\,dr
&\rightharpoonup \delta_{(0,e^{1/3})}
\end{aligned}
\qquad
\text{in }\mathcal{P}(\mathbb{R}^2_+).
\end{equation}

\noindent
(iv) (Convergence of circulation). For any sequence $p_k\to1^+$ and any sequence of maximizers $\zeta_k\in S_{p_k}$, we have 
\begin{align}
\lim_{k\to\infty }  \Gamma[\zeta_k] = \frac{1}{2\pi e^{2/3}}.  \label{eq:L1masslimit}
\end{align}
\end{theorem}

In particular, the limiting filament radius is not prescribed in advance but is uniquely selected by the energy maximization principle. This selection follows from the sharp leading-order energy asymptotics as $p\to1^+$. 

\subsubsection{Kelvin--Hicks formula}

Finally, we derive the Kelvin--Hicks formula by analyzing the asymptotic behavior of the vortex-core profile $\omega=r\zeta$ and the traveling speed $W_p$ as $p\to1^+$.
In what follows, we extend the core profile $\omega$ by zero to $\mathbb R^2$ and denote the coordinates on $\mathbb R^2$ by $X=(z,r)$.

A key step in this analysis is to identify the limiting profile of the concentrating vortex core.
The logarithmic singularity of the axisymmetric Green kernel leads to a limiting variational problem involving logarithmic interaction energy and entropy.
The optimizers of this problem are precisely those of the logarithmic Hardy--Littlewood--Sobolev (log-HLS) inequality.

For a probability density $\omega$ on $\mathbb R^2$, we define the internal energy, entropy, and free energy by
\begin{align}
U[\omega]
&=2\iint_{\mathbb R^2\times\mathbb R^2}
\log|X-Y|\,\omega(X)\omega(Y)\,dX\,dY
+1+\log\pi, \label{eq:U}\\
S[\omega]
&=-\int_{\mathbb R^2}
\omega(X)\log\omega(X)\,dX, \label{eq:S}\\
H[\omega]
&=U[\omega]-S[\omega]. \label{eq:logHLSfunc}
\end{align}
The logarithmic Hardy--Littlewood--Sobolev (log-HLS) inequality \cite{CL1992,Beckner93} states that
\begin{align}
H[\omega]\ge 0, \label{eq:HlogHLS}
\end{align}
with equality attained by elements of 
\begin{align}
\mathcal S_*
=
\left\{\frac1{s^2}
\omega_*\left(\frac{X-X_0}{s}\right)
:
s>0,\ X_0\in\mathbb R^2
\right\},\quad \omega_*(X)
=\frac{1}{\pi(1+|X|^2)^2}.
\label{eq:optimizer}
\end{align}
The log-HLS inequality also plays a fundamental role in the critical Keller--Segel model \cite{CF2013,Carlen2025} and in the variational theory of planar vortices
\cite{GS24}.
For our purposes, we use the following quantitative stability estimate of Carlen \cite[Theorem 1.3]{Carlen2025}:
\begin{align}
\frac18
\inf_{g\in\mathcal S_*}
\|\omega-g\|_{L^1(\mathbb R^2)}^2
\le H[\omega].
\label{eq:logHLSstability}
\end{align}
See also \cite[Theorem 1.9]{CF2013}.
We further introduce the entropy-based scale functional
\begin{align}
\epsilon[\omega]
=\frac{1}{e\sqrt{\pi}}
\exp\left(\frac12 S[\omega]\right).
\label{eq:scalefactor}
\end{align}
We show that, as $p\to1^+$, the rescaled vortex-core profiles converge to a log-HLS optimizer and the traveling speeds satisfy the Kelvin--Hicks formula with the constant correction $-1$, expressed in terms of the entropy-based core scale.

\begin{theorem}[Core profile and traveling speed]
\label{t:MT5}
Let $\zeta_p\in S_p$ be a Steiner symmetric maximizer, and set
\begin{align}
\omega_p=r\zeta_p,\qquad
\Gamma_p=\Gamma[\zeta_p],\qquad
\hat{\omega}_p=\frac{\omega_p}{\Gamma_p}.  \label{eq:hatomega}
\end{align}
We extend $\omega_p$ by zero to $\mathbb R^2$, so that $\hat{\omega}_p$ is a probability density on $\mathbb R^2$.
Then the following statements hold as $p\to1^+$.

\medskip
\noindent
(i) (Core-scale asymptotics).
The effective radius $R_p=(2\pi\Gamma_p)^{-1/2}$ and the effective core scale $\epsilon_p =\epsilon[\hat{\omega}_p]$ satisfy
\begin{align}
R_p
&=e^{1/3}
\left[ 1+\left(\frac32-\log(8\sqrt2\pi)\right)(p-1) +o(p-1) \right], \label{eq:Rp_asymptotics}\\
\epsilon_p
&=\frac{8}{e^{3}} \exp\left(-\frac{1}{3(p-1)}\right) (1+o(1)).
\label{eq:GRe}
\end{align}

\medskip
\noindent
(ii) (Core-profile convergence).
The free energy satisfies
\begin{equation}
\frac18\inf_{g\in\mathcal S_*}
\|\hat{\omega}_p-g\|_{L^1(\mathbb R^2)}^2
\le H[\hat{\omega}_p]\longrightarrow0.
\label{eq:stabilityapplogHLS}
\end{equation}
In particular, there exist $X_p\in \mathbb{R}^{2}$ and $s_p>0$ such that 
\begin{equation}
s_p^2
\hat{\omega}_p(X_p+s_p X)
\longrightarrow\omega_*(X)
\quad\text{strongly in }L^1(\mathbb R^2).
\label{eq:p1_core_prof}
\end{equation}

\medskip
\noindent
(iii) (Kelvin--Hicks formula).
The traveling speed $W_p$ satisfies
\begin{equation}
W_p
=\frac{\Gamma_p}{4\pi R_p}
\left( \log\frac{8R_p}{\epsilon_p} -1+o(1) \right).
\label{eq:W_Kelvin0}
\end{equation}
\end{theorem}

\begin{remark}[$2$D mutual induction vs. $3$D self-induction]
The logarithmic divergence of the traveling speed in the classical Kelvin--Hicks formula \eqref{eq:Kelvin-Hicks} reflects the self-induction of vortex rings in three dimensions.
This contrasts with the motion of point vortices in two dimensions, where an isolated point vortex remains stationary, while a pair of point vortices with opposite circulations $\Gamma$ and $-\Gamma$, located at $(0,R)$ and $(0,-R)$, translates with the finite speed $\Gamma/(4\pi R)$ due to mutual induction.

In three dimensions, however, a circular vortex filament of radius $R$ has a singular self-induced velocity in the axial direction.
Consequently, as the core size $\epsilon\to0$, the traveling speed of a desingularized vortex pair in two dimensions converges to a finite limit, whereas that of a thin-cored vortex ring in three dimensions diverges logarithmically, as described by the Kelvin--Hicks formula \eqref{eq:Kelvin-Hicks}.

In Theorem \ref{t:MT4}, the Casimir exponent $p$ can also be viewed as a desingularization parameter for a circular vortex filament of radius $R=e^{1/3}$.
In joint work with Choi, Jeong, and Sim \cite{ACJSW2026}, we studied energy maximization under Casimir constraints for two-dimensional vortex pairs and obtained the following analogue of \eqref{eq:speed}:
\begin{equation}
w_p
= \frac{4(p-1)}{3p-2}i_{p}.
\label{eq:2D-speed}
\end{equation}
In contrast to the three-dimensional case, this identity shows that the traveling speed $w_p$ of the vortex pair remains bounded as $p\to1^+$, even though the maximal energy $i_{p}$ diverges at the rate $(p-1)^{-1}$. 
\end{remark}

\begin{remark}[Prescribed vs. selected parameters]\label{r:Kelvin-Hicks}
De Valeriola and Van Schaftingen \cite[Theorem 1]{dVVS13} constructed, for given $1<p<2$, $\Gamma>0$, $R>0$, and sufficiently small $\epsilon>0$, a traveling wave solution of \eqref{eq:Euler0} whose vorticity is concentrated in a thin toroidal region of major radius $R$ and satisfies the profile  
\begin{align*}
\zeta=\epsilon^{-2}\left( \mathcal{G}\zeta-\frac{W}{2}r^{2} -\gamma\right)_{+}^{\frac{1}{p-1}},\quad W=\frac{\Gamma}{4\pi R}\log\frac{1}{\epsilon}, \quad \gamma = \frac{3\Gamma R}{8\pi}  \log \frac{1}{\epsilon}.
\end{align*}
In contrast, the formula \eqref{eq:W_Kelvin0} is obtained without prescribing the circulation $\Gamma_p$, the major radius $R_p$, or the core scale $\epsilon_p$.
All three quantities are determined solely by the energy--Casimir variational principle \eqref{eq:variational}.
A key feature of our approach is the definition of an effective core scale through the entropy of the vortex-core profile.
The functional $\epsilon$ in \eqref{eq:scalefactor} is normalized so that $\epsilon[\omega_*]=1$, using $S[\omega_*] = \log\pi + 2$.
More generally, it recovers the dilation parameter of the family of log-HLS optimizers:
\begin{align*}
\epsilon\left[ \frac{1}{s^{2}}\omega_*\left(\frac{\cdot-X_0}{s}\right) \right]=s.
\end{align*}
This choice is motivated by the convergence of the rescaled core profiles to a log-HLS optimizer in \eqref{eq:p1_core_prof}.
With this choice of effective core scale, the constant correction in the Kelvin--Hicks formula is $-1$.
Its relation to the classical value $-1/4$, obtained using a different definition of the core radius, is explained in the following remark.
\end{remark}

\begin{remark}[Correction term in the Kelvin--Hicks formula]\label{r:correction}
For a thin-cored vortex ring with a radially symmetric core profile $\omega$ normalized to have unit mass on the cross-section, the Biot--Savart law yields the following formal asymptotic formula for its traveling speed \cite[(1.3)]{FM2000}:
\begin{align}
W = \frac{\Gamma}{4\pi R}
\left( \log\frac{8R}{\epsilon} -\frac{1}{2} + E_{\textrm{log}}[\omega]+o(1) \right)\quad \textrm{as}\ \frac{\epsilon}{R}\to 0. \label{eq:KelvinHick2}
\end{align}
Here, the logarithmic interaction energy is given by
\begin{align*}
E_{\mathrm{log}}[\omega] =\frac{1}{2}(1+\log\pi-U[\omega]).
\end{align*}
The correction term $-1$ in \eqref{eq:W_Kelvin0} is consistent with this formal expansion.
Indeed, the log-HLS optimizer \eqref{eq:optimizer} yields $U[\omega_*]=S[\omega_*]=\log{\pi}+2$ and 
\[
E_{\textrm{log}}[\omega_{*}]= - \frac{1}{2},
\] 
whereas the classical formula \eqref{eq:Kelvin-Hicks} with the correction $-1/4$ is recovered by the patch profile $\omega_{\textrm{patch}} = \pi^{-1} \mathbf{1}_{B_1}$ with  
\[
E_{\textrm{log}}[\omega_{\textrm{patch}}] = \int_0^1 (1 - r^2)rdr = \frac{1}{4}.
 \]
\end{remark}

\subsection{Comparison with previous work}
\label{subsec:FN}

\begin{figure}[htbp]
\centering
\includegraphics[width=0.60\textwidth]{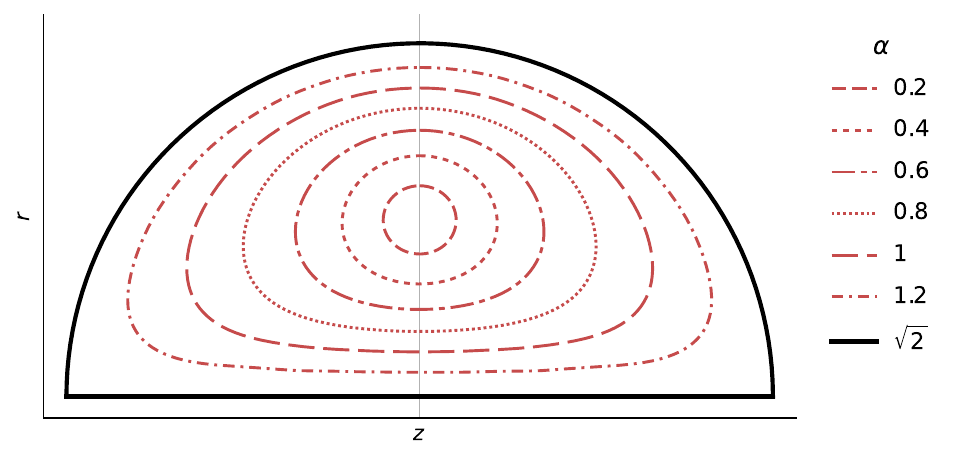}
\caption{Boundaries of the vortex cores of the Fraenkel--Norbury family
in the meridional plane, reconstructed from the numerical
Fourier coefficients in \cite{Nor73}. For every $\alpha<\sqrt{2}$, the vortex core remains
strictly separated from the symmetry axis.
}
\label{fig:FNF}
\end{figure}

The energy-maximizing vortex rings constructed in this paper provide a complementary perspective on the classical Fraenkel--Norbury family of steady vortex rings \cite{Fra70,Fra72,Nor72,Nor73}.
The Fraenkel--Norbury family connects thin-core vortex rings to Hill's spherical vortex and is parametrized by the non-dimensional mean core radius
\begin{align}\label{eq:FN_alpha}
\alpha=\frac{\sqrt{|\Omega|}}{\sqrt{\pi}R_{\mathrm{av}}},
\qquad
0<\alpha\le\sqrt{2},
\end{align}
where $\Omega$ is the meridional cross-section of the vortex core and $R_{\mathrm{av}}$ is the arithmetic mean of its inner and outer radii at $z=0$. Unlike the Casimir exponent $p$, the parameter $\alpha$ directly characterizes the geometry of the vortex core.
For $0<\alpha<\sqrt{2}$, the vortex core is toroidal and remains separated from the symmetry axis, whereas $\alpha=\sqrt{2}$ corresponds to Hill's spherical vortex; see Figure \ref{fig:FNF}. The relative vorticity of the Fraenkel--Norbury family has the patch-type profile
\begin{align}
\zeta=\lambda\mathbf{1}_{\Omega},
\qquad
\Omega=\left\{
\mathcal{G}\zeta-\frac{W}{2}r^2-\gamma>0
\right\},
\end{align}
where $\lambda>0$ is the vorticity strength, $W$ is the propagation speed, and $\gamma$ is the flux constant.
For each prescribed $\alpha$, Norbury \cite{Nor73} numerically determined the vortex-core shape, propagation speed, and flux constant, together with the kinetic energy, impulse, and circulation.
Although global existence results for steady vortex rings have been established \cite{FB74,FT81,AT88,AS89}, the existence of a vortex ring in the Fraenkel--Norbury family for every prescribed $\alpha\in(0,\sqrt{2})$ remains open.

While the Fraenkel--Norbury family describes vortex-ring geometry through the geometric parameter $\alpha$, our variational family reveals how the Casimir exponent $p$ selects the geometry of energy-maximizing vortex rings.
Furthermore, as $p\to1^+$, the leading-order energy asymptotics determine the limiting filament radius, while the next-order asymptotics characterize the microscopic vortex-core structure through the log-HLS inequality; see Table \ref{tab:comparison}.

A related variational problem for two-dimensional vortex dipoles was studied in \cite{ACJ2025,ACJSW2026}.
Unlike the present three-dimensional vortex rings, whose cores are toroidal for sufficiently small $p>1$, the energy-maximizing two-dimensional vortex dipoles touch the symmetry axis for all $1 <  p  \le \infty$.

The variational principle studied in this paper can be viewed as a coarse-grained version of the Kelvin--Benjamin principle \cite{Kelvin1875,Ben76} for vortex rings or Arnold's variational principle \cite{Arnold1965} for steady Euler flows.
In the Kelvin--Benjamin and Arnold principles, variations are taken on isovortical surfaces or coadjoint orbits, on which all Casimir invariants are preserved; see \cite{Burton87a,Bur03} for rearrangement formulations and \cite{FT81} for a coarse-grained formulation.

Arnold's energy--Casimir method provides a classical variational framework for studying the structure and stability of steady Euler flows through Casimir invariants.
This method has found numerous applications to two-dimensional Euler flows \cite{AK21}; see also \cite{CDG21b,RS2026} for axisymmetric Euler flows without swirl.
In the latter setting, steady solutions can be constructed from solutions of the semilinear elliptic equation
\begin{equation*}
-\frac{1}{r^2}L\Psi =  f(\Psi).
\end{equation*}
For two-dimensional steady Euler flows satisfying analogous vorticity--stream function relations, various structural properties have been investigated, including the local structure of steady solutions \cite{CS12,CEW23,DS24} and rigidity and flexibility under domain deformations \cite{CDG21b}.
For surveys on Arnold's energy--Casimir method and recent developments in the PDE theory of steady Euler flows, see \cite{KMS23,Elg26}.

\begin{landscape}
\thispagestyle{plain}

\vspace*{\fill}

\begin{center}
\renewcommand{\arraystretch}{1.8}
\setlength{\tabcolsep}{10pt}
\small

\begin{tabularx}{0.94\linewidth}{
|>{\raggedright\arraybackslash}m{2.5cm}|
>{\raggedright\arraybackslash}m{3.0cm}|
>{\centering\arraybackslash}X|
>{\centering\arraybackslash}X|}
\hline

Category
& Property
& Fraenkel--Norbury family
& Present variational family
\\ \hline

\multirow{4}{2.5cm}[-0.5em]{\textbf{General framework}}
& Parameter
& Mean core radius
  $\alpha\in(0,\sqrt{2}]$
& Casimir exponent
  $p\in(1,\infty]$
\\ \cline{2-4}

& Construction
& Analytical and numerical
& Energy--Casimir maximization
\\ \cline{2-4}

& Vorticity profile
& $\displaystyle
\zeta=\lambda\mathbf{1}_{\Omega}$
& $\displaystyle
\zeta=
\kappa_p^{-\frac{1}{p-1}}
\left(
\mathcal{G}\zeta-\frac{W_p}{2}r^2
\right)_+^{\frac{1}{p-1}}$
\\ \cline{2-4}

& Orbital stability
& Not established for the entire family
& Established for all energy maximizers
\\ \hline

\multirow{2}{2.5cm}[-2.3em]{\textbf{Macroscopic structure}}
& Limiting geometry
& \makecell[c]{
  Thin-core limit: $\alpha\to0$\\[0.4em]
  Hill's vortex: $\alpha=\sqrt{2}$}
& \makecell[c]{
  Thin-core limit: $p\to1^+$\hspace{3pt}
  ($R_p\to e^{1/3}$)\\[0.4em]
  Hill's vortex: $p=\infty$}
\\ \cline{2-4}

& Characterization
& \makecell[c]{
  Vortex geometry\\[0.6em]
  \hspace{-5pt}$\longrightarrow$ Numerical solutions\\[0.6em]
  $\longrightarrow$ Conserved quantities}
& \makecell[c]{
  Casimir exponent\\[0.6em]
  \hspace{-17pt} $\longrightarrow$ Energy maximizers\\[0.6em]
  \hspace{-28pt} $\longrightarrow$ Vortex geometry}
\\ \hline

\multirow{3}{2.5cm}[-0.5em]{\textbf{Microscopic structure}}
& Limiting core profile
& Uniform vortex patch
& Log-HLS optimizer as $p\to1^+$
\\ \cline{2-4}

& Core radius
& $\epsilon=\displaystyle\sqrt{\frac{|\Omega|}{\pi}}$
& Entropy-based effective core radius $\epsilon_p$
\\ \cline{2-4}

& Thin-core speed
& $\displaystyle
W=\frac{\Gamma}{4\pi R}
\left(
\log\frac{8R}{\epsilon}
-\frac14+o(1)
\right)$
& $\displaystyle
W_p=\frac{\Gamma_p}{4\pi R_p}
\left(
\log\frac{8R_p}{\epsilon_p}
-1+o(1)
\right)$
\\ \hline

\end{tabularx}

\vspace{10pt}

\captionsetup{hypcap=false}
\captionof{table}{Comparison of the Fraenkel--Norbury family and the present variational family of steady vortex rings. The energy--Casimir principle determines the macroscopic structure of vortex rings, with the leading-order asymptotics of $I_p$ as $p\to1^+$ selecting the limiting filament radius.
The next-order asymptotics, governed by the log-HLS inequality, reveal the microscopic vortex-core structure.}
\label{tab:comparison}
\end{center}

\vspace*{\fill}
\end{landscape}

\subsection{Ideas of the proof}

We outline the main ideas of the proofs in three regimes: $1<p<\infty$, $p=\infty$, and $p\to1^+$, including the derivation of the Kelvin--Hicks formula in the last regime.

\subsubsection{$1<p<\infty$}

The starting point of our analysis is the scale-invariant energy inequality \eqref{eq:variational-inequality}.
For arbitrary constants $a,b>0$, consider the rescaling
\begin{equation}
\zeta_{a,b}(z,r)=a\zeta(bz,br).  \label{eq:scalinglaw}
\end{equation}
Then
\begin{align}
E[\zeta_{a,b}]=\frac{a^2}{b^7}E[\zeta],\quad 
P[\zeta_{a,b}]=\frac{a}{b^5}P[\zeta],\quad
Q_p[\zeta_{a,b}]=\frac{a}{b^{3/p}}Q_p[\zeta].  \label{eq:scalingconservation}
\end{align}
Consequently, the exponents $\alpha_p$ and $\beta_p$ for which the identity
\[
\mathscr{E}_p[\zeta_{a,b}]
=\mathscr{E}_p[\zeta]
\]
holds are uniquely determined by the scaling relations and are precisely those given in \eqref{eq:exponent}.
At the endpoint $p=\infty$, the variational inequality \eqref{eq:variational-inequality} is equivalent to the following Sobolev-type interpolation inequality on $\mathbb R^5$:
\[
\|\nabla\varphi\|_{L^2(\mathbb R^5)}
\le C_s
\|\Delta\varphi\|_{L^1(\mathbb R^5)}^{7/10}
\|\Delta\varphi\|_{L^\infty(\mathbb R^5)}^{3/10}.
\]
This equivalence follows from the fact that the transformation
\[
\psi\longmapsto\varphi=\frac{\psi}{r^2}
\]
is an isometry from the weighted homogeneous Hilbert space $\dot{H}^1_0(\mathbb R^2_+;r^{-1})$ onto the axisymmetric homogeneous Sobolev space $\dot H^1_{\mathrm{axi}}(\mathbb R^5)$.

We begin by establishing energy estimates \eqref{eq:variational-inequality} based on the axisymmetric Green function.
We then formulate the variational problem \eqref{eq:variational}, derive the Euler--Lagrange equation \eqref{eq:EL_intro0} and the Poho\v{z}aev identities \eqref{eq:speed}, and establish structural properties of energy maximizers.

A distinctive feature of our variational formulation is that the maximal energy $I_{p,\mu,\nu}$ depends on two constraint parameters $\mu$ and $\nu$.
We establish a strict superadditivity property with respect to the distribution of these parameters, which allows us to apply Lions' concentration--compactness principle to prove the existence of energy maximizers and the orbital stability of the set of maximizers.

\subsubsection{$p=\infty$}

At the endpoint $p=\infty$, energy maximizers exist and are unique up to axial translations, coinciding with Hill's spherical vortex.
The uniqueness of Hill's vortex can be established by the moving plane method applied to the associated elliptic equation \cite[Theorem 1.1]{AF86}.
Burton \cite{Bur03} gave a variational proof without relying on the elliptic equation.
Here, we provide an elementary proof, also based solely on the variational inequality.

Orbital stability at $p=\infty$ requires a different approach from the case of finite $p$, since stability in $L^\infty$ cannot be expected for vortex patches.
Variational compactness yields only weak-* convergence in $L^\infty$ and strong convergence in $L^1_{*,r}$.
By incorporating the $L^q$ norm into the stability distance, we first establish orbital stability in $L^q\cap L^1_{*,r}$ for every $1<q<\infty$, under the constraint $Q_\infty[\zeta_0]\le 1$.
This follows from the conservation of the $L^q$ Casimir and the uniform convexity of $L^q$.
We then remove the $L^\infty$ constraint on the initial data by exploiting the conservation of Casimir invariants.
To remove this constraint, we use the decomposition
\[
\zeta(t)=\min\{\zeta(t),1\}+(\zeta(t)-1)_+.
\]
The conservation of Casimir invariants for both truncated parts allows us to control the excess above $1$ and apply variational compactness to the bounded part.

Finally, to prove Theorem~\ref{t:MT3}(i)--(ii), we combine variational compactness with the uniqueness of Hill's vortex to obtain convergence of maximizers as $p\to\infty$.
Estimates for the five-dimensional Newtonian potential and the non-degeneracy of Hill's free boundary then yield $C^1$ convergence of the free boundaries, showing that the vortex cores touch the symmetry axis and are homeomorphic to a ball for sufficiently large $p$.

\subsubsection{$p\to1^+$}\label{s:1.6.3}

The main mathematical challenge of this paper is to establish the sharp energy asymptotics \eqref{eq:p1_exact_en}.
Although these asymptotics concern the energy, they also determine the limiting geometry of energy maximizers, yielding convergence to a single circular vortex filament of radius $R=e^{1/3}$ in \eqref{eq:p1_mu_concen}.

The Green function estimates give $I_p=O((p-1)^{-1})$.
Testing the variational problem with vortex patches concentrated near $(0,R)$ yields
\begin{align}
\liminf_{p\to1^+}(p-1)I_p \ge \frac{\log R}{8\pi^2R^3}.  \label{eq:SLBR0}
\end{align}
Optimizing over $R\ge1$ gives the sharp lower bound $(24\pi^2e)^{-1}$, attained at $R=e^{1/3}$.

The main difficulty is to establish the matching upper bound.
We choose a sequence $p_k\to1^+$ such that
\[
(p_k-1)I_{p_k}\to \mathcal{I}_*,
\]
and let $\zeta_k$ be an energy maximizer for $I_{p_k}$.
We introduce the probability measures
\begin{align}
d\mathfrak{p}_k=2\pi r^3\zeta_k\,dz\,dr,
\qquad
d\mathfrak{q}_k=2\pi r\zeta_k^{p_k}\,dz\,dr. \label{eq:meas_cons0}
\end{align}
By concentration--compactness, these measures converge, up to subsequences and axial translations, to probability measures $\mathfrak{p}$ and $\mathfrak{q}$.
We then establish an atomic energy upper bound for the limiting measures.
For each atom $X_j=(z_j,r_j)$, let $\mathsf{p}_j=\mathfrak{p}(\{X_j\})$ and $\mathsf{q}_j=\mathfrak{q}(\{X_j\})$.
We show that
\begin{align}
\frac{1}{24\pi^2e}
\le \mathcal{I}_*
\le \sum_j e_{r_j}(\mathsf{p}_j,\mathsf{q}_j), \label{eq:Lambdabound0}
\end{align}
where
\begin{align}
e_r(\mathsf{p},\mathsf{q})
=\frac{\mathsf{p}^2}{16\pi^2r^3}
\log\frac{r^2 \mathsf{q}}{\mathsf{p}} \label{eq:e_Rbound0}
\end{align}
for $\mathsf{p}>0$, and $e_r(0,\mathsf{q})=0$.

To establish the bound \eqref{eq:Lambdabound0}, we decompose the maximizer $\zeta_k$ into the part concentrating near the large atoms and two remainder terms.
We first choose an increasing sequence of compact sets $K_n\Subset\bbR^2_+$ expanding to $\mathbb{R}^{2}_{+}$ and decompose $\bbR^2_+$ into $K_n$ and its complement. We then consider the atoms in $K_n$ whose Casimir weights $\mathsf{q}_j$ exceed $\tau>0$.
We set
\begin{align}
J_{n,\tau}
=
\{j\in \mathbb{N}\ :\ X_j\in K_n,\ \mathsf{q}_j>\tau\}.
\label{eq:Jneta0}
\end{align}
We surround these atoms by balls $B_{\delta_j}=B_{\delta_j}(X_j)\Subset K_n^\circ$ with pairwise disjoint closures and decompose $K_n$ into
\begin{align}
\bigcup_{j\in J_{n,\tau}} B_{\delta_j}(X_j)
\label{eq:cupball0}
\end{align}
and its complement in $K_n$.
We define
\begin{align}
\zeta_k^{\mathrm{atom}}=\sum_{j\in J_{n,\tau}}\zeta_{k,j}^{\textrm{atom}},\quad \zeta_{k,j}^{\textrm{atom}}=\zeta_k\mathbf 1_{B_{\delta_j}},\quad \zeta_k^{\mathrm{rem}}=
\zeta_k\mathbf 1_{K_n\setminus
\bigcup_{j\in J_{n,\tau}}B_{\delta_j}},
\quad
\zeta_k^{\mathrm{far}}=
\zeta_k\mathbf 1_{\bbR^2_+\setminus K_n}.  \label{eq:fgh0}
\end{align}
so that 
\begin{align}
\zeta_k
=
\zeta_k^{\mathrm{atom}}
+\zeta_k^{\mathrm{rem}}
+\zeta_k^{\mathrm{far}}.  \label{eq:fghsum0}
\end{align}
Here, $\zeta_k^{\mathrm{atom}}$ is the part concentrated near the finitely many large atoms, $\zeta_k^{\mathrm{rem}}$ contains the small atoms and the diffuse part inside $K_n$, and $\zeta_k^{\mathrm{far}}$ is the part outside $K_n$.
By estimating the energies of these terms separately and then taking the limits as $\tau\to0$ and $n\to\infty$, we obtain the upper bound \eqref{eq:Lambdabound0} for $\mathcal{I}_*$.
The main contribution to the limiting energy comes from $\zeta_k^{\mathrm{atom}}$.
Using a logarithmic-type HLS inequality derived from the classical HLS inequality and Jensen's inequality, we bound the energy arising from the logarithmic singularity of the Green kernel in terms of the atomic energies.

The sharp asymptotics and limiting geometry follow from optimizing the atomic energies.
Indeed,
\[
e_{r_j}(\mathsf{p}_j,\mathsf{q}_j)
\le\frac{1}{24\pi^2e}
\mathsf{p}_j^{1/2}\mathsf{q}_j^{3/2},
\]
with equality precisely when $r_j^2\mathsf{q}_j/\mathsf{p}_j=e^{2/3}$.
Since $\sum_j\mathsf{p}_j\le1$ and $\sum_j\mathsf{q}_j\le1$, we find
\[
\frac{1}{24\pi^2e}
\le \mathcal{I}_*
\le\frac{1}{24\pi^2e}
\sum_j\mathsf{p}_j^{1/2}\mathsf{q}_j^{3/2}
\le\frac{1}{24\pi^2e}.
\]
Equality forces both measures to concentrate at a single atom $X_{j_*}=(z_{j_*},r_{j_*})$, with
\[
\mathsf{p}_j=\mathsf{q}_j=\delta_{j,j_*}.
\]
Moreover, the equality condition
\[
\frac{r_{j_*}^2 \mathsf{q}_{j_*}}{\mathsf{p}_{j_*}}
=e^{2/3}
\]
implies $r_{j_*}=e^{1/3}$.
Thus, up to axial translations,
\[
\mathfrak{p}=\mathfrak{q}=\delta_{(0,e^{1/3})}.
\]
This proves both the sharp energy asymptotics \eqref{eq:p1_exact_en} and the convergence to a single circular vortex filament \eqref{eq:p1_mu_concen}.
Finally, combining this concentration with the Euler--Lagrange equation and the sharp energy inequality, we exclude residual vorticity near the symmetry axis.
Consequently, every energy maximizer has a vortex core separated from the axis for $p$ sufficiently close to $1$.

\subsubsection{Kelvin--Hicks formula}

We further derive a refined asymptotic expansion for $I_p-\mathcal{I}_*(p-1)^{-1}$ as $p\to1^+$.
We begin by constructing a test function whose core profile is given by the log-HLS optimizer, obtaining a lower bound for $\liminf_{p\to 1^{+}}(I_p-\mathcal{I}_*(p-1)^{-1})$.

Next, using a lower bound for the logarithmic potential, we show that the diameter $d_p$ of the vortex core decays exponentially:
\[
d_p\le C\exp\left(-\frac{c}{p-1}\right).
\]
Combining this estimate with a log-type HLS inequality, we derive a sharp upper bound for the logarithmic interaction energy.
We then use the exponential decay estimate and the asymptotic expansion of $F(s)$ to express the energy associated with the Green kernel in terms of the logarithmic interaction energy, with an error of $o(1)$.
This yields a matching sharp upper bound for $\limsup_{p\to 1^{+}}(I_p-\mathcal{I}_*(p-1)^{-1})$ and establishes the desired asymptotic expansion.

The resulting sharp asymptotics imply that the log-HLS deficit $H[\hat{\omega}_p]$ tends to zero.
By Carlen's stability theorem, it follows that the normalized vortex core profile approaches the family of log-HLS optimizers strongly in $L^1(\mathbb R^2)$.
Finally, applying the Poho\v{z}aev identity to the asymptotic expansion of $I_p$, we derive the asymptotic expansion of the traveling speed $W_p$, which yields the Kelvin--Hicks formula.

\subsection{Organization of the paper}

The rest of the paper is organized as follows.

\begin{itemize}
\item In \S\ref{s:2}, we establish energy estimates using bounds for the Green function and prove strict superadditivity of the maximal energy \eqref{eq:variational} with respect to the two parameters $\mu,\nu>0$.
We also show that Hill's vortex is the unique energy maximizer for $p=\infty$, up to axial translations.
\item In \S \ref{s:3}, we derive the Euler--Lagrange equation satisfied by energy maximizers for $1<p<\infty$, establish identities for the Lagrange multipliers, and investigate properties of energy maximizers, including compact support and Steiner symmetry in the $z$-direction.
\item In \S \ref{s:4}, we apply the concentration--compactness principle to establish the existence of energy maximizers and compactness of maximizing sequences. This completes the proof of Theorem \ref{t:MT1}.
\item In \S \ref{s:5}, we use the compactness of maximizing sequences to establish orbital stability of the set of energy maximizers, completing the proof of Theorem \ref{t:MT2}.
We further establish orbital stability of Hill's spherical vortex using Casimir invariants.
\item In \S \ref{s:6}, we prove that energy maximizers for $1<p<\infty$ converge to Hill's spherical vortex as $p\to\infty$.
We further show that, for sufficiently large $p$, the vorticity supports of energy maximizers touch the symmetry axis and are homeomorphic to a ball.
This completes the proof of Theorem \ref{t:MT3}.
\item In \S \ref{s:7}, we establish sharp asymptotics for the maximal energy $I_p$ as $p\to1^+$ and prove that energy maximizers converge to a circular vortex filament.
This completes the proof of Theorem \ref{t:MT4}.
\item In \S\ref{s:8}, we derive the asymptotic expansion of $I_p$ up to the constant-order term and show that the log-HLS deficit converges to zero. This yields the convergence of the vortex-core profile to the log-HLS optimizer and the Kelvin–Hicks formula, completing the proof of Theorem~\ref{t:MT5}.
\end{itemize}

\section{Variational framework and energy estimates}\label{s:2}

For $1<p\le\infty$, we establish the energy inequality \eqref{eq:variational-inequality} and formulate the associated variational problem \eqref{eq:variational}.
We then prove that the maximal energy is strictly superadditive with respect to the two parameters $\mu,\nu>0$.
In the case $p=\infty$, we show that the energy maximizer is given by Hill's spherical vortex.

\subsection{The five-dimensional isometry}\label{sec:iso}

Throughout this paper, we often identify axisymmetric stream functions $\psi=\mathcal{G}\zeta$ in $\mathbb{R}^{2}_{+}$ with axisymmetric functions $\varphi$ on $\mathbb{R}^5$
via the following mapping
\begin{align}
\psi(z,r)\longmapsto \varphi(y)=\frac{\psi(z,r)}{r^{2}},\quad y=(y_1,y'),\ y_1=z,\ |y'|=r.     \label{eq:5dtrans}
\end{align}
The Dirichlet problem in the two-dimensional half-plane \eqref{eq:elliptichalf} is transformed into the Poisson equation
in the class of axisymmetric functions on $\mathbb{R}^5$:
\begin{align}
-\Delta_{y}\varphi =\zeta,\quad y\in \mathbb{R}^{5}.  \label{eq:5D}
\end{align}
This mapping defines an isometry from
the weighted Hilbert space on the half-plane
to the homogeneous Sobolev space on $\mathbb{R}^5$. Let $\dot{H}^{1}_{0}(\mathbb{R}^{2}_{+};r^{-1})$ denote a Hilbert space endowed with the inner product 
\begin{equation}
\begin{aligned}
(\psi,\tilde{\psi})_{\dot{H}^{1}_0(\mathbb{R}^{2}_{+};r^{-1})}=2\pi^{2}\int_{\mathbb{R}^{2}_{+}}\nabla \psi\cdot \nabla \tilde{\psi}\frac{1}{r}  \, dz  dr.
\end{aligned}
\label{eq:ip}
\end{equation}
This space is isometrically isomorphic to $\dot{H}^{1}_{\textrm{axi}}(\mathbb{R}^{5})= \{ \varphi\in \dot{H}^{1}(\mathbb{R}^{5})\ | \ \varphi\text{ is axisymmetric } \}$ \cite[Lemma 2.2]{AF88} in the sense that 
\begin{equation}
\begin{aligned}
(\psi,\tilde{\psi})_{\dot{H}^{1}_0(\mathbb{R}^{2}_{+};r^{-1})}=\int_{\mathbb{R}^{5}}\nabla_y\varphi \cdot \nabla_y\tilde{\varphi} \, d y,\quad \varphi=\frac{\psi}{r^{2}},\ \tilde{\varphi}=\frac{\tilde{\psi}}{r^{2}}.
\end{aligned}
\label{eq:isometry}
\end{equation}
In the five-dimensional setting, the conserved quantities \eqref{eq:kineticenergy}-\eqref{eq:casimir} are expressed as  
\begin{align}
E[\zeta]
&=
\frac{1}{16\pi^{3}}\iint_{\mathbb{R}^{5}\times \mathbb{R}^{5}}\frac{1}{|y-w|^{3}}\zeta(y)\zeta(w)dydw,
\label{eq:kineticenergy5D}\\
P[\zeta]
&=\frac{1}{\pi} \int_{\mathbb R^5}
\zeta\,dy,
\label{eq:impulse5D}\\
Q_p^p[\zeta]
&=\frac{1}{\pi}\int_{\mathbb R^5}
\frac{1}{|y'|^{2}}|\zeta|^p\,dy,
\quad 1\le p<\infty.
\label{eq:casimir5D}
\end{align}
For $p=\infty$, we have $Q_\infty[\zeta]=\|\zeta\|_{L^\infty(\mathbb R^5)}$. We show that the identity \eqref{eq:kineticenergy5D} holds for $\zeta\in (L^{1}_{*}\cap L^{p})(\mathbb{R}^{2}_{+};r)$ and $1<p\le \infty$. In particular, the kinetic energy is expressed as the Dirichlet energy
\begin{align}
E[\zeta]
=\frac{1}{2\pi} \int_{\mathbb R^5}
|\nabla_y\varphi|^{2}dy
=\pi \int_{\mathbb{R}^{2}_{+}}|\nabla \psi|^{2}\frac{1}{r}dzdr. \label{eq:ke2}
\end{align}
The variational inequality \eqref{eq:variational-inequality} for $p=\infty$ follows from the Sobolev embedding $\dot{W}^{1,10/7}(\mathbb{R}^{5})\subset L^{2}(\mathbb{R}^{5})$. Indeed, by H\"{o}lder's inequality, we have 
\begin{align*}
||\nabla \varphi||_{L^{2}(\mathbb{R}^{5})}\le C_s||\Delta \varphi||_{L^{\frac{10}{7}}(\mathbb{R}^{5})}
\le 
C_s||\Delta \varphi||_{L^{1}(\mathbb{R}^{5})}^{\frac{7}{10}}||\Delta \varphi||_{L^{\infty}(\mathbb{R}^{5})}^{\frac{3}{10}}
=C_s||\zeta||_{L^{1}(\mathbb{R}^{5})}^{\frac{7}{10}}||\zeta||_{L^{\infty}(\mathbb{R}^{5})}^{\frac{3}{10}}.
\end{align*}
Using \eqref{eq:5dtrans} yields 
\begin{align}
E[\zeta]\le \frac{C_s^{2}\pi^{\frac{2}{5}}}{2}P[\zeta]^{\frac{7}{5}}Q_\infty[\zeta]^{\frac{3}{5}}
=\frac{C_s^{2}\pi^{\frac{2}{5}}}{2} \|r^2 \zeta\|_{L^1(\bbR^2_+;r)}^{\frac{7}{5}} ||\zeta||_{L^\infty(\bbR^2_+;r)}^{\frac{3}{5}}. \label{eq:kee0}
\end{align}

\subsection{Green function estimates}

We show the variational inequality \eqref{eq:variational-inequality} for $1<p\le \infty$ using Green function estimates.
In what follows, we write the Green function as $G(X,Y)=G(z,r,z',r')$, using the coordinates $X=(z,r)$ and $Y=(z',r')$ in the half-plane $\mathbb{R}^2_+$.

\begin{proposition}\label{lem:p1_kernel}
The function $F$ in \eqref{eq:Green} is strictly decreasing for $s>0$. As $s \to 0$, we have
\begin{equation}\label{eq:p1_F_log}
F(s) = \log\frac8{\sqrt s} - 2 + O(s|\log s|), \quad
sF'(s) = - \frac{1}{2} + O(s|\log s|).
\end{equation}
Moreover, for every $\tau \in (0,3/2]$ and $s>0$,
\begin{equation}\label{eq:p1_F_pow}
    0 < F(s) \lesssim \frac{1}{\tau  s^{\tau}}.
\end{equation}
\end{proposition}

\begin{proof}
The asymptotic expansion \eqref{eq:p1_F_log} is derived in \cite[Lemma~2.7]{FS2015}, \cite[p. 94]{SverakLec}.
Differentiating $F(s)$ and making a change of variables, we have 
\begin{align*}
-2F'(s) = \int_0^{\pi} \frac{\cos \theta}{(2(1 - \cos \theta) + s)^{\frac{3}{2}}} d\theta
= \int_0^{\pi/2} \cos \theta\left(\frac{1}{(2(1 - \cos \theta) + s)^{\frac{3}{2}}}- \frac{1}{(2(1 +\cos \theta) + s)^{\frac{3}{2}}} \right)d\theta>0.
\end{align*}
Thus $F(s)$ is strictly decreasing for $s>0$. It follows that $F(s) \lesssim 1-\log{s}$ for $0<s\le1$. Since $\sup_{0 < s \le 1} s^\tau\log s^{-1} = (e\tau)^{-1}$, we obtain \eqref{eq:p1_F_pow} for $0<s\le1$.
For $s\ge1$, we have $F(s) \lesssim s^{-3/2} \le s^{-\tau}$.
This proves \eqref{eq:p1_F_pow}. 
\end{proof}

\begin{proposition}
Let $J$ be a compact interval in $(0,\infty)$. Then
\begin{equation}\label{eq:p1_loc_ker}
\pi G(X,Y)X_2Y_2
=
\frac12 (X_2Y_2)^{3/2}
\log\frac{(X_2Y_2)^{1/2}}{|X-Y|}
+ O_J(1)
\quad\text{as } |X-Y|\to 0,
\end{equation}
uniformly for $X_2,Y_2\in J$.

Moreover, there exists a constant $C>0$ such that for every $\tau\in(0,3/2]$,
\begin{equation}\label{eq:Green_bound}
0<G(X,Y)
\le
\frac{C}{\tau}
\frac{(X_2Y_2)^{1/2+\tau}}{|X-Y|^{2\tau}}
\end{equation}
for all $X,Y\in\mathbb{R}^2_+$ with $X\ne Y$.
\end{proposition}

\begin{proof}
The properties \eqref{eq:p1_loc_ker} and \eqref{eq:Green_bound} follow from \eqref{eq:p1_F_log}, \eqref{eq:Green}, and \eqref{eq:p1_F_pow}.
\end{proof}

\begin{proposition}\label{lem:green_strip}
For $X,Y\in\mathbb{R}^2_+$ with $X\ne Y$,
the Green function admits the representations
\begin{align}
G(X,Y)
&=\frac{X_2Y_2}{2\pi}
\int_0^\pi\frac{\cos\theta}{|x-y|}\,d\theta
=\frac{X_2^2Y_2^2}{2\pi}
\int_0^\pi\frac{\sin^2\theta}{|x-y|^3}\,d\theta,
\label{eq:Greenrep}
\end{align}
where $x=(X_2,0,X_1)$ and $y=(Y_2\cos\theta,Y_2\sin\theta,Y_1)$. Moreover, we have 
\begin{align}
|\partial_{X_1}G(X,Y)|+|\nabla_Y G(X,Y)|&\lesssim X_2^{2}\int_{0}^{2\pi}\frac{d\theta}{|x-y|^{2}},\quad X, Y\in\mathbb{R}^{2}_{+}, \label{eq:DanchinGreen}\\
|\partial_{X_2}G(X,Y)|&\lesssim Y_2^{2}\int_{0}^{2\pi}\frac{d\theta}{|x-y|^{2}},\quad X, Y\in\mathbb{R}^{2}_{+}, \label{eq:DanchinGreen2}\\
\int_{\mathbb{R}}G(X,Y)\,dX_1 &= \frac{1}{2} \min\{X_2^{2},Y_2^{2}\}, \quad X_2>0,\quad Y\in\mathbb{R}^{2}_{+},\label{eq:green_z} \\
\int_{\bbR^2_+ \cap \{ X_2 < h \} } G(X,Y)\,X_2 dX &\le \frac{1}{4}Y_2h^{3},\quad h>0,\quad Y\in\mathbb{R}^{2}_{+}. \label{eq:green_strip} 
\end{align}
\end{proposition}

\begin{proof}
We set 
\begin{align*}
D^{2}=|x-y|^{2}
&= 
|X_1-Y_1|^{2}+|X_2-Y_2\cos\theta|^{2} + Y_2^{2}\sin^{2}\theta \\
&=|X_1-Y_1|^{2}+X_2^{2} + Y_2^{2} - 2X_2Y_2\cos\theta.
\end{align*}
It follows that 
\begin{align*}
|X_1-Y_1|,\ X_2\sin\theta,\ Y_2\sin\theta,\ |X_2-Y_2\cos\theta|,\ |Y_2-X_2\cos\theta| \le D,\\
\partial_{X_1}D
=\frac{X_1-Y_1}{D},\quad 
\partial_{Y_1}D
=-\frac{X_1-Y_1}{D},\quad 
\partial_{X_2}D
=\frac{X_2-Y_2\cos\theta}{D},\quad 
\partial_{Y_2}D
=\frac{Y_2-X_2\cos\theta}{D}.
\end{align*}
By integration by parts in $\theta$, we have
\[
\int_{0}^{\pi}\frac{\cos\theta}{D} \,d\theta = X_2Y_2 \int_{0}^{\pi} \frac{\sin^{2}\theta}{D^{3}}\,d\theta,\quad \int_0^\pi \frac{\cos\theta}{D^3}\,d\theta
=3X_2Y_2\int_0^\pi
\frac{\sin^2\theta}{D^5}\,d\theta.
\]
This yields the representation \eqref{eq:Greenrep}. Differentiating $G(X,Y)$ in $X_1$, it follows that 
\begin{align*}
\partial_{X_1}G(X,Y)=\frac{3(X_2Y_2)^{2}}{2\pi}(Y_1-X_1)\int_{0}^{\pi}\frac{\sin^{2}\theta}{D^{5}}d\theta.
\end{align*}
We estimate $|X_1-Y_1|$ and $Y_2\sin\theta$ in terms of $D$. Since $D(\theta)=D(-\theta)$, we obtain the desired estimate for $\partial_{X_1}G$. The same estimate holds for $\partial_{Y_1}G=-\partial_{X_1}G$. Similarly, we estimate $\partial_{Y_2}G(X,Y)$. This proves \eqref{eq:DanchinGreen}. Since $\partial_{X_2}G(X,Y)=\partial_{Y_2}G(Y,X)$, we obtain \eqref{eq:DanchinGreen2}.

Using the identity 
\begin{align*}
\int_{\bbR} \frac{1}{(t^2 + s^2)^{3/2}} dt  = \frac{2}{s^{2}},
\end{align*}
it follows that  
\begin{align*}
\int_{\mathbb{R}} G(X,Y)\,dX_1 &= \frac{X_2^{2}Y_2^{2}}{\pi} \int_{0}^{\pi} \frac{\sin^2\theta}{X_2^{2}+Y_2^{2}-2X_2Y_2\cos\theta}\,d\theta.
\end{align*}
We further use the identity
\[
\int_{0}^{\pi}\frac{\sin^{2}\theta}{a-b\cos\theta}\,d\theta = \frac{\pi}{b^{2}} \left(a-\sqrt{a^{2}-b^{2}}\right),\quad  a \ge b > 0.
\]
Applying this for $a=X_2^{2} + Y_2^{2}$ and $b = 2X_2Y_2$, we obtain \eqref{eq:green_z}.
Then \eqref{eq:green_strip} follows from
\begin{align*}
\int_{0}^{h}\int_{\mathbb{R}} G(X,Y)\, X_2dX &= \frac{1}{2}\int_{0}^{h} \min\{X_2^{2},Y_2^{2}\}\,X_2\,dX_2  
=
\begin{cases}
\dfrac{h^{4}}{8}, & h\le Y_2,\\[4pt]
\dfrac{Y_2^{2}h^{2}}{4} - \dfrac{Y_2^{4}}{8},
& h>Y_2.
\end{cases}
\end{align*}
\end{proof}

\subsection{Stream function estimates}

We first establish pointwise estimates for the stream function
\eqref{eq:streamfunction} associated with relative vorticity
$\zeta$ satisfying
\[
P[|\zeta|]
=\|r^2\zeta\|_{L^1(\mathbb{R}^2_+;r)}=\|\zeta\|_{L^1_{*,r}}<\infty,
\qquad
Q_p[\zeta]
=\|\zeta\|_{L^p(\mathbb{R}^2_+;r)}=\|\zeta\|_{L^p_r}<\infty,
\]
for $1<p\le\infty$.
Equivalently,
$\zeta\in L^1_{*,r}\cap L^p_r$.
To this end, we use the weighted $L^q$ norm
$\|X_2^{2\theta}\zeta\|_{L^q_r}$.

\begin{proposition}\label{prop:IE}
Let $1<p\le\infty$ and $1\le q\le p$, and define $\theta\in[0,1]$ by
\begin{align}
\frac1q=\theta+\frac{1-\theta}{p},  \label{eq:thetapq}
\end{align}
with the convention $1/\infty=0$. Then, for all
$\zeta\in L^1_{*,r}\cap L^p_r$,
\begin{align}
\|X_2^{2\theta}\zeta\|_{L^q_r}
\le
\|\zeta\|_{L^1_{*,r}}^\theta
\|\zeta\|_{L^p_r}^{1-\theta}.
\label{eq:IE}
\end{align}
\end{proposition}

\begin{proof}
Applying H\"older's inequality for $l=1/(\theta q)$ and $l'=p/(q(1-\theta))$ yields 
\begin{align*}
\int_{\mathbb{R}^{2}_{+}}\left|X_2^{2\theta}\zeta\right|^{q}
=\int_{\mathbb{R}^{2}_{+}}\left|X_2^{2}\zeta\right|^{\theta q}\left|\zeta\right|^{(1-\theta)q} 
\le \left(\int_{\mathbb{R}^{2}_{+}}\left|X_2^{2}\zeta\right|^{\theta ql}\right)^{\frac{1}{l}}
\left(\int_{\mathbb{R}^{2}_{+}}\left|\zeta\right|^{q(1-\theta)l'}\right)^{\frac{1}{l'}}
=||X_2^{2}\zeta||_{L^{1}_r}^{\theta q}||\zeta||_{L^{p}_r}^{(1-\theta)q},
\end{align*}
where the measure $2\pi X_2 dX$ is suppressed. The condition $1/l+1/l'=1$ implies $1/q=\theta+(1-\theta)/p$.
\end{proof}

The following estimates for the stream function provide upper bounds of order $(p-1)^{-1}$ as $p\to 1^+$.

\begin{lemma}\label{lem:str}
There exists a constant $C>0$, independent of $p$,
such that, for all $1<p\le\infty$ and
$\zeta\in L^1_{*,r}\cap L^p_r$,
\begin{align}
|\mathcal{G}\zeta(X)|
&\le C\frac{p}{p-1}
X_2^{\frac{7p-5}{8p-4}}
\|\zeta\|_{L^1_{*,r}}^{
\frac{5p-4}{8p-4}}
\|\zeta\|_{L^p_r}^{
\frac{3p}{8p-4}},
\quad X=(z,r)\in\mathbb{R}^2_+,
\label{eq:str_1}\\
\|\mathcal{G}\zeta\|_{L^\infty_r}
&\le C\frac{p}{p-1}
\|\zeta\|_{L^1_{*,r}}^{
\frac{4p-3}{5p-3}}
\|\zeta\|_{L^p_r}^{
\frac{p}{5p-3}}.
\label{eq:str_2}
\end{align}
For $p=\infty$, the exponents and the factor
$p/(p-1)$ are understood in the limiting sense.
\end{lemma}

\begin{proof}
We may assume that $\zeta$ is nonnegative. For any $L > 0$, we decompose the stream function $\mathcal{G}\zeta$ into two parts:
\[
\mathcal{G}\zeta(X) = \left( \int_{|X-Y| \ge L} + \int_{|X-Y| < L} \right) G(X,Y)  \zeta(Y)Y_2 d Y  =: I_1 + I_2.
\]
Using $\tau = 3/2$ in \eqref{eq:Green_bound}, we estimate 
\[
I_1 \lesssim X_2^2 \int_{|X-Y| \ge L} \frac{Y_2^2}{|X-Y|^3} \zeta(Y) Y_2 d Y \le X_2^2 L^{-3}||r^{2}\zeta||_{L^{1}(\mathbb{R}^{2}_{+};r)}.
\]
Let $1< q< p$. We set $0<\theta<1$ by \eqref{eq:thetapq} and $\tau$ by $\tau+3/2-1/q=2\theta$. Namely, 
\begin{align}
\theta 
= \frac{p - q}{q(p - 1)},\quad 
\tau =
\frac{3p-2q-1}{q(p-1)}-\frac{3}{2}.  \label{eq:thetatau}
\end{align}
The condition $0<\tau<3/2$ holds if and only if 
\begin{align*}
1 < q < \frac{6p-2}{3p+1}. 
\end{align*}
The condition $\tau q'<1$ holds if and only if 
\begin{align*}
\frac{8p-4}{5p-1} < q,
\end{align*}
where $q'$ is the conjugate exponent to $q$. Thus, for the constant $q$ satisfying   
\begin{align}
\frac{8p-4}{5p-1} < q < \frac{6p-2}{3p+1}, \label{eq:cq}
\end{align}
the constant $\tau$ satisfies $0<\tau <3/2$ and $\tau q'<1$.
We choose $q$ so that $\tau q' = 1/2$. Then,
\begin{align}
\tau=\frac{3(p-1)}{2(7p-3)},\quad 
\frac{1}{\tau}=\frac{2(7p-3)}{3(p-1)}\le \frac{14}{3}\left(\frac{p}{p-1}\right). \label{eq:tau}
\end{align}
By the Green function estimate \eqref{eq:Green_bound} and H\"older's inequality, we estimate 
\begin{align*}
\tau I_2 \lesssim X_2^{\tau+\frac{1}{2}} \int_{|X-Y| < L} \frac{Y_2^{\frac{3}{2}+\tau}}{|X-Y|^{2\tau}} \zeta \, d Y &\le X_2^{\tau+\frac{1}{2}} \left( \int_{|X-Y| < L} \frac{1}{|X-Y|^{2\tau q'}} dY \right)^{\frac{1}{q'}} \left( \int_{|X-Y| < L} |Y_2^{\tau+\frac{3}{2}} \zeta |^q dY \right)^{\frac{1}{q}} \\
&\lesssim X_2^{\tau+\frac{1}{2}} L^{\frac{2}{q'}-2\tau} \left\| X_2^{\tau+\frac{3}{2}-\frac{1}{q}} \zeta  \right\|_{L^q_r}
=X_2^{\tau+\frac{1}{2}} L^{\frac{2}{q'}-2\tau} \| X_2^{2\theta} \zeta  \|_{L^q_r}.
\end{align*}
Applying the interpolation inequality \eqref{eq:IE} yields
\begin{align*}
\tau I_2\lesssim X_2^{\tau+\frac{1}{2}} L^{\frac{2}{q'}-2\tau} \|  \zeta  \|_{L^1_{*,r}}^{\theta}\|\zeta  \|_{L^p_r}^{1-\theta}.
\end{align*}
We thus obtain 
\begin{align*}
\tau(I_1+I_2)\lesssim X_2^2 L^{-3}||\zeta||_{L^{1}_{*,r}}+X_2^{\tau+\frac{1}{2}} L^{\frac{2}{q'}-2\tau} \|  \zeta  \|_{L^1_{*,r}}^{\theta}\|\zeta  \|_{L^p_r}^{1-\theta}.
\end{align*}
We optimize the upper bound by choosing $L > 0$ so that the two terms are equal. Namely, 
\begin{align*}
L=X_2^{\frac{\frac{3}{2}-\tau}{\frac{2}{q'}-2\tau+3}}\left(\frac{||\zeta||_{L^{1}_{*,r}}}{||\zeta||_{L^{p}_r}}\right)^{\frac{1-\theta}{\frac{2}{q'}-2\tau+3}}.
\end{align*}
Using 
\begin{align}
\frac{2}{q'}-2\tau+3=\frac{4(2p-1)(q-1)}{q(p-1)},\quad 1-\theta=\frac{p(q-1)}{q(p-1)},\quad \frac{3}{2}-\tau =\frac{(3p-1)(q-1)}{q(p-1)},    \label{eq:thetataupq}
\end{align}
we observe that  
\begin{align*}
\frac{\frac{3}{2}-\tau}{\frac{2}{q'}-2\tau+3}=\frac{3p-1}{4(2p-1)},\quad \frac{1-\theta}{\frac{2}{q'}-2\tau+3}=\frac{p}{4(2p-1)}, 
\end{align*}
and $L$ is independent of the choice of $q$. We obtain the estimate \eqref{eq:str_1}.

We choose $L = X_2/2$ and estimate $I_1$ and $I_2$. We first interpolate two inequalities 
\begin{align*}
I_1&\lesssim X_2^{-1}||\zeta||_{L^{1}_{*,r}}\le \frac{p}{p-1} X_2^{-1}||\zeta||_{L^{1}_{*,r}},\\
I_1&\lesssim \frac{p}{p-1} X_2^{\frac{7p-5}{8p-4}} \| \zeta  \|_{L^1_{*,r}}^{\frac{5p-4}{8p-4}} \|\zeta\|_{L^p_r}^{\frac{3p}{8p-4}}.
\end{align*}
For $0<s<1$, it follows that 
\begin{align*}
I_1=I_1^{s}I_1^{1-s}\lesssim \frac{p}{p-1} X_2^{-1+\frac{15p-9}{8p-4}s}|| \zeta||_{L^1_{*,r}}^{1-\frac{3p}{8p-4}s}||\zeta||_{L^p_r}^{\frac{3p}{8p-4}s}.
\end{align*}
We take $s=\frac{8p-4}{15p-9}$ so that $-1+\frac{15p-9}{8p-4}s=0$ and obtain 
\begin{align*}
I_1\lesssim \frac{p}{p-1} || \zeta||_{L^1_{*,r}}^{\frac{4p-3}{5p-3}}||\zeta||_{L^p_r}^{\frac{p}{5p-3}}.
\end{align*}
Let $1<q<p$. We choose $0<\tau\le 3/2$ sufficiently small so that $\tau q'<1$ and apply H\"older's inequality to estimate 
\begin{align*}
\tau I_2 \lesssim X_2^{\tau+\frac{1}{2}} \int_{|X-Y| < L} \frac{Y_2^{\frac{1}{2}+\tau}}{|X-Y|^{2\tau}} \zeta \, Y_2dY 
\le X_2^{\tau+\frac{1}{2}} \left( \int_{|X-Y| < X_2/2} \frac{Y_2^{(\frac{1}{2}+\tau)q'}}{|X-Y|^{2\tau q'}} Y_2 dY \right)^{\frac{1}{q'}}
||\zeta||_{L^{q}(|X-Y|<X_2/2; r)}
\end{align*}
The right-hand side is a constant multiple of $X_2^{1+\frac{3}{q'}}||\zeta||_{L^{q}(|X-Y|<X_2/2; r)}$. We take $0<\theta<1$ such that $1/q=\theta+(1-\theta)/p$. Using $1/2<Y_2/X_2$ for $|X-Y|<X_2/2$ and applying \eqref{eq:IE}, we estimate  
\begin{align*}
\tau I_2 \lesssim X_2^{1+\frac{3}{q'}}||\zeta||_{L^{q}(|X-Y|<X_2/2; r)}
\lesssim X_2^{1+\frac{3}{q'}-2\theta}||X_2^{2\theta}\zeta||_{L^{q}_r}
\lesssim X_2^{1+\frac{3}{q'}-2\theta} || \zeta ||_{L^{1}_{*,r}}^{\theta}||\zeta||_{L^{p}_r}^{1-\theta}.
\end{align*}
We take $q$ so that $1+3/q'-2\theta=0$.
Namely, $q=\frac{5p-3}{4p-2}$ for which $\theta=\frac{4p-3}{5p-3}$.
We obtain 
\begin{align*}
\tau I_2
\lesssim  || \zeta||_{L^1_{*,r}}^{\frac{4p-3}{5p-3}}||\zeta||_{L^p_r}^{\frac{p}{5p-3}}.
\end{align*}
We choose $\tau$ so that $\tau q'=1/2$. Namely, $\tau=\frac{p-1}{2(5p-3)}$ and 
\begin{align*}
\frac{1}{\tau}\le 2\frac{5p-3}{p-1}\le 10\frac{p}{p-1}.
\end{align*}
This completes the proof of \eqref{eq:str_2}.
\end{proof}
	
\subsection{Kinetic energy estimates}

We define the interaction energy 
\begin{equation}\label{eq:en_int}
E_{\mathrm{int}}[\zeta,\eta]
=\pi\iint_{\mathbb{R}^2_+\times\mathbb{R}^2_+}
G(X,Y)\,
\zeta(X)\eta(Y)\,
X_2Y_2\,dX\,dY,
\end{equation}
so that $E[\zeta]=E_{\mathrm{int}}[\zeta,\zeta]$.
We establish the variational inequality \eqref{eq:variational-inequality}.

\begin{lemma}\label{lem:energy}
Let $1<p\le\infty$, and let $\alpha_p$ and
$\beta_p$ be defined by \eqref{eq:exponent}. There exists a constant $C>0$, independent of $p$,
such that for all
$\zeta,\eta\in L^1_{*,r}\cap L^p_r$,
\begin{align}
E[\zeta]
&\le C\frac{p}{p-1}
\|\zeta\|_{L^1_{*,r}}^{\alpha_p}
\|\zeta\|_{L^p_r}^{\beta_p},
\label{eq:en_bound}\\
|E[\zeta]-E[\eta]|
&\le C\frac{p}{p-1}
\|\zeta+\eta\|_{L^1_{*,r}}^{\frac{\alpha_p}{2}}
\|\zeta+\eta\|_{L^p_r}^{\frac{\beta_p}{2}}
\|\zeta-\eta\|_{L^1_{*,r}}^{\frac{\alpha_p}{2}}
\|\zeta-\eta\|_{L^p_r}^{\frac{\beta_p}{2}},
\label{eq:en_diff}\\
|E_{\mathrm{int}}[\zeta,\eta]|
&\le C\frac{p}{p-1}
\|\zeta\|_{L^1_{*,r}}^{\frac{\alpha_p}{2}}
\|\zeta\|_{L^p_r}^{\frac{\beta_p}{2}}
\|\eta\|_{L^1_{*,r}}^{\frac{\alpha_p}{2}}
\|\eta\|_{L^p_r}^{\frac{\beta_p}{2}}.
\label{eq:en_cross}
\end{align}
For $p=\infty$, the exponents and the factor
$p/(p-1)$ are understood in the limiting sense.
\end{lemma}

\begin{proof}
Using the symmetry $G(X, Y) = G(Y, X)$, we observe that
\begin{align*}
E[\zeta] - E[\eta]  
&= \iint_{\mathbb{R}^{2}_{+}\times\mathbb{R}^{2}_{+}} G(X, Y) \zeta(X)\zeta(Y)
- \iint_{\mathbb{R}^{2}_{+}\times\mathbb{R}^{2}_{+}} G(X, Y) \eta(X)\eta(Y) \\
&= \iint_{\mathbb{R}^{2}_{+}\times\mathbb{R}^{2}_{+}} G(X, Y) \zeta(X)(\zeta(Y)-\eta(Y))
+ \iint_{\mathbb{R}^{2}_{+}\times\mathbb{R}^{2}_{+}} G(X, Y) (\zeta(X)-\eta(X))\eta(Y) \\
&=\iint_{\mathbb{R}^{2}_{+}\times\mathbb{R}^{2}_{+}} G(X, Y) (\zeta(X)+\eta(X)) (\zeta(Y)-\eta(Y))=E_{\textrm{int}}[\zeta+\eta,\zeta-\eta],
\end{align*}
where the measure $\pi X_2Y_2dXdY$ is omitted for brevity. Thus the estimate \eqref{eq:en_diff} and \eqref{eq:en_bound} follow from \eqref{eq:en_cross}.   

We may assume that $\zeta$ and $\eta$ are nonnegative and defined on $\mathbb{R}^{2}$ by the zero extension. We take an arbitrary $L>0$ and split $E_{\textrm{int}}[\zeta, \eta]$ into the far-field and near-field regions:
\begin{align*}
\frac{1}{\pi}E_{\textrm{int}}[\zeta, \eta]&=
\iint_{|X-Y| \ge L} X_2Y_2G(X,Y) \zeta(X) \eta(Y)   dX dY+\iint_{|X-Y| < L} X_2Y_2G(X,Y) \zeta(X) \eta(Y)  dX dY\\
&=:I_1+I_2.
\end{align*}
Using the Green function estimate \eqref{eq:Green_bound} for $\tau=3/2$,
\[
I_1 \lesssim \iint_{|X-Y| \ge L} \frac{(X_2Y_2)^3}{|X-Y|^3} \zeta(X) \eta(Y)  dX  dY \le \frac{1}{L^{3}} ||\zeta||_{L^{1}_{*,r}}||\eta||_{L^{1}_{*,r}}.
\]
Let $1<q<p$.
We define $\theta$ and $\tau$ by \eqref{eq:thetatau}, so that $1/q = \theta+(1-\theta)/p$ and $\tau+3/2-1/q=2\theta$, as in the proof of \eqref{eq:str_1}.
For $q$ satisfying \eqref{eq:cq}, $0<\tau<3/2$ and $q'\tau<1$ for the conjugate exponent $q'$ of $q$.
We choose $\tau$ by \eqref{eq:tau}. Using \eqref{eq:Green_bound} and H\"older's inequality,
\begin{align*}
\tau I_2 \lesssim \iint_{|X-Y|< L} \frac{(X_2Y_2)^{\frac{3}{2}+\tau }}{|X-Y|^{2\tau}} \zeta(X)\eta(Y)  dX dY 
&=\int_{\mathbb{R}^{2}_{+}}\left(\int_{|X-Y|< L}\frac{(Y_2)^{\frac{3}{2}+\tau }}{|X-Y|^{2\tau}} \eta(Y)   dY \right)(X_2)^{\frac{3}{2}+\tau }\zeta(X)dX \\
&=\int_{\mathbb{R}^{2}_{+}}\left(\left(1_{\{|X|<L\}}\frac{1}{|X|^{2\tau}}\right)*(X_2)^{\frac{3}{2}+\tau }\eta\right)(X_2)^{\frac{3}{2}+\tau }\zeta(X)dX \\
&\le \left\|\left(1_{\{|X|<L\}}\frac{1}{|X|^{2\tau}}\right)*(X_2)^{\frac{3}{2}+\tau }\eta\right\|_{L^{q'}(\mathbb{R}^{2}_{+})}\left\|  (X_2)^{\frac{3}{2}+\tau }\zeta\right\|_{L^{q}(\mathbb{R}^{2}_{+})}.
\end{align*}
We apply Young's convolution inequality and estimate
\begin{align*}
\left\|\left(1_{\{|x|<L\}}\frac{1}{|X|^{2\tau}}\right)*(X_2)^{\frac{3}{2}+\tau }\eta\right\|_{L^{q'}(\mathbb{R}^{2}_{+})}
\le \left\|1_{\{|X|<L\}}\frac{1}{|X|^{2\tau}}\right\|_{L^{k}(\mathbb{R}^{2})}
\left\|(X_2)^{\frac{3}{2}+\tau }\eta\right\|_{L^{l}(\mathbb{R}^{2}_{+})},
\end{align*}
for $1/q'=1/k+1/l-1$.
We choose $l=q$ and $k=q'/2$ for which the right-hand side is a constant multiple of $L^{\frac{4}{q'}-2\tau}||X_2^{2\theta}\eta||_{L^{q}_r}$.
Applying the inequality \eqref{eq:IE} for $\zeta$ and $\eta$ yields 
\begin{align*}
\tau I_2\lesssim L^{\frac{4}{q'}-2\tau}||\eta||_{L^{1}_{*,r}}^{\theta}||\eta||_{L^{p}_r}^{1-\theta}
||\zeta||_{L^{1}_{*,r}}^{\theta}||\zeta||_{L^{p}_r}^{1-\theta}.
\end{align*}
We thus obtain 
\begin{align*}
\tau (I_1+I_2)\lesssim \frac{1}{L^{3}} ||\zeta||_{L^{1}_{*,r}}||\eta||_{L^{1}_{*,r}}+L^{\frac{4}{q'}-2\tau}||\eta||_{L^{1}_{*,r}}^{\theta}||\eta||_{L^{p}_{r}}^{1-\theta}
||\zeta||_{L^{1}_{*,r}}^{\theta}||\zeta||_{L^{p}_r}^{1-\theta}.
\end{align*}
We take $L>0$ so that the two terms are equal.
Namely, 
\begin{align*}
L=\left(\frac{||\zeta||_{L^{1}_{*,r}}||\eta||_{L^{1}_{*,r}}} {||\zeta||_{L^{p}_r}||\eta||_{L^{p}_r}} \right)^{\frac{1-\theta}{\frac{4}{q'}-2\tau+3}}.
\end{align*}
It follows from \eqref{eq:thetataupq} and 
\begin{align*}
\frac{4}{q'}-2\tau+3=\frac{2}{q(p-1)}(5p-3)(q-1)
\end{align*}
that  
\begin{align*}
\frac{1-\theta}{\frac{4}{q'}-2\tau+3}=\frac{p}{2(5p-3)}.
\end{align*}
Thus $L$ is independent of the choice of $q$. We obtain the estimate \eqref{eq:en_cross}. 
\end{proof}

We show the identity \eqref{eq:ke2} for $\zeta \in L^{1}_{*,r}\cap L^p_r$.

\begin{lemma}\label{lem:str_energy}
Let $1<p\le \infty$. 
For $ \zeta,\eta \in L^{1}_{*,r}\cap L^p_r$, we have
\begin{align}
\sup_{|X| \ge L} |\mathcal{G}\zeta(X)| &\to 0 \quad \text{as } \,\, L\to\infty, \label{eq:vanishing} \\
E[\zeta] &= \pi \int_{\mathbb{R}^{2}_{+}} \frac{1}{X_2} \left|\nabla_{X} \mathcal{G}\zeta\right|^2  \,dX, \label{eq:L2energy}\\
E_{\rm int}[\zeta,\eta] &= \pi \int_{\mathbb{R}^{2}_{+}} \frac{1}{X_2} \nabla_{X} \mathcal{G}\zeta\cdot \nabla_{X} \mathcal{G}\eta   \,dX, \label{eq:L2energy2}\\
|E_{\rm int}[\zeta,\eta]| &\le \sqrt{E[\zeta]}\sqrt{E[\eta]},\label{eq:p1_en_CS}\\
\sqrt{E[\zeta+\eta]} &\le \sqrt{E[\zeta]} + \sqrt{E[\eta]}. \label{eq:p1_en_tri}
\end{align}
\end{lemma}

\begin{proof}
By the definition of the energy, we have 
\begin{align*}
E[\zeta+\eta]=E[\zeta]+E[\eta]+2E_{\textrm{int}}[\zeta,\eta].    
\end{align*}
The representation \eqref{eq:L2energy2} follows from \eqref{eq:L2energy} since  
\begin{align*}
E[\zeta+\eta]=E[\zeta]+E[\eta]+2\pi \int_{\mathbb{R}^{2}_{+}} \frac{1}{X_2} \nabla_{X} \mathcal{G}\zeta\cdot \nabla_{X} \mathcal{G}\eta   \,dX.
\end{align*}
The estimates \eqref{eq:p1_en_CS}-\eqref{eq:p1_en_tri} follow from \eqref{eq:L2energy2} and the Cauchy--Schwarz inequality. 

It suffices to show \eqref{eq:vanishing} and \eqref{eq:L2energy}. We take an arbitrary $\varepsilon>0$. Using a nonnegative cutoff function $\chi_N(X) = \chi(|X|/N)\le 1$, we set $\zeta_N=\zeta\chi_{N}$ so that $|\zeta_N| \le |\zeta|$, and $\zeta_N \to \zeta$ in $L_{*,r}^1$. We fix a sufficiently large $N=N_{\varepsilon}$ such that 
\[
\| \zeta - \zeta_N  \|_{L^1_{*,r}} \le \varepsilon.
\]
Since $\zeta_N$ has compact support, we have $\supp \zeta_N \subset B_{L} = \{X\in \mathbb{R}^{2}_{+}: |X| < L\}$ for large $L >0$. For $|X| > 2L$ and $Y \in \supp \zeta_N$, we have $|X-Y| \ge |X|/2$.
Thus, 
\[
|\mathcal{G}\zeta_N(X)| = \left| \int_{\mathbb{R}^{2}_{+}} G(X, Y) Y_2 \zeta_N(Y) \, dY \right| \lesssim \int_{B_{L}} \frac{X_2^2 Y_2^3}{|X|^3} |\zeta_N(Y)| \, dY \le \frac{1}{2L} ||\zeta||_{L^{1}_{*,r}},
\]
where we used \eqref{eq:Green_bound} with $\tau = 3/2$. Applying \eqref{eq:str_2} to $\mathcal{G}(\zeta - \zeta_N)$, we obtain
\begin{align*}
|\mathcal{G}\zeta(X)| &\le |\mathcal{G}(\zeta - \zeta_N)(X)| + |\mathcal{G}\zeta_N(X)| \lesssim_p \varepsilon^{\frac{4p-3}{5p-3}} + \frac{1}{2L} ||\zeta||_{L^{1}_{*,r}},
\end{align*}
for $|X| \ge 2 L$. Taking the supremum over $|X|\ge 2L$ and letting $L \to \infty$ and $\varepsilon\to0$, we obtain \eqref{eq:vanishing}. 

It remains to show \eqref{eq:L2energy}.
We define a sequence of bounded, compactly supported functions by
\begin{align*}
\tilde{\zeta}_N=T_N(\zeta_N),
\qquad
T_N(s)=\operatorname{sgn}(s)\min\{|s|,N\}.
\end{align*}
The identity \eqref{eq:ke2} holds for $\tilde{\zeta}_{N}$.
Applying \eqref{eq:ke2} for $\tilde{\zeta}_{N}-\tilde{\zeta}_{M}$, it follows that $\mathcal{G}\tilde{\zeta}_{N}$ converges in $\dot{H}^1_0(\mathbb{R}^{2}_{+}; r^{-1})$. Since $\tilde{\zeta}_{N}\to \zeta$ in $L^{1}_{*,r}\cap L^{p}_{r}$, we have $\mathcal{G}\tilde{\zeta}_{N}\to \mathcal{G}\zeta$ in $L^{\infty}(\mathbb{R}^{2}_{+})$ by \eqref{eq:str_2}. Thus $\mathcal{G}\tilde{\zeta}_{N}$ converges to $\mathcal{G}\zeta$ in $\dot{H}^1_0(\mathbb{R}^{2}_{+}; r^{-1})$. Applying \eqref{eq:ke2} for $\tilde{\zeta}_{N}$ yields  
\begin{align*}
E[\tilde{\zeta}_{N}]=\pi \int_{\mathbb{R}^{2}_{+}}\frac{1}{X_2}|\nabla_X\mathcal{G}\tilde{\zeta}_N|^{2}dX.
\end{align*}
Letting $N\to\infty$ implies \eqref{eq:L2energy}.
\end{proof}

\subsection{Strict superadditivity}

By the variational inequality \eqref{eq:variational-inequality}, for every $1<p\le\infty$ and $\mu,\nu>0$, the constant $I_{p,\mu,\nu}$ defined by \eqref{eq:variational} is finite and positive.
We deduce the following strict superadditivity from the scaling property \eqref{eq:scaling}.

\begin{lemma}[Sharp energy inequality]
Let $1<p\le\infty$, and let $\alpha_p$ and $\beta_p$ be defined by \eqref{eq:exponent}.
Then the optimal constant $I_p=I_{p,1,1}\in(0,\infty)$ defined in \eqref{eq:variational} is the best constant in the inequality
\begin{equation}\label{eq:variational-inequality2}
E[\zeta]
\le I_p P[|\zeta|]^{\alpha_p}
Q_p[\zeta]^{\beta_p}
\end{equation}
for every $\zeta\in L^1_{*,r}\cap L^p_r$.
\end{lemma}

\begin{lemma}[Strict superadditivity]\label{lem:superadd}
Let $1<p\le \infty$ and $\mu_i, \nu_i > 0$ for $i \in \{1, 2\}$.
The optimal constant satisfies
\begin{align}
I_{p,\mu_1, \nu_1} + I_{p,\mu_2, \nu_2} < I_{p,\mu_1+\mu_2, (\nu_1^p + \nu_2^p)^{1/p}}, \quad & p \in (1, \infty), \label{eq:superadd_p} \\
I_{\infty,\mu_1, \nu_1} + I_{\infty,\mu_2, \nu_2} < I_{\infty,\mu_1+\mu_2, \max\{\nu_1, \nu_2\}}, \quad & p = \infty. \label{eq:superadd_inf}
\end{align}
\end{lemma}

\begin{proof}
We use the scaling property \eqref{eq:scaling}. For $p=\infty$, $\alpha_{\infty}>1$ and  
\begin{align*}
\mu_{1}^{\alpha_{\infty}}\nu_{1}^{\beta_{\infty}}
+\mu_{2}^{\alpha_{\infty}}\nu_{2}^{\beta_{\infty}}
\le (\mu_{1}^{\alpha_{\infty}}+
\mu_{2}^{\alpha_{\infty}})\max\{\nu_{1},\nu_{2}\}^{\beta_{\infty}}
<(\mu_{1}+\mu_{2})^{\alpha_{\infty}}\max\{\nu_{1},\nu_{2}\}^{\beta_{\infty}}.
\end{align*}
Multiplying this inequality by $I_{\infty,1,1}$ yields \eqref{eq:superadd_inf}.

For $1<p< \infty$, we set $\alpha=\alpha_p$ and $\beta=\beta_p$ by \eqref{eq:exponent}.
For any $t, s \in (0, 1)$ and an exponent $q \in (1, \infty)$ with its conjugate $q'$, we use Young's inequality 
\[
t^\alpha s^{\frac{\beta}{p}} \le \frac{1}{q} t^{\alpha q} + \frac{1}{q'} s^{ \frac{\beta q'}{p}}.
\]
Since 
\[
1-\frac{\beta}{p}
=\frac{5p-6}{5p-3} < \min\{\alpha,1\},
\]
there exists $q\in (1,\infty)$ such that 
\[
1 - \frac{\beta}{p} < \frac{1}{q} < \alpha.
\]
Using convexity of $t^{\sigma}$ for $\sigma=\alpha q > 1$ and $\beta q'/p > 1$, it follows that
\[
t^\alpha s^{\frac{\beta}{p}} + (1-t)^\alpha (1-s)^{\frac{\beta}{p}} \le \frac{t^{\alpha q} + (1-t)^{\alpha q}}{q}  + \frac{s^{ \frac{\beta q'}{p}} + (1-s)^{\frac{\beta q'}{p}}}{q'} < \frac{1}{q} + \frac{1}{q'} = 1.
\]
By substituting $t = \mu_1/(\mu_1+\mu_2)$ and $s = \nu_1^p/(\nu_1^p+\nu_2^p)$, we arrive at
\[
\mu_1^\alpha \nu_1^\beta + \mu_2^\alpha \nu_2^\beta < (\mu_1+\mu_2)^\alpha (\nu_1^p+\nu_2^p)^{\frac{\beta}{p}}.
\]
Multiplying both sides by $I_{p,1, 1}$ completes the proof of \eqref{eq:superadd_p}.
\end{proof}

\begin{remark}\label{r:strictp_1}
The last inequality also holds for $p=1$.
Namely, 
\begin{align*}
 \mu_1^{1/2}\nu_1^{3/2} + \mu_2^{1/2}\nu_2^{3/2}<(\mu_1+\mu_2)^{1/2}(\nu_1+\nu_2)^{3/2}. 
\end{align*}
The inequality remains strict for nonnegative $\mu_i,\nu_i$, provided that $\mu_i+\nu_i>0$ for $i=1,2$ and $(\mu_1+\mu_2)(\nu_1+\nu_2)>0$. We will use this inequality to study the asymptotics of $I_p$ as $p\to 1$ in \S \ref{s:7}. 
\end{remark}

\begin{remark}[Saturation of the constraints]\label{sec:saturation}
Any maximizer of the variational problem \eqref{eq:variational} must saturate both constraints. Namely, for the set of maximizers $S_{p,\mu, \nu}$, we have 
\begin{equation*}
S_{p,\mu, \nu} = \{ \zeta \in K_{p,\mu, \nu} : E[\zeta] = I_{p,\mu, \nu} \} = \{ \zeta \ge 0 : E[\zeta] = I_{p,\mu, \nu}, \,\, P[\zeta] = \mu, \,\, Q_p[\zeta] = \nu \}.
\end{equation*}
Indeed, a maximizer $\zeta$ satisfying
\begin{align*}
E[\zeta]=I_{p,\mu,\nu},\quad P[\zeta]=\mu_1\le \mu,\quad Q_p[\zeta]=\nu_1\le \nu,
\end{align*}
must obey  
\begin{align*}
I_{p,\mu,\nu}=E[\zeta]\le I_{p,\mu_1,\nu_1}\le I_{p,\mu,\nu}. 
\end{align*}
Since $I_{p,\mu,\nu}$ is strictly increasing for $\mu>0$ and $\nu>0$, we find that $\mu_1=\mu$ and $\nu_1=\nu$. 

Similarly, any nonnegative sequence $\{\zeta_n\}$ satisfying
\[
I_{p,\mu, \nu} \le \liminf_{n \to \infty} E[\zeta_n], \quad \limsup_{n \to \infty} P[\zeta_n] \le \mu, \quad \limsup_{n \to \infty} Q_p[\zeta_n] \le \nu,
\]
necessarily obeys
\[
I_{p,\mu, \nu} = \lim_{n \to \infty} E[\zeta_n],\quad 
\lim_{n \to \infty} P[\zeta_n] = \mu, \quad \lim_{n \to \infty} Q_p[\zeta_n] = \nu.
\]
\end{remark}

\subsection{The case \texorpdfstring{$p=\infty$}{p=infinity}: Hill's vortex}

We show that all maximizers of the variational problem \eqref{eq:variational} for $p=\infty$ are translations of Hill's vortex \eqref{eq:normalizedHill0}.

\begin{proof}[Proof of Theorem \ref{t:MT3}(iii)]
Let $\zeta\in S_\infty$ be a maximizer.
By Remark \ref{sec:saturation}, $P[\zeta] = 1$ and $ Q_{\infty}[\zeta]= 1$.
It follows from \eqref{eq:kineticenergy5D} and \eqref{eq:impulse5D} that 
\[
E[\zeta] = \frac{1}{16\pi^3} \iint_{\bbR^5\times\bbR^5} \frac{\zeta(w)\zeta(y)}{|w-y|^3}\,dwdy = I_\infty,\quad 
\int_{\bbR^5}\zeta\,dy = \pi,\quad 
0 \le \zeta \le 1.
\]
Let $\zeta^*(y)=\zeta^*(|y|)$ be the symmetric decreasing rearrangement of $\zeta(y)$ in $\bbR^5$. By the Riesz rearrangement inequality \cite[Theorem 3.7]{LL},
\begin{align}
\iint_{\bbR^5\times\bbR^5} \frac{\zeta(w)\zeta(y)}{|w-y|^3}\,dwdy \le \iint_{\bbR^5\times\bbR^5} \frac{\zeta^*(w)\zeta^*(y)}{|w-y|^3}\,dwdy.
\label{eq:riesz}
\end{align}
Since $\zeta^*$ and $N(y)=(8\pi^{2})^{-1}|y|^{-3}$ are radial, the Newtonian
potential $N*\zeta^*$ is also radially symmetric. We therefore define radially symmetric functions
\[
m(\varrho)=\int_{B_\varrho}\zeta^*(y)\,dy = \frac{8\pi^2}{3}\int_0^\varrho\zeta^*(t)t^4\,dt, \quad \varphi(\varrho)=N*\zeta^{*}(\varrho)=\frac{1}{8\pi^{2}}\int_{\bbR^5}\frac{\zeta^*(y)}{|w-y|^3}\,dy \quad\textrm{for}\,\,|w|=\varrho,
\]
where $B_\varrho=B_\varrho(0)$ denotes the open ball in $\bbR^5$ of radius $\varrho$ centered at the origin. Since $-\Delta_y \varphi = \zeta^*$ and $\zeta^*\in L^1\cap L^\infty(\bbR^5)$, $\varphi \in W^{2,q}_{\rm loc}(\bbR^5)$ for every $1 < q < \infty$.
Thus, for a.e. $\varrho>0$,
\[
m(\varrho) = -\int_{B_\varrho}\Delta \varphi\,dy = -\int_{\partial B_\varrho}\partial_\nu \varphi\,d\sigma = - \frac{8\pi^2}{3} \varrho^4\varphi'(\varrho).
\]
Moreover, $\varphi(\varrho) \to 0$ as $\varrho \to \infty$, and hence
\[
\varphi(\varrho) = -\int_\varrho^\infty \varphi'(t)\,dt = \frac{3}{8\pi^{2}}\int_\varrho^\infty\frac{m(t)}{t^4}\,dt.
\]
Therefore, by Tonelli's theorem with $m(0)=0$,
\begin{equation}\label{eq:radial}
\frac{1}{8\pi^{2}}\iint_{\bbR^5\times\bbR^5} \frac{\zeta^*(w)\zeta^*(y)}{|w-y|^3}\,dwdy = \int_0^\infty \varphi(\varrho)\,dm(\varrho)  = \frac{3}{8\pi^{2}}\int_0^\infty\frac{m(\varrho)^2}{\varrho^4}\,d\varrho.
\end{equation}
Since $0\le \zeta^* \le 1$, $\int_{\bbR^5} \zeta^*\,dy \le \pi$, and $|B_1| = 8\pi^2/15$,
\begin{equation}\label{eq:rad_mass}
m(\varrho) \le \min\left\{ \pi,\frac{8\pi^2}{15}\varrho^5 \right\} = \int_{B_\varrho}\mathbf{1}_{B_{\varrho_H}}(y)\,dy = \int_{B_\varrho} \zeta_H(y)\,dy =: m_H(\varrho).
\end{equation}
It follows from \eqref{eq:radial} that
\[
\iint_{\bbR^5\times\bbR^5} \frac{\zeta^*(w)\zeta^*(y)}{|w-y|^3}\,dwdy \le \iint_{\bbR^5\times\bbR^5} \frac{\zeta_H(w)\zeta_H(y)}{|w-y|^3}\,dwdy.
\]
This gives $E[\zeta]\le E[\zeta_H]$ for any $\zeta \in S_{\infty}$.
By $\zeta_H\in K_{\infty}$, we have $\zeta_H\in S_{\infty}$ and 
\begin{align*}
I_\infty=E[\zeta_H] = \frac{3}{16\pi^3}\int_0^\infty\frac{m_H(\varrho)^2}{\varrho^4}\,d\varrho=  \frac{3}{16\pi^3} \left(\frac{8\pi^2}{15}\right)^2 \left( \int_0^{\varrho_H}\varrho^6\,d\varrho + \varrho_H^{10}\int_{\varrho_H}^\infty \varrho^{-4}\,d\varrho \right) =\frac{1}{21}\varrho_H^{2}.
\end{align*}
It remains to prove uniqueness. Let $\zeta\in S_{\infty}$, and let $m$ be defined as above.
Then equality holds in all of the preceding inequalities.
In particular, \eqref{eq:radial} yields
\[
E[\zeta_H] - E[\zeta^*] = \frac{3}{16\pi^3}\int_0^\infty \frac{1}{\varrho^4} \left( m_H(\varrho)^2 - m(\varrho)^2 \right)d\varrho = 0.
\]
The integrand is nonnegative by \eqref{eq:rad_mass}, so $m(\varrho) = m_H(\varrho)$ for a.e. $\varrho > 0$, and hence for every $\varrho > 0$ by continuity.
Therefore, for every $\varrho<\varrho_H$,
\[
\int_{B_\varrho}(1-\zeta^*)\,dy=0.
\]
It follows that $\zeta^*=1$ almost everywhere in $B_{\varrho_H}$.
On the other hand, $m(\varrho) = m_H(\varrho) = \pi$ for every $\varrho>\varrho_H$, and hence
\[
\int_{\bbR^5}\zeta^*\,dy = \lim_{\varrho \to \infty}m(\varrho) = \pi = |B_{\varrho_H}|.
\]
Therefore, $\zeta^*=0$ almost everywhere in $\bbR^5\setminus B_{\varrho_H}$, and thus $\zeta^* = \mathbf{1}_{B_{\varrho_H}}$.
Since $\zeta \in S_{\infty}$, equality also holds in \eqref{eq:riesz}.
By the equality condition for rearrangement inequality \cite[Theorem 3.9]{LL} applied to $\zeta$ and $|\cdot|^{-3}$, there exists $y_0 \in \bbR^5$ such that  $\zeta(y) = \zeta^*(y-y_0) = \mathbf{1}_{B_{\varrho_H}(y_0)}$.
Since $\zeta$ is invariant under rotations in the $y'$-variables, the center $y_0$ must lie on the $z$-axis.
Hence $y_0=(c,0)$ for some $c\in\bbR$, and therefore $\zeta=\zeta_{H}(\cdot+c,\cdot)$.
Conversely, since translations in the $z$-direction preserve $P[\cdot]$, $Q_{\infty}[\cdot] $, and $E[\cdot]$, for every $c \in \bbR$ we have $\zeta_{H}(\cdot+c,\cdot) \in S_{\infty}$.
This completes the proof.
\end{proof}

\begin{remark}\label{r:hill_EL}
Hill's vortex $\zeta_H$ is expressed as 
\begin{equation}\label{eq:EL_hill}
\zeta_H = \mathbf{1}_{\left\{\mathcal{G}\zeta_H - \frac{W_{\infty}}{2}r^2>0\right\}} \quad \textrm{a.e. in}\,\,\bbR^2_+,
\end{equation}
by the constant $W_{\infty}$ in \eqref{eq:speedHill}. Moreover, its stream function is expressed as 
\begin{align}
\frac{\mathcal{G}\zeta_H}{r^{2}}=
\begin{cases}
\displaystyle
\frac{\varrho_H^2}{6} - \frac{\varrho^2}{10}, & \varrho<\varrho_H,
\\[2mm]
\displaystyle
\frac{\varrho_H^5}{15 \varrho^3}, & \varrho \ge \varrho_H,
\end{cases}
\label{eq:5Dstreamhill}
\end{align}
where $\varrho^{2}=z^{2}+r^{2}$. These representations follow from the five-dimensional potential representation $\mathcal{G}\zeta_H/r^{2}=N*\zeta_{H}$.
\end{remark}

\section{Structure of energy maximizers}\label{s:3}

In this section, we show that energy maximizers of the variational problem \eqref{eq:variational} for $1<p<\infty$ satisfy the Euler--Lagrange equation \eqref{eq:EL_intro0} with the constants determined by the maximum energy \eqref{eq:speed}.
We further prove that energy maximizers have compact support and are Steiner symmetric up to translation in the $z$-direction.
By scaling \eqref{eq:scaling}, we reduce \eqref{eq:variational} to the case $\mu = \nu = 1$ and use a shorthand notation $I_p=I_{p,1,1}$, $K_p = K_{p, 1, 1}$, and $S_p=S_{p,1,1}$.
We denote by $L^{\infty}_{c}(\mathbb{R}^{2}_{+})$ the space of compactly supported functions in $L^{\infty}(\mathbb{R}^{2}_{+})$.

\subsection{Euler--Lagrange Equation}

We use the symbol 
\begin{align}
(f,g)=2\pi \int_{\mathbb{R}^{2}_{+}} f(z,r)\,g(z,r)\, rdzdr,   \label{eq:product}
\end{align}
for two functions $f$ and $g$ and express the conserved quantities as
\begin{align}
E[\zeta]=\frac{1}{2}(\mathcal{G}\zeta,\zeta),\quad P[\zeta]=(r^{2},\zeta),\quad Q_p^{p}[\zeta]=(\zeta^{p-1},\zeta). \label{eq:innerproductexpression}
\end{align}
We first find the Lagrange multipliers $W \ge 0$ and $\kappa > 0$ for each maximizer in $S_p$.

\begin{lemma}\label{lem:EL_p}
Let $1<p<\infty$.
For any $\zeta \in S_p$, there exist unique constants $W \ge 0$ and $\kappa > 0$ such that \eqref{eq:EL_intro0} holds.
\end{lemma}

\begin{proof}
By the energy estimates \eqref{eq:en_bound} and \eqref{eq:en_diff}, the functional $E:L^{1}_{*,r}\cap L^{p}_r \to \mathbb{R}$ is continuous. Furthermore, for $\zeta,\eta\in L^{1}_{*,r}\cap L^{p}_r$, we have 
\begin{equation}\label{eq:frechet0}
 E[\zeta + \tau \eta] 
 =E[\zeta ]+2\tau E_{\textrm{int}}[\zeta,\eta]+\tau^{2}E[\eta ].
\end{equation}
In particular, the directional derivative is evaluated as
\begin{equation}\label{eq:frechet}
\langle E'[\zeta], \eta \rangle = 2E_{\textrm{int}}[\zeta,\eta].
\end{equation}
By the interaction energy estimate \eqref{eq:en_cross}, $E':L^{1}_{*,r}\cap L^{p}_r\to (L^{1}_{*,r}\cap L^{p}_r)^{*}$ is continuous. Thus $E$ is Fr\'echet differentiable. 

Let $\zeta\in S_p$ be a maximizer. Since $E[\zeta]=I_p>0$, there exists a small constant $\delta_0 > 0$ such that 
\begin{align*}
|\{\zeta > \delta_0\}|=|\{(z,r)\in \mathbb{R}^{2}_{+}: \zeta(z,r) > \delta_0\}| > 0.
\end{align*}
We choose a nonnegative function $h_* \in L^\infty_c(\bbR^2_+)$ supported in $\{\zeta > \delta_0\}$ and satisfying $P[h_*] = 1$. We set a projection of $h \in L^\infty_c(\bbR^2_+)$ onto a zero impulse space by  
\begin{equation}\label{eq:proj}
\mathbb{P}h = h - (r^2, h) h_*.
\end{equation}
If $\supp{h}\subset \{\zeta>\delta_0\}$, then $\supp{\mathbb{P}h}\subset \{\zeta>\delta_0\}$.

\noindent
(i) \textit{Existence of multipliers.}
We first show a weaker assertion
\begin{equation}\label{eq:EL_p_2}
    \begin{cases}
        \mathcal{G}\zeta - \kappa \zeta^{p-1} - \frac{W}{2} r^2 = 0&  \textrm{ on }\,\, \{\zeta > 0\}, \\
        \mathcal{G}\zeta - \frac{W}{2} r^2 \le 0&  \textrm{ on }\,\, \{\zeta = 0\},
    \end{cases}
\end{equation}
for some constants $W$, $\kappa$, assuming the existence of a nonnegative function $h_0 \in L^\infty_c(\bbR^2_+)$ satisfying 
\begin{equation}
\begin{aligned}
&\supp{h_0}\subset \{\zeta>\delta_0\},\\
&( \zeta^{p-1}, \mathbb{P}h_0)=( \zeta^{p-1}-C_{**}r^{2}, h_0) \neq 0,\quad C_{**}=(\zeta^{p-1},h_{*}).
\end{aligned}
\label{eq:nonzero}
\end{equation}
For such $h_0$, the projection $\eta_0=\mathbb{P}h_0$ also satisfies $\supp{\eta_0}\subset \{\zeta>\delta_0\}$.

We take arbitrary $\delta \in (0, \delta_0)$ and $h \in L^\infty_c(\bbR^2_+)$ satisfying $h \ge 0$ on $\{0 \le \zeta \le \delta\}$. For $\tau,\ s\in \mathbb{R}$ and $\eta =\mathbb{P}h$, we set 
\[
j(\tau, s) = 2\pi \int_{\bbR^2_+} |\zeta + \tau \eta + s \eta_0|^p r dzdr.
\]
Since $j(0,0)=1$ and $\partial_s j(0,0) = p(\zeta^{p-1}, \eta_0)\neq 0$, the implicit function theorem ensures the existence of a $C^1$ function $s(\tau)$ for small $|\tau|$ such that $s(0)=0$ and $j(\tau, s(\tau)) = 1$.
Differentiating this identity at $\tau=0$ yields
\[
s'(0) = - \frac{(\zeta^{p-1}, \eta )}{(\zeta^{p-1}, \eta_0 )}.
\]
The function 
\begin{align*}
\zeta(\tau) = \zeta+\tau \eta+s(\tau)\eta_0
\end{align*}
satisfies both $P[\zeta (\tau)]=1$ and $Q_p[\zeta(\tau)]=1$.
Moreover,  $\zeta(\tau)=\zeta+\tau h+s(\tau)\eta_0\ge 0$ on $\{0\le \zeta\le \delta\}$.
Thus $\zeta(\tau)$ is nonnegative for sufficiently small $\tau>0$. 
Since $\zeta(\tau)=\zeta+\tau(\eta+s'(0)\eta_0)+O(\tau^{2})$ as $\tau\to 0$, it follows that  
\begin{align*}
E[\zeta(\tau)]=E[\zeta]+2\tau E_{\textrm{int}}[\zeta,\eta+s'(0)\eta_0]+O(\tau^{2}).
\end{align*}
The maximality of $\zeta$ requires 
\[
0 \ge 2 E_{\textrm{int}}[\zeta,\eta+s'(0)\eta_0] =(\mathcal{G}\zeta,\eta+s'(0)\eta_0).
\]
Substituting the expression for $s'(0)$ into the preceding inequality, 
\begin{align}
(\mathcal{G}\zeta - \kappa \zeta^{p-1}, \eta)
\le 0,\quad 
\kappa = \frac{(\mathcal{G}\zeta,\eta_0)}{(\zeta^{p-1},\eta_0)}.  \label{eq:kappa}
\end{align}
By expanding $\eta=\mathbb{P}h$, 
\begin{align}
\left(\mathcal{G}\zeta - \kappa \zeta^{p-1}- \frac{W}{2} r^2, h\right)\le 0,\quad \frac{W}{2} =(\mathcal{G}\zeta - \kappa \zeta^{p-1}, h_*).  \label{eq:W}
\end{align}
Since $h \ge 0$ on $\{0 \le \zeta \le \delta\}$ and is arbitrary, taking $\delta \to 0$ leads to \eqref{eq:EL_p_2}.

\medskip

\noindent
(ii) \textit{Existence of an auxiliary function.}
We next show the existence of a nonnegative function $h_0\in L^{\infty}_{c}(\mathbb{R}^{2}_{+})$ satisfying \eqref{eq:nonzero}, so that \eqref{eq:EL_p_2} holds without additional conditions.

If the set $\{\zeta^{p-1}-C_{**}r^{2}\neq 0\}\cap \{\zeta>0\}$ has positive measure, we can construct the desired function. Indeed, for $A=\{\zeta^{p-1}-C_{**}r^{2}>0\}$ and $B=\{\zeta^{p-1}-C_{**}r^{2}<0\}$, at least one of the sets
$A\cap\{\zeta>0\}$ and
$B\cap\{\zeta>0\}$ has positive measure. We may assume that $|A\cap \{\zeta>0\}|>0$. We choose $\delta_0>0$ sufficiently small so that $|A\cap \{\zeta>\delta_0\}|>0$. Then, the function  
\begin{align*}
h_0=\mathbf{1}_{A\cap \{\zeta>\delta_0\}}
\end{align*}
satisfies the desired property \eqref{eq:nonzero}.

Thus it suffices to show that the set $\{\zeta^{p-1}-C_{**}r^{2}\neq 0\}\cap \{\zeta>0\}$ has a positive measure. Suppose, for contradiction, that $\{\zeta^{p-1}-C_{**}r^{2}\neq 0\}\cap \{\zeta>0\}=\emptyset$. Then, $\zeta^{p-1}=C_{**}r^{2}$ on $\{\zeta>0\}$. By $P[\zeta]=Q_p[\zeta]=1$, we have $C_{**}=1$. We take an arbitrary $\delta\in (0,\delta_0)$ and $h\in L^{\infty}_{c}(\mathbb{R}^{2}_{+})$ satisfying $h=0$ on $\{0\le \zeta\le \delta\}$ and set $\eta=\mathbb{P}h$. For arbitrary $\varepsilon>0$ and $\tau>0$, we set 
\begin{align*}
\zeta(\tau)=\zeta+\tau \eta-\tau \varepsilon h_{*}. 
\end{align*}
The function $\zeta(\tau)$ satisfies $P[\zeta(\tau)]=1-\varepsilon \tau<1$ and $\zeta(\tau)\ge 0$ for small $\tau>0$. By Taylor expansion,
\begin{align*}
\zeta(\tau)^{p}=\zeta^{p}+p\zeta^{p-1}(\eta-\varepsilon h_{*})\tau+O(\tau^{2})\quad \textrm{as}\ \tau\to 0.
\end{align*}
Using $(\zeta^{p-1},h_{*})=1$ and $(\zeta^{p-1},\eta)=(\zeta^{p-1}-r^{2},h)=0$, we find that 
\begin{align*}
Q_p^{p}[\zeta(\tau)]=2\pi \int_{\mathbb{R}^{2}_{+}}\zeta(\tau)^{p}rdzdr=1-p\varepsilon \tau+O(\tau^{2})<1.
\end{align*}
We thus differentiate the energy $E[\zeta(\tau)]$ at $\tau=0$. Repeating the computation above (replacing $s(\tau)$ by $-\tau$ and $\eta_0$ by $\varepsilon h_{*}$), we find that 
\begin{align*}
0\ge (\mathcal{G}\zeta,\eta-\varepsilon h_{*}).
\end{align*}
Letting $\varepsilon\to 0$ yields 
\begin{align*}
0\ge (\mathcal{G}\zeta,\eta)=(\mathcal{G}\zeta-C_{*}r^{2},h),\quad C_{*}=(\mathcal{G}\zeta,h_{*}).
\end{align*}
Since $h$ is arbitrary on $\{\zeta>\delta\}$ for each $\delta>0$, we deduce that $\mathcal{G}\zeta=C_{*}r^{2}$ a.e. on $\{\zeta>0\}$. Using a property of Sobolev functions, we differentiate $\mathcal{G}\zeta$ on $\{\zeta>0\}$ and deduce $\zeta=-r^{-2}L\mathcal{G}\zeta=0$. This is a contradiction. Thus, $|\{\zeta^{p-1}-C_{**}r^{2}\neq 0\}\cap \{\zeta>0\}|>0$.

\medskip

\noindent 
(iii) \textit{Signs of multipliers.} We now prove \eqref{eq:EL_intro0} by showing the signs of the constants $W \ge 0$ and $\kappa > 0$. 

By \eqref{eq:EL_p_2}, we have $\kappa \zeta^{p-1} \ge \mathcal{G}\zeta - \frac{W}{2} r^2$ in $\bbR^2_+$.
Using $\mathcal{G}\zeta > 0$, it follows that 
\begin{align*}
\kappa \zeta^{p-1} \ge \mathcal{G}\zeta - \frac{W}{2} r^2 > -\frac{W}{2} r^2.
\end{align*}
By $\zeta \in L^p_r$, $W$ is nonnegative. 

If $\kappa = 0$, we have $\mathcal{G}\zeta= \frac{W}{2}r^2$ a.e. on $\{\zeta > 0\}$ by \eqref{eq:EL_p_2}.
Differentiating $\mathcal{G}\zeta$, we find that $r^{2}\zeta = -L \mathcal{G}\zeta = 0$ a.e. on $\{\zeta > 0\}$. This is a contradiction.
Thus $\kappa \ne 0$. 

By \eqref{eq:nonzero}, we may choose $h \in \{ \eta_0, -\eta_0\}$   such that $(\zeta^{p-1}, h) < 0$.
Since $(r^2, h) = 0$, and $h \in L^{\infty}$ is supported in $\{ \zeta > \delta_0 \}$, we have $\zeta + \tau h \ge 0$ and $P[\zeta + \tau h] = 1$ for sufficiently small $\tau > 0$.
Moreover,
\[
Q_p^p[\zeta + \tau h] = 1 + p\tau(\zeta^{p-1}, h) + o(\tau) < 1.
\]
Since $\zeta \in S_p$, $(\mathcal{G}\zeta, h) \le 0$.
It follows from \eqref{eq:kappa} that $\kappa = (\mathcal{G}\zeta,h)/(\zeta^{p-1},h) \ge 0$.
Since $\kappa \ne 0$, we conclude that $\kappa > 0$, and \eqref{eq:EL_intro0} follows from \eqref{eq:EL_p_2}.

\medskip

\noindent 
(iv) \textit{Uniqueness of multipliers.} 
It remains to show the uniqueness of the constants $W$ and $\kappa$ satisfying \eqref{eq:EL_intro0}. We take another pair of constants $\overline{W} \ge 0$ and $\overline{\kappa} > 0$ satisfying \eqref{eq:EL_intro0}. Testing \eqref{eq:EL_intro0} against $\eta_0$ and using $(r^{2},\eta_0)=0$, we obtain $\overline{\kappa}={\kappa}$ from \eqref{eq:kappa}. Testing \eqref{eq:EL_intro0} against $h_*$ and using $(r^{2},h_*)=1$, we then obtain $\overline{W} = W$ from \eqref{eq:W}.
\end{proof}

\subsection{Poho\v{z}aev identity}
We establish the identities \eqref{eq:speed} and show that the Lagrange multipliers $W$ and $\kappa$ in Lemma \ref{lem:EL_p} are independent of the choice of maximizer in $S_p$.

\begin{proposition}\label{lem:reg_loc}
Let $1<p<\infty$ and $\zeta\in S_p$.
Let $\varphi$ be the function defined in \eqref{eq:5dtrans}.

\noindent
(i) We have $\zeta\in L^\infty(\bbR^5)$ and
$\varphi\in W^{2,q}_{\mathrm{loc}}(\bbR^5)$
for every $q\in(1,\infty)$.
In particular,
$\varphi\in C^{1,\beta}_{\mathrm{loc}}(\bbR^5)$
for every $\beta\in(0,1)$.

\noindent
(ii) Let $W\ge0$ and $\kappa>0$ be the constants
appearing in \eqref{eq:EL_intro0}.
Then $\varphi$ satisfies
\begin{equation}\label{eq:EL_5D_2}
-\Delta_y\varphi=\lambda\Psi_+^q,
\qquad
\Psi=r^2\left(\varphi-\frac{W}{2}\right),
\end{equation}
where
\begin{align}
q=\frac{1}{p-1},
\qquad
\lambda=\kappa^{-q}.
\end{align}
\end{proposition}

\begin{proof}
By $\kappa \zeta^{p-1} = \left(\mathcal{G}\zeta - \frac{W}{2} r^2 \right)_+ \le \mathcal{G}\zeta$ and \eqref{eq:str_2}, we have $\zeta \in L^\infty(\mathbb{R}^{2}_{+})$.
By the isometry \eqref{eq:ke2}, we have $\varphi\in \dot{H}^{1}_{\textrm{axi}}(\mathbb{R}^{5})$.
Applying the standard elliptic regularity theory to the Poisson equation \eqref{eq:5D} in $\mathbb{R}^{5}$, we deduce $\varphi \in W^{2,q}_{\mathrm{loc}}(\mathbb{R}^5)$ for any $1<q < \infty$, which embeds into $C^{1,\beta}_{\mathrm{loc}}(\mathbb{R}^5)$ for any $\beta \in (0, 1)$. The equation \eqref{eq:EL_5D_2} follows from \eqref{eq:EL_intro0} and \eqref{eq:5D}. 
\end{proof}

\begin{lemma}[Poho\v{z}aev identity]\label{lem:pohozaev}
The identities \eqref{eq:speed} hold. In particular, $W=W_p>0$.  
\end{lemma}

\begin{proof}
We apply the following identities for axisymmetric functions on $\mathbb{R}^5$:
\begin{align}
    \nabla \cdot \left( \nabla \varphi (y \cdot \nabla \varphi) - \frac{1}{2}y |\nabla \varphi|^2 \right) &= \Delta \varphi (y \cdot \nabla \varphi) - \frac{3}{2}|\nabla \varphi|^2,  \label{eq:0422_0} \\
    \frac{1}{r^{2}}y\cdot \nabla \varphi&=y\cdot \nabla \left(\frac{\varphi}{r^{2}}\right)+2\frac{\varphi}{r^{2}}.  \label{eq:0422_10}
\end{align}
Let $\zeta\in S_p$ and $\varphi=\mathcal{G}\zeta/r^{2}$. It follows from \eqref{eq:impulse5D}-\eqref{eq:ke2} that 
\begin{align}
I_p&=E[\zeta]=\frac{1}{2\pi}\int_{\mathbb{R}^5} |\nabla \varphi|^2 \, dy,  \label{eq:0909_5D} \\
1&=P[\zeta] =
    \frac{1}{\pi}\int_{\mathbb{R}^5} \zeta \, dy, \label{eq:5d_impulse1} \\
1&=Q_p^p[\zeta]
=\frac{1}{\pi}\int_{\mathbb R^5}
\frac{1}{|y'|^{2}}\zeta^p\,dy.
\label{eq:casimir5D1}
\end{align}
Observe that $\Psi=r^{2}(\varphi-W/2)$ satisfies $y \cdot \nabla \Psi = 2\Psi + r^2 (y \cdot \nabla \varphi)$ and  
\begin{align*}
    -\Delta \varphi (y \cdot \nabla \varphi) 
    = \lambda \Psi_+^q (y \cdot \nabla \varphi) 
    =\frac{\lambda \Psi_{+}^{q}}{r^{2}}(y \cdot \nabla \Psi-2\Psi)
=\frac{\lambda}{r^{2}}\left(y\cdot \nabla \frac{1}{q+1}\Psi_{+}^{q+1}-2\Psi_{+}^{q+1}\right).
\end{align*}
Applying the identity \eqref{eq:0422_10} yields
\begin{equation}
\begin{aligned}
    -\Delta \varphi (y \cdot \nabla \varphi) 
    &=\frac{\lambda}{q+1}y\cdot \nabla \left(\frac{\Psi_{+}^{q+1}}{r^{2}}\right)- \frac{2\lambda q}{q+1}\left(\frac{\Psi_{+}^{q+1}}{r^{2}}\right) \\
    &=\frac{\lambda}{q+1}
    \left( \nabla \cdot  \left(y\frac{\Psi_{+}^{q+1}}{r^{2}}\right)-(5+2q)\left(\frac{\Psi_{+}^{q+1}}{r^{2}}
    \right)\right).
\end{aligned}
\label{eq:0422_1}
\end{equation}
Using the identity \eqref{eq:0422_0}, we have
\begin{equation}
    \label{eq:poho2}
    \nabla \cdot \left(  \frac{\lambda}{q+1} \left(y\frac{\Psi_{+}^{q+1}}{r^{2}}\right) + \nabla \varphi (y \cdot \nabla \varphi) - \frac{1}{2}y |\nabla \varphi|^2 \right) = \lambda\frac{5+2q}{q+1} \left(\frac{\Psi_{+}^{q+1}}{r^2}\right) - \frac{3}{2}|\nabla \varphi|^2.
\end{equation}
We express $\Psi$, $q$, and $\lambda$ in terms of $\zeta$, $p$, and $\kappa$ using
\begin{equation*}
\Psi_{+}^{q+1}=\left(\frac{\zeta}{\lambda}\right)^{p},\quad \frac{5+2q}{q+1}=\frac{3}{\beta_p},\quad \kappa=\frac{1}{\lambda^{p-1}}
\end{equation*}
Integrating both sides of \eqref{eq:poho2} over $\mathbb{R}^5$, we obtain
\[
    \frac{3}{2}\int_{\mathbb{R}^5} |\nabla \varphi|^2 \, dy = \frac{3\kappa}{\beta_p}\int_{\mathbb{R}^5}\frac{1}{r^2}\zeta^{p}dy.
\]
Substituting \eqref{eq:casimir5D1} and \eqref{eq:0909_5D} implies the identity \eqref{eq:speed} for $\kappa=\kappa_p$. By multiplying $\kappa \zeta^{p-1}=r^{2}(\varphi-W/2)_+$ by $r^{-2}\zeta$ and integrating it over $\mathbb{R}^{5}$, we have
\begin{align*}
\kappa\int_{\mathbb{R}^5}\frac{1}{r^2}\zeta^{p}dy
=\int_{\mathbb{R}^5}\varphi \zeta dy-\frac{W}{2}\int_{\mathbb{R}^5}\zeta dy.
\end{align*}
It follows from \eqref{eq:0909_5D}--\eqref{eq:casimir5D1} that 
\begin{align*}
\kappa_p =  2I_p-\frac{W}{2}.
\end{align*}
Substituting $\kappa_p=\beta_pI_p$ and $\alpha_p+\beta_p=2$ into the above implies the identity \eqref{eq:speed} for $W=W_p$.
\end{proof}

\begin{remark}\label{r:poho}
The identities \eqref{eq:speed} can be written as   
\begin{align*}
\frac{\partial}{\partial \mu} I_{p,\mu,\nu} \Big|_{(\mu,\nu)=(1,1)} = \frac{W_p}{2}, \quad \frac{\partial}{\partial \nu} I_{p,\mu,\nu} \Big|_{(\mu,\nu)=(1,1)} = \kappa_p.
\end{align*}
\end{remark}

\begin{proposition}\label{p:uniformbound}
Let $1<p< \infty$. The following estimates hold: 
\begin{align}
\frac{1}{2}&<\alpha_p\le \frac{7}{5},\quad 
\frac{3}{5}\le \beta_p<\frac{3}{2}, \label{eq:alphabeta}\\
\left(\frac{1}{2\pi \Gamma[\zeta_H]}\right)^{\frac{3}{5p-3}}I_\infty &\le I_p=\frac{W_p}{2\alpha_p}=\frac{\kappa_p}{\beta_p}\le C\frac{p}{p-1}  \label{eq:unifrombound2},
\end{align}
where $\zeta_H$ is Hill's vortex \eqref{eq:normalizedHill0}. The constant $C$ is independent of $p$.
\end{proposition}

\begin{proof}
The estimate \eqref{eq:alphabeta} follows from \eqref{eq:exponent}.
The upper bound in \eqref{eq:unifrombound2} follows from \eqref{eq:en_bound}. Hill's vortex $\zeta_H \in S_{\infty}$ satisfies $P[\zeta_H]=1$ and $E[\zeta_H]=I_\infty$. Moreover, we have $Q_{p}^{p}[\zeta_H]=2\pi \Gamma[\zeta_H]$ by $\zeta_H=\mathbf{1}_{\Omega}$. We rescale $\zeta=\zeta_H$ by \eqref{eq:scalinglaw} and the constants 
\begin{align*}
a=b^{5},\quad b=\left(\frac{1}{2\pi \Gamma[\zeta_H]}\right)^{\frac{1}{5p-3}},
\end{align*}
so that the rescaled $\zeta_{a,b}$ satisfies 
\begin{align*}
E[\zeta_{a,b}]=b^{3}I_\infty,\quad P[\zeta_{a,b}]=1,\quad Q_{p}[\zeta_{a,b}]=1.
\end{align*}
The lower bound in \eqref{eq:unifrombound2} follows from the maximality of $I_p$.
\end{proof}

\subsection{Compactness of vortex core}\label{sec:max_pro}

The positivity of the traveling wave speed $W_p>0$ implies that energy maximizers have compact support.
Consequently, the circulation of every energy maximizer is finite.
We further use the positivity of $W_p$ to establish a uniform bound on the extent of the support in the $r$-direction.

\begin{lemma}[Compact support]\label{lem:compact}
Let $1<p< \infty$ and $\zeta \in S_p$.
Then, the support of $\zeta(z,r)$ is compact in $\overline{\bbR^2_+}$.
\end{lemma}

\begin{proof}
We show that $\textrm{supp}\ \zeta$ is compact in $\mathbb{R}^{5}$. By \eqref{eq:EL_5D_2}, we have 
\[
    \supp \zeta = \left\{ y \in \mathbb{R}^5 : \varphi(y) \ge \frac{W_p}{2} \right\}.
\]
By Proposition \ref{lem:reg_loc}, $\zeta\in L^{\infty}\cap L^{1}(\mathbb{R}^{5})$ and $\varphi=N*\zeta$ for the Newtonian potential $N(y) = (8\pi^{2})^{-1} |y|^{-3}\in L^{5/3,\infty}(\mathbb{R}^{5})$ by the uniqueness of the Poisson equation. By the H\"{o}lder inequality for Lorentz spaces, e.g., \cite[Theorem 1.4.16(v)]{Gra14}, $||\varphi||_{L^{\infty}(\mathbb{R}^{5})}\le C||N||_{L^{5/3,\infty}(\mathbb{R}^{5})}||\zeta||_{L^{5/2,1}(\mathbb{R}^{5})}$. We take a sequence $\{\zeta_n\}\subset C_c^{\infty}(\mathbb{R}^{5})$ such that $\zeta_n\le \zeta$, $\zeta_n\to \zeta$ in $L^{5/2,1}(\mathbb{R}^{5})$, and $\textrm{supp}\ \zeta_n\subset B_{L_n}(0)$.
For $|y|\ge L$,
\begin{align*}
N*\zeta_{n}(y)=\frac{1}{8\pi^{2}}\int_{B_{L_n}(0)}\frac{1}{|y-w|^{3}}\zeta_n(w)dw\le \frac{1}{8\pi^{2}(L-L_n)^{3}}||\zeta||_{L^{1}(\mathbb{R}^{5})}.
\end{align*}
We thus estimate $\varphi=N*\zeta_n+N*(\zeta-\zeta_n)$ by 
\begin{align*}
\sup_{|y|\ge L}\varphi(y)\le \frac{1}{8\pi^{2}(L-L_n)^{3}}||\zeta||_{L^{1}(\mathbb{R}^{5})}+C||N||_{L^{5/3,\infty}(\mathbb{R}^{5})}||\zeta-\zeta_n||_{L^{5/2,1}(\mathbb{R}^{5})}.
\end{align*}
Taking $n$ sufficiently large and then $L$ sufficiently large,
we see that the right-hand side is smaller than $W_p/2$. Thus, the support of $\zeta$ is contained within some ball in $\mathbb{R}^5$.
\end{proof}

\begin{lemma}[Uniform bounds]
\label{lem:unif_mass}
Let $1<p<\infty$. Let $W_p,\kappa_p>0$ be defined by \eqref{eq:speed}. Then, the following statements hold for every maximizer $\zeta\in S_p$:

\begin{align}
&\sup_{\zeta\in S_p}\|\zeta\|_{L^{\infty}(\mathbb{R}^{2}_{+})} \le M_p,\quad M_p=\left(\frac{pC}{\kappa_p(p-1)}\right)^{\frac{1}{p-1}},  \label{eq:vortexsupbound}\\
&\operatorname{supp}\zeta\subset \{(z,r)\in \overline{\mathbb{R}^{2}_{+}} :\ r\le \bar{R}_p \},\quad \bar{R}_p = \sqrt{\frac{2pC}{W_p(p-1)}},  \label{eq:rbounded}\\
&\frac{1}{2\pi \bar{R}_p^{2}}\le \Gamma[\zeta]\le \frac{8}{\pi}\left(\frac{\bar{R}_pM_p}{W_p}\right)^{2}.\label{eq:unif_mass}
\end{align}
The constant $C$ is independent of $p$. Moreover, for any $1<p_0<\infty$, the constants $M_p$ and $\bar{R}_p$ are uniformly bounded for $p_0\le p<\infty$.
\end{lemma}

\begin{proof}
The boundedness of the constants $M_p$ and $\bar{R}_p$ follows from \eqref{eq:unifrombound2}. The estimate \eqref{eq:vortexsupbound} follows from \eqref{eq:str_2} and \eqref{eq:EL_intro0}. We use the coordinates $X=(z,r)$. Since $\{\mathcal{G}\zeta-W_pX_2^{2}/2>0\}=\{\zeta > 0\}$, we have $\zeta(X) = 0$ for $X_2 > \bar{R}_p$. Thus \eqref{eq:rbounded} holds. This implies 
\[
1 = 2\pi\int_{\mathbb{R}^{2}_{+}}X_2^{3}\zeta\,dX \le 2\pi \bar{R}^2_p \int_{\mathbb{R}^{2}_{+}}X_2\zeta\,dX = 2\pi \bar{R}^2_p \Gamma[\zeta].
\]
Thus the lower bound in \eqref{eq:unif_mass} holds. We set
\[
h_{n}=2^{-n}\bar{R}_p, \quad H_n= \left\{X\in\mathbb{R}^{2}_{+}: \frac{h_{n}}{2}\le X_2<h_{n}\right\}.
\]
Since $\zeta$ vanishes for $X_2 > \bar{R}_p$,
\[
\Gamma[\zeta] = \sum_{n=0}^{\infty} \int_{H_n}\zeta(X)\,X_2 dX.
\]
On $H_n \cap \{\zeta>0\}$, we have 
\[
1 \le \frac{2\mathcal{G}\zeta(X)}{W_p X_2^{2}} \le \frac{8\mathcal{G}\zeta(X)}{W_p h_{n}^{2}}.
\]
It follows that 
\[
\Gamma[\zeta] \le \frac{8}{W_p} \sum_{n=0}^{\infty} h_{n}^{-2} \int_{H_n} \zeta(X) \mathcal{G}\zeta(X)\,X_2 dX.
\]
Using the Green function, we write   
\begin{align*}
\int_{H_n} \zeta(X) \mathcal{G}\zeta(X)\,X_2 dX = \int_{\mathbb{R}^{2}_{+}} \zeta(Y)Y_2 \left( \int_{H_n} G(X,Y)\zeta(X) \, X_2dX \right) dY.
\end{align*}
We apply \eqref{eq:green_strip} and \eqref{eq:vortexsupbound} to estimate 
\begin{align*}
\int_{H_n} G(X,Y)\zeta(X) \, X_2dX \le M_p \int_{0}^{h_{n}}\int_{\mathbb{R}} G(X,Y)\, X_2 dX \le \frac{M_p}{4}  Y_2 h_{n}^3.
\end{align*}
This yields 
\[
\int_{H_n} \zeta(X)\mathcal{G}\zeta(X) \,X_2 dX \le \frac{M_p}{4}h_{n}^3 \int_{\mathbb{R}^{2}_{+}}X_2^{2}\zeta\,dX.
\]
Using this upper bound, we estimate $\Gamma[\zeta]$. By using $\sum_{n=0}^{\infty}h_{n}=2\bar{R}_p$, it follows that 
\begin{align*}
\Gamma[\zeta]\le \frac{4\bar{R}_pM_p}{W_p}\int_{\mathbb{R}^{2}_{+}}X_2^{2}\zeta\,dX.
\end{align*}
By H\"older's inequality and $P[\zeta]=1$, we estimate 
\begin{align*}
\int_{\mathbb{R}^{2}_{+}}X_2^{2}\zeta\, dX   \le  \left( \int_{\mathbb{R}^{2}_{+}}X_2^{3}\zeta\,dX \right)^{1/2} \left( \int_{\mathbb{R}^{2}_{+}}X_2\zeta\,dX \right)^{1/2} = \frac{1}{\sqrt{2\pi}} \Gamma[\zeta]^{1/2}.
\end{align*}
Thus the upper bound in \eqref{eq:unif_mass} holds.
\end{proof}

\begin{corollary}\label{rmk:W_formula}
Let $1<p<\infty$. For any $\zeta\in S_p$, we have 
\begin{equation}\label{eq:W_formula_2}
    W_p = \frac{1}{8\pi \Gamma[\zeta]}\iint_{\bbR^2_+\times \bbR^2_+ } \frac{(X_2^2+Y_2^2)\bigl(F(s)-2sF'(s)\bigr)
-4(X_2-Y_2)^2F'(s)}{\sqrt{X_2Y_2}}
  \zeta(X) \zeta(Y) \,dX\,dY,
\end{equation}    
where $F(s)$ is the function in \eqref{eq:Green}.
\end{corollary}

\begin{proof}
Since the circulation is finite, we integrate $\zeta \partial_r \Psi$ in $\mathbb{R}_{+}^{2}$ for $\Psi = \mathcal{G}\zeta - W_pr^2/2$ and obtain   
\begin{equation*}
W_p = \frac{1}{\Gamma[\zeta]}\int_{\bbR^2_+} \zeta \partial_{X_2} \mathcal{G}\zeta \, dX.
\end{equation*}
This representation \eqref{eq:W_formula_2} follows from the expansion of $\mathcal{G}\zeta$ using the Green function. 
\end{proof}

We establish global regularity of maximizers using the compactness of their supports. For $\alpha \in (0,1)$ and a positive integer $m$, we denote by $C^{m,\alpha}(\mathbb{R}^5)$ the standard H\"older space of functions whose partial derivatives up to order $m$ are bounded and continuous, and the $m$-th order derivatives are $\alpha$-H\"older continuous. We use the ceiling function $\lceil s \rceil$ to denote the smallest integer greater than or equal to $s \in \bbR$.

\begin{lemma}[Global regularity]
Let $1<p<\infty$ and $\zeta\in S_p$. We have $\zeta\in C^{m,\alpha}(\overline{\bbR^2_+})$ for 
\begin{equation}\label{eq:m_alpha}
m=\left\lceil\frac{1}{p-1}-1\right\rceil,
\qquad
\alpha=\frac{1}{p-1}-m\in(0,1].
\end{equation}
\end{lemma}

\begin{proof}
We first consider the case $m=0$. In this case, $\varphi$ satisfies \eqref{eq:EL_5D_2} for $q=\alpha\in (0,1]$. By $\Psi \in C^{1,\beta}_{\mathrm{loc}}(\mathbb{R}^5)$ and $s_+^{q}\in C^{0,q}_{\textrm{loc}}(\mathbb{R})$, we have $\zeta = \lambda \Psi_+^{q}\in C^{0,q}_{\mathrm{loc}}(\mathbb{R}^5)$. Since $\zeta$ is compactly supported, $\zeta \in C^{0,q}(\mathbb{R}^5)$. By axisymmetry, $\zeta(z,r) \in C^{0, \alpha}(\overline{\bbR^2_+})$.

For the case $m=1$, $\zeta\in C^{0,\beta}(\mathbb{R}^{5})$ for any $\beta\in (0,1]$. By \eqref{eq:EL_5D_2} and Schauder estimates, we have $\varphi \in C^{2,\beta}_{\mathrm{loc}}(\mathbb{R}^5)$. Since $\Psi \in C^{2,\beta}_{\mathrm{loc}}(\mathbb{R}^5)$ and $s_+^{q}\in C^{1,\alpha}_{\textrm{loc}}(\mathbb{R})$, we obtain $\zeta \in C^{1,\alpha}_{\mathrm{loc}}(\mathbb{R}^5)$. Since $\zeta$ is compactly supported, $\zeta(z,r) \in C^{1, \alpha}(\overline{\bbR^2_+})$.

For $m\ge 1$, iterating the same argument yields
$\zeta\in C^{m,\alpha}(\overline{\bbR^2_+})$.
\end{proof}

\subsection{Steiner symmetry}

Using a polarization argument as in \cite{BT1976, BS2000}, we establish the Steiner symmetry of energy maximizers up to translation in the $z$-direction and derive a uniform bound on the extent of their supports in the $z$-direction.
For $1<p<\infty$, Steiner symmetry further implies that the boundary of the vortex core can be represented by a graph $z=f(r)$.
We also prove that, for $1<p<2$, the vortex core is simply connected by adapting the second-variation approach \cite[Theorem~3G]{FB74}.

\begin{lemma}[Steiner symmetry]\label{lem:steiner}
Let $1 < p < \infty$ and $\zeta\in S_p$.
There exists $c \in \bbR $ such that $\zeta(z + c,r)$ is even in $z$ and nonincreasing for $z > 0$.
\end{lemma}

\begin{proof}
For $\lambda \in \bbR$, let $\Sigma_{\lambda}= \{X\in \bbR^2_+ : X_1 > \lambda\}$ and $\sigma_{\lambda}X = (2\lambda - X_1, X_2)$. We set 
\[
\zeta^{\lambda}(X) =
\begin{cases}
\max\{\zeta(X), \,\, \zeta(\sigma_{\lambda}X) \}, & X \in \Sigma_{\lambda},\\
\min\{\zeta(X), \,\, \zeta(\sigma_{\lambda}X)\}, & X\notin \Sigma_{\lambda}.
\end{cases}
\]
Since $\zeta^{\lambda}$ is defined by interchanging the values of $\zeta$ at two points having the same $r$-coordinate, it preserves the impulse and $L^p$-norm.
Thus $\zeta^{\lambda} \in K_p$.

If $X, Y \in \Sigma_{\lambda}$, we see $|X_1 - Y_1| <  |X_1 - (2\lambda-Y_1)|$ and then $G(X,Y) > G(X,\sigma_{\lambda}Y)$. We set
\[
D_{\lambda}(X) = \zeta(X) - \zeta(\sigma_{\lambda}X), \quad D_{\lambda}^{+} = \max\{D_{\lambda},0\}, \quad D_{\lambda}^{-} = \max\{-D_{\lambda},0\}.
\]
Using symmetry, we have
\begin{align*}
E[\zeta^{\lambda}] - E[\zeta] =  \pi \iint_{\Sigma_{\lambda}\times \Sigma_{\lambda}} \left(G(X,Y)-G(X,\sigma_{\lambda}Y)\right) \left( D_{\lambda}^{+}(X)D_{\lambda}^{-}(Y) + D_{\lambda}^{-}(X)D_{\lambda}^{+}(Y) \right) X_2Y_2dXdY \ge 0.
\end{align*}
Since $\zeta \in S_p$  and $\zeta^{\lambda} \in K_p$, equality must hold.
Since $G(X, Y) - G(X, \sigma_{\lambda}Y) > 0$ on $\Sigma_{\lambda}\times \Sigma_{\lambda}$, the equality implies that, for every $\lambda \in \bbR$, either $D_{\lambda}\ge 0$ a.e. in $\Sigma_{\lambda}$ or $D_{\lambda}\le 0$ a.e. in $\Sigma_{\lambda}$. We choose $c \in \bbR$ such that $P[\zeta \bfone_{X_1 < c}]=P[\zeta \bfone_{X_1 > c}]$. It follows that
\begin{align*}
\int_{\Sigma_c}X_2^3D_c(X)\,dX &= \int_{ \{X_1 > c\} }X_2^3\zeta(X)\,dX - \int_{ \{X_1 < c\} }X_2^3\zeta(X)\,dX = 0.
\end{align*}
However, since $D_c$ has definite sign in $\Sigma_c$, it must be $D_c=0$ a.e. in $\Sigma_c$.
This implies  $\zeta(c + X_1,X_2) = \zeta(c - X_1,X_2)$ a.e. in $\bbR^2_+$.

To prove the nonincreasing property, we claim that $D_{\lambda}\le 0$ a.e. in  $\Sigma_{\lambda}$ for every $\lambda > c$.
Indeed, when $\lambda > c$, by the choice of $c$ we readily have
\[
\int_{ \{X_1 > \lambda\} }X_2^3\zeta(X)\,dX \le \int_{ \{X_1 > c\} }X_2^3\zeta(X)\,dX =  \int_{ \{X_1 < c\} }X_2^3\zeta(X)\,dX \le \int_{ \{X_1 < \lambda\} }X_2^3\zeta(X)\,dX.
\]
If $D_{\lambda}\ge 0$ a.e. in $\Sigma_{\lambda}$ for this $\lambda$, we then see
\begin{align*}
0 \le \int_{\Sigma_{\lambda}}X_2^3D_{\lambda}(X)\,dX = \int_{ \{X_1 > \lambda \}}X_2^3\zeta(X)\,dX - \int_{\{ X_1 < \lambda \} }X_2^3\zeta(X)\,dX \le 0.
\end{align*}
Thus $D_{\lambda} = 0$ a.e. in $\Sigma_{\lambda}$. This means that $\zeta$ is symmetric with respect to both $X_1 = c$ and $X_1 = \lambda$. For such $\zeta$, the nonnegative function $X_1 \longmapsto \int_{0}^{\infty}X_2^3\zeta(X)\,dX_2$ is periodic. This contradicts $P[\zeta]=1$. Therefore, $D_{\lambda}\le 0$ a.e. in $\Sigma_{\lambda}$ for every $\lambda > c$. This holds everywhere in $\Sigma_{\lambda}$ by continuity of $\zeta \in S_p$. For any $c < z_1 < z_2$, taking $\lambda = (z_1 + z_2)/2$ yields $\zeta(z_2,r) \le \zeta(z_1,r)$ for every $r > 0$.
Together with the symmetry about $z = c$, this proves the lemma.
\end{proof}

\begin{proposition}\label{p:p2_simply}
Let $1<p<\infty$ and $\zeta\in S_p$.
There exist $c\in\mathbb R$, a compact set $J\subset[0,\infty)$, and a nonnegative function $f:J\to[0,\infty)$ such that
\[
\operatorname{supp}\zeta = \{(z,r):r\in J,\ |z-c|\le f(r)\}.
\]
\end{proposition}

\begin{proof}
We show that the support of $\zeta$ is connected in the $z$-direction. Suppose that $\zeta=\zeta_1+\zeta_2$ for some functions $\zeta_1$ and $\zeta_2$ such that $\operatorname{supp}\zeta_1\subset \{X_1<a\}$ and $\operatorname{supp}\zeta_2\subset \{X_1>b\}$ for some constants $a<b$.
We set 
\begin{align*}
\tilde{\zeta}(X)=\zeta_1(X_1-t,X_2)+\zeta_2(X_1+t,X_2),\quad t>0.
\end{align*}
By the translational invariance of the impulse, the $L^p$-norm, and the self-interaction energies, we have  
\begin{align*}
P[\tilde{\zeta}]&=P[\zeta_1] + P[\zeta_2]=P[\zeta]=1,\\
Q_p^{p}[\tilde{\zeta}]&=Q_p^{p}[\zeta_1]+Q_p^{p}[\zeta_2]=Q_p^{p}[\zeta]=1,\\
E[\tilde{\zeta}]&=E[\zeta_1(\cdot - t,\cdot)+\zeta_2(\cdot+t,\cdot)]=E[\zeta_1]+E[\zeta_2]+2E_{\textrm{int}}[\zeta_1(\cdot-t,\cdot),\zeta_2(\cdot+t,\cdot)].
\end{align*}
In particular, $\tilde{\zeta}\in K_{p}$. Since the function $F(s)$ is strictly decreasing for $s>0$ by Lemma \ref{lem:p1_kernel}, the Green function $G(X,Y)=\sqrt{X_2Y_2}F(s)/2\pi$ satisfies 
\begin{align*}
G(X_1+t,X_2,Y_1-t,Y_2)>G(X,Y). 
\end{align*}
It follows that, for all sufficiently small $t>0$,
\begin{align*}
E_{\textrm{int}}[\zeta_1(\cdot-t,\cdot),\zeta_2(\cdot+t,\cdot)]
&=\pi \iint_{\mathbb{R}^{2}_{+}\times \mathbb{R}^{2}_{+}}G(X,Y) \zeta_1(X_1-t,X_2)\zeta_2(Y_1+t,Y_2)X_2Y_2 dXdY \\
&=\pi \iint_{\mathbb{R}^{2}_{+}\times \mathbb{R}^{2}_{+}}G(X_1+t,X_2,Y_1-t,Y_2) \zeta_1(X)\zeta_2(Y)X_2Y_2 dXdY\\
&>\pi \iint_{\mathbb{R}^{2}_{+}\times \mathbb{R}^{2}_{+}}G(X,Y) \zeta_1(X)\zeta_2(Y)X_2Y_2 dXdY=E_{\textrm{int}}[\zeta_1,\zeta_2].
\end{align*}
Hence $E[\tilde{\zeta}]>E[\zeta]$. This contradicts the maximality of $\zeta$. 
\end{proof}

\begin{lemma}\label{lem:p2_simply}
Let $1<p<2$ and $\zeta \in S_p$.
The support of $\zeta$ is simply connected.
\end{lemma}

\begin{proof}
We first show that $\mathcal{G}\zeta$ is a maximizer of the functional 
\[
J[\psi] = \frac{2\pi}{q+1}\int_{\bbR^2_+}\left(\psi-\frac{W_p}{2}X_2^2\right)_+^{q+1} X_2\,dX
\]
among all $\psi \in \dot H^1_0(\bbR^2_+;r^{-1})$ satisfying the constraint $E[\xi ] = I_p$ for $\xi=-X_2^{-2}L\psi$. Let $\alpha_p$ and $\beta_p$ be the constants in \eqref{eq:exponent}.
By concavity of $t^{\alpha_p/2}$, we have 
\[
2t^{\alpha_p/2}-\alpha_p t
\le 2-\alpha_p,
\qquad t\ge0.
\]
Let $\eta$ be an arbitrary nonnegative function supported on $\{\psi-W_pX_2^{2}/2>0\}$.
We may assume $Q_p[\eta] = 1$ by normalization.
Then, it follows from $\alpha_p+\beta_p=2$ and the identity \eqref{eq:speed} that 
\begin{align*}
\left(\psi-\frac{W_p}{2}X_2^2,\eta\right) &= 2E_{\textrm{int}}[\xi,\eta]-\frac{W_p}{2}P[\eta] 
\le 2\sqrt{E[\xi]}\sqrt{E[\eta]} -\frac{W_p}{2}P[\eta] \\
&\le 2I_p P[\eta]^{\alpha_p/2}-\alpha_pI_pP[\eta] 
\le (2-\alpha_p)I_p =  \kappa_p.
\end{align*}
By duality $1/p+1/(q+1)=1$, we obtain 
\[
\left\|\left(\psi - \frac{W_p}{2}X_2^2\right)_+ \right\|_{L^{q+1}(\mathbb{R}^{3})} \le \kappa_p.
\]
The maximizer $\zeta$ satisfies $\kappa_p \zeta^{p-1}=(\mathcal{G}\zeta-\frac{W_p}{2}X_2^{2})_{+}$ and $Q_p[\zeta] = 1$. Multiplying this equation by $\zeta$ and integrating it over $\mathbb{R}^{2}_{+}$, we have
\[
\left\|\left(\mathcal{G}\zeta-\frac{W_p}{2}X_2^{2}\right)_{+}\right\|_{L^{q+1}(\mathbb{R}^{3})} = \kappa_p.
\]
Thus equality holds for $\psi=\mathcal{G}\zeta$, and hence $\mathcal{G}\zeta$ maximizes $J$ under the above constraint.

We take a variation of $J$ around $\mathcal{G}\zeta$. We consider the half-plane formulation of the elliptic problem \eqref{eq:EL_5D_2}:
\begin{align*}
-\frac{1}{X_2^{2}}L\Psi=\lambda \Psi_{+}^{q},\quad \Psi=\mathcal{G}\zeta-\frac{W_p}{2}X_2^{2}.
\end{align*}
Suppose that $\{\zeta> 0\}=\{\Psi> 0\}$ has two distinct connected components $\Omega_1$ and $\Omega_2$. We set 
\[
\Psi_j = \Psi \mathbf1_{\Omega_j}, \quad j = 1,2.
\]
Testing the equation for $\Psi$ against $\Psi_j$, we have 
\[
\int_{\bbR^2_+}\nabla\mathcal{G}\zeta \cdot \nabla\Psi_j\frac{1}{X_2}\,dX = \int_{\bbR^2_+}|\nabla\Psi_j|^2\frac{1}{X_2}\,dX = \lambda \int_{\Omega_j}\Psi^{q+1}X_2\,dX.
\]
We set $\Phi= c_1\Psi_1 - c_2\Psi_2$ and choose the constants $c_1, c_2 > 0$ so that 
\[
\int_{\bbR^2_+}\nabla\mathcal{G}\zeta\cdot\nabla \Phi \frac{1}{X_2}\,dX = 0.
\]
By the above identity, we have  
\[
\pi\int_{\bbR^2_+}|\nabla \Phi|^2\frac{1}{X_2}\,dX
=\lambda \pi\sum_{j=1}^{2}c_j^{2}\int_{\Omega_j}\Psi^{q+1}X_2\,dX=:A. 
\]
We consider a variation of $J$ around $\mathcal{G}\zeta$ using the function  
\[
\psi(t) = \mathcal{G}\zeta\cos t+\left(\frac{I_p}{A}\right)^{1/2}\Phi\sin t.
\]
By the orthogonal condition of $\psi$ and $\Phi$, the function $\psi(t)$ satisfies the constraint 
\[
\pi\int_{\bbR^2_+}|\nabla \psi(t)|^2\frac{1}{X_2}\,dX = I_p.
\]
By the condition $1<p<2$, we have $q>1$. The second variation of $J$ at $t=0$ is well-defined for such $q$, and we have 
\[
\frac{d^2}{dt^2}J[\psi(t)]\Bigg|_{t=0} =
\frac{qI_p}{A}\int_{\mathbb{R}^{2}_{+}} \Psi_+^{q-1}\Phi^2-\int_{\mathbb{R}^{2}_{+}}\Psi^{q}_{+}\mathcal{G}\zeta,
\]
where the measure $2\pi X_2dX$ is suppressed. Expanding $\Phi=c_1\Psi_1-c_2\Psi_2$, the integral of the first term is expressed as 
\[
\int_{\bbR^2_+}\Psi_+^{q-1}\Phi^2 
=\sum_{j=1}^{2}c_j^{2}\int_{\Omega_j}\Psi_+^{q+1}=\frac{2A}{\lambda}.
\]
Using $\zeta=\lambda\Psi_+^{q}$ and $P[\zeta]=1$, the second term is expressed as 
\[
\int_{\mathbb{R}^{2}_{+}}\Psi^{q}_{+}\mathcal{G}\zeta=\frac{1}{\lambda}\int_{\mathbb{R}^{2}_{+}}\zeta\mathcal{G}\zeta=\frac{2E[\zeta]}{\lambda}=\frac{2I_p}{\lambda}.
\]
It follows that
\begin{align*}
\frac{d^2}{dt^2}J[\psi(t)]\Bigg|_{t=0} =\frac{2I_p}{\lambda}(q-1)>0.
\end{align*}
This contradicts the maximality of $\mathcal{G}\zeta$.
Therefore, $\{\Psi > 0\}$ is connected. After a translation in the $z$-direction, every section of $\{\zeta>0\}$ at fixed $r$ is an interval centered at $z=0$ by Lemma \ref{lem:steiner}. Thus $\{\Psi > 0\}$ is simply connected.
\end{proof}

\begin{lemma}[Vortex core bound]
\label{lem:unif_supp}
Let $1<p<\infty$. Let $W_p$, $M_p$, and $\bar{R}_p$ be the constants appearing in \eqref{eq:speed}, \eqref{eq:vortexsupbound}, and \eqref{eq:rbounded}, respectively. Let $C_1,C_2>0$ be suitable absolute constants independent of $p$, and set
\begin{align}
\bar{Z}_p=C_1a_p\max\left\{
1,\frac{M_p^{3/2}}{W_p^{5/2}}
\right\},
\qquad
a_p=\left(\frac{4C_2}{W_p}\right)^{1/3}.
\label{eq:zbounded}
\end{align} 
Then, for every $\zeta\in S_p$, there exists $c\in\mathbb{R}$ such that 
\begin{equation}
\label{eq:unif_supp}
\operatorname{supp}\zeta(\cdot+c,\cdot) \subset \left\{(z,r)\in\overline{\mathbb{R}^{2}_{+}}: |z| \le \bar{Z}_p, \,\, r\le \bar{R}_p\right\}.
\end{equation}
Moreover, for any $1<p_0<\infty$, the constants $\bar{R}_p$ and $\bar{Z}_p$ are uniformly bounded for $p_0\le p<\infty$.
\end{lemma}

\begin{proof}
The uniform boundedness of $\bar{R}_p$ and $\bar{Z}_p$ follows from \eqref{eq:unifrombound2}-\eqref{eq:rbounded}.

It suffices to show a horizontal bound for $\operatorname{supp}\zeta$. By Lemma~\ref{lem:steiner}, we may assume that $\zeta$ is even in $z$ and nonincreasing for $z>0$. In particular, the function 
\[
g(z)=\int_{\mathbb{R}^4}\zeta(z,|w'|)\,dw'
\]
is even and nonincreasing for $z>0$. By $P[\zeta]=1$, we have $||g||_{L^{1}(\mathbb{R})}=\pi$. For $a>0$, we take $z>2a$ so that $z - a > z/2$. It follows that 
\begin{align*}
\int_{\{|w_1-z|<a\}}\zeta(w)\,dw = \int_{z-a}^{z+a}g(s)\,ds \le 2a g(z-a) \le \frac{2a}{z-a}\int_{0}^{z-a}g(s)\,ds \le \frac{4\pi a}{z}.
\end{align*}
By a standard potential estimate, for $f\in L^1\cap L^\infty(\mathbb{R}^5)$, we have 
\[
\int_{\mathbb{R}^5}\frac{f(w)}{|w-y|^3}\,dw \lesssim  \|f\|_{L^1(\bbR^5)}^{2/5}\|f\|_{L^\infty(\bbR^5)}^{3/5}.
\]
Applying this estimate for $f=\mathbf{1}_{\{|w_1-z|<a\}}\zeta$ and \eqref{eq:vortexsupbound}, we estimate
\begin{align*}
\varphi(y)=N*\zeta(y) 
&= \frac{1}{8\pi^{2}}\int_{\{|w_1-z|
\ge a\}}\frac{\zeta(w)}{|w-y|^3}\,dw + \frac{1}{8\pi^{2}}\int_{\{|w_1-z|<a\}}\frac{\zeta(w)}{|w-y|^3}\,dw  \\
&\le C\left(\frac{1}{a^{3}} + \left(\frac{a}{z}\right)^{2/5}M_p^{3/5}\right).
\end{align*}
We choose $a=a_p$ as in \eqref{eq:zbounded}
with $C_2=C$, so that the first term equals $W_p/4$.
We then choose $\bar{Z}_p$ as in \eqref{eq:zbounded}
with $C_1=\max\{(4C)^{5/2},2\}$,
so that the second term is strictly less than
$W_p/4$ for all $z>\bar{Z}_p$. Since $\{\zeta>0\}=\{\varphi-W_p/2>0\}$, this yields $\zeta(z,r) = 0$ for $z > \bar{Z}_p$. By the evenness of $\zeta$, the same holds for $z<-\bar{Z}_p$. Thus \eqref{eq:unif_supp} holds.
\end{proof}

\section{Existence and variational compactness}\label{s:4}

In this section, we establish the compactness of maximizing sequences for the variational problem \eqref{eq:variational} for $1 < p \le\infty$ using the concentration-compactness principle.
In particular, this yields the existence of energy maximizers for $1 < p < \infty$.
Combining this existence result with the properties of energy maximizers established in \S\ref{s:3}, we complete the proof of Theorem \ref{t:MT1}.

For $p=\infty$, energy maximizers are given by Hill's vortex. Nevertheless, we also establish the compactness of maximizing sequences in this case, as it will be used to prove the stability result for Hill's vortex in Theorem \ref{t:MT3}(iv). We apply the following concentration-compactness principle \cite{Lions1984}; see also \cite{BNL13,AC2022}.

\begin{lemma}[Concentration-compactness]\label{lem:lions}
Let $0<\mu<\infty$. Let $\{\rho_n\}_{n \ge 1} \subset L^1(\bbR^2_+)$ be a sequence of nonnegative functions satisfying
\[
\int_{\bbR^2_+} \rho_n \, dX = \mu_n \to \mu \quad \text{as } n \to \infty.
\]
Then, there exists a subsequence $\{\rho_{n_k}\}$ satisfying one of the following three conditions:

\noindent (i) \textit{Compactness}: There exists a sequence
$\{X_k\}\subset\overline{\bbR^2_+}$
such that, for every $\varepsilon>0$,
there exists $L>0$ satisfying
\[
\liminf_{k \to \infty} \int_{B_L(X_k) \cap \bbR^2_+} \rho_{n_k} \, dX \ge \mu - \varepsilon.
\]

\noindent (ii) \textit{Vanishing}: For any $L > 0$,
\[
\lim_{k \to \infty} \left( \sup_{X \in \bbR^2_+} \int_{B_L(X) \cap \bbR^2_+} \rho_{n_k} \, dX \right) = 0.
\]

\noindent (iii) \textit{Dichotomy}: There exists $\bar{\mu} \in (0, \mu)$ such that for any $\varepsilon > 0$, there exist $k_0 \ge 1$ and sequences $\{\rho_{k, 1}\}, \{\rho_{k, 2}\} \subset L^1(\bbR^2_+)$ satisfying $0 \le \rho_{k, 1}, \rho_{k, 2} \le \rho_{n_k}$ for $k \ge k_0$, $\supp\rho_{k, 1} \cap \supp\rho_{k, 2} = \emptyset$, $\dist(\supp\rho_{k, 1}, \supp\rho_{k, 2}) \to \infty$, and
\[
\limsup_{k \to \infty} \left( \|\rho_{n_k} - (\rho_{k, 1} + \rho_{k, 2})\|_{L^1(\bbR^2_+)} + \left| \|\rho_{k, 1}\|_{L^1(\bbR^2_+)} - \bar{\mu} \right| + \left| \|\rho_{k, 2}\|_{L^1(\bbR^2_+)} - (\mu - \bar{\mu}) \right| \right) \le \varepsilon.
\]
\end{lemma}

\begin{theorem}\label{thm:existence}
Let $1<p\le \infty$ and $\mu, \nu > 0$.
For any sequence of nonnegative functions $\{\zeta_n\}_{n \ge 1}$ on $\bbR^2_+$ satisfying
\begin{equation}\label{eq:max_seq}
\liminf_{n \to \infty} E[\zeta_n] \ge I_{p,\mu, \nu}, \quad \limsup_{n \to \infty} P[\zeta_n] \le \mu, \quad \limsup_{n \to \infty} Q_p[\zeta_n] \le \nu,
\end{equation}
there exist $\zeta\in S_{p,\mu,\nu}$
and a subsequence of $\{\zeta_n\}$
which, after suitable translations in the $z$-direction, satisfy
\begin{align}
&\hspace{10pt}\zeta_n \to \zeta\hspace{11pt} \textrm{in} \hspace{7pt} L^1_*(\bbR^2_+; r),  \label{eq:conv1}\\
&\begin{cases}
\zeta_n \to \zeta \,\,\hspace{4pt} \text{ in }\,\, L^p(\bbR^2_+; r)  &\text{ for }\,\, 1 < p < \infty,\\
\zeta_n \rightharpoonup^* \zeta \,\, \text{ in }\,\, L^\infty(\bbR^2_+; r)  &\text{ for }\,\, p = \infty.
\end{cases}
\label{eq:conv2}
\end{align}
\end{theorem}

\begin{proof}
We reduce to the case $\mu =\nu = 1$ by the scaling \eqref{eq:scaling}. The sequence satisfying \eqref{eq:max_seq} must obey 
\[
\lim_{n \to \infty} E[\zeta_n] = I_{p,1,1}, \quad \lim_{n \to \infty} P[\zeta_n] = 1, \quad \lim_{n \to \infty} Q_{p}[\zeta_n] = 1,
\]
by Remark \ref{sec:saturation}. We apply Lemma~\ref{lem:lions} with $\rho_n = 2\pi X_2^3 \zeta_n$ and obtain a subsequence (still denoted by $\{\rho_n\}$) that falls into one of three scenarios: vanishing, dichotomy, or compactness.
We proceed by ruling out the first two possibilities and then deduce the compactness of the maximizing sequence.\\

\noindent \textit{(ii) Vanishing}.
Suppose that vanishing occurs.
For any $L > 0$,
\begin{equation*}
\lim_{n \to \infty} \left( \sup_{X \in \bbR^2_+} \int_{B_L(X) \cap \bbR^2_+} Y_2^3 \zeta_n(Y) \, dY \right) = 0.
\end{equation*}
We show that this implies $\lim_{n \to \infty} E[\zeta_n] = 0$, which contradicts $I_{p,1,1}> 0$. For $L > 1$, we decompose the energy into far-field and near-field parts:
\begin{equation*}
E[\zeta_n] = \pi \left[ \iint_{\{|X - Y| \ge L\}} + \iint_{\{|X - Y| < L\}} \right] G(X, Y) \zeta_n(X) \zeta_n(Y) \, X_2 Y_2 dX dY.
\end{equation*}
For the far-field term, we apply the bound \eqref{eq:Green_bound} with $\tau = 3/2$, which gives $G(X, Y) \lesssim X_2^2 Y_2^2/|X-Y|^3$.
This yields
\begin{equation*}
\iint_{\{|X - Y| \ge L\}} G(X, Y) \zeta_n(X) \zeta_n(Y) \, X_2 Y_2 dX dY \lesssim \frac{1}{L^3}  P[\zeta_n]^2 \lesssim \frac{1}{L^3}.
\end{equation*}
For the near-field term, we further split the integral as
\begin{equation*}
\pi \iint_{\{|X - Y| < L\}} G(X, Y) \zeta_n(X) \zeta_n(Y) \, X_2 Y_2 dX dY = I_n(L) + II_n(L),
\end{equation*}
where $I_n(L)$ and $II_n(L)$ are the integrals on $\{|X-Y|<L\}$ restricted to the regions $\{G \ge L X_2^2 Y_2^2\}$ and $\{G < L X_2^2 Y_2^2\}$, respectively. The condition $G \ge L X_2^2 Y_2^2$ and \eqref{eq:Green_bound} with $\tau = 3/2$ imply $|X-Y| \lesssim L^{-1/3}$.
Applying \eqref{eq:Green_bound} with $\tau \in (0, 3/2)$ yields
\[
G(X, Y) X_2 Y_2 \lesssim \frac{(X_2Y_2)^{3/2 + \tau}}{ |X-Y|^{2\tau}}.
\]
We thus estimate
\begin{equation*}
I_n(L) \lesssim_\tau \iint_{\{|X - Y| \lesssim L^{-1/3}\}} \frac{X_2^{3/2 + \tau} \zeta_n(X) \cdot Y_2^{3/2 + \tau} \zeta_n(Y)}{|X - Y|^{2\tau}} \, dX dY.
\end{equation*}
For sufficiently small $\tau > 0$, we can choose $\theta \in (0, 1)$ such that $3\theta + (1-\theta)/p = 3/2 + \tau$. For $q > 1$ satisfying $1/q = \theta + (1-\theta)/p$, we have $3/2 + \tau - 1/q = 2\theta$. With this choice, \eqref{eq:IE} ensures that $\|X_2^{2\theta}\zeta_n\|_{L^q_r}$ is uniformly bounded.
Applying Young's convolution inequality, we get
\begin{equation*}
I_n(L) \lesssim_\tau \left( \int_{\{|X| \lesssim L^{-1/3}\}} \frac{1}{|X|^{\frac{\tau q}{q-1}}} \, dX \right)^{\frac{2q-2}{q}} \|X_2^{2\theta}\zeta_n\|_{L^q_r}^2 \lesssim L^{-\frac{4(1-1/q) - 2\tau}{3}},
\end{equation*}
where we choose $\tau > 0$ small enough such that $\tau < 2(1-1/q)$.
Thus, $\sup_n I_n(L) \to 0$ as $L \to \infty$.
Using $G < L X_2^2 Y_2^2$, we estimate
\begin{align*}
II_n(L) \lesssim L \iint_{\{|X - Y| < L\}} X_2^3 \zeta_n(X) Y_2^3 \zeta_n(Y) \, dX dY 
\lesssim  L \left( \sup_{X \in \bbR^2_+} \int_{B_L(X) \cap \bbR^2_+} Y_2^3 \zeta_n(Y) \, dY \right).
\end{align*}
By the vanishing assumption, we have $\lim_{n \to \infty} II_n(L) = 0$ for any fixed $L > 0$.
Hence, letting $n \to \infty$ and then $L \to \infty$, we conclude $\lim_{n \to \infty} E[\zeta_n] = 0$.
Thus, vanishing cannot occur.\\

\noindent \textit{(iii) Dichotomy}.
Suppose that dichotomy occurs for some $\bar{\mu} \in (0, 1)$.
For any $\varepsilon > 0$, we decompose $\zeta_n$ by $\zeta_n = \zeta_{1,n} + \zeta_{2,n} + \tilde{\zeta}_n$ such that $0 \le \zeta_{1,n}, \zeta_{2,n} \le \zeta_n$, $d_n = \dist(\supp\zeta_{1,n}, \supp\zeta_{2,n}) \to \infty$, and the remainder term $\tilde{\zeta}_n$ satisfies $\limsup_{n \to \infty} P[\tilde{\zeta}_n] \le \varepsilon$.
We denote $\mu_{i,n} = P[\zeta_{i,n}]$ and $\nu_{i,n} = Q_p[\zeta_{i,n}]$ for $i=1, 2$.
The total energy expands as
\[
E[\zeta_n] = E[\zeta_{1,n}] + E[\zeta_{2,n}] + I_n + II_n,
\]
where $I_n = 2E_{\mathrm{int}}[\zeta_{1,n}, \zeta_{2,n}]$ and $II_n$ contains all cross-interactions involving the remainder $\tilde{\zeta}_n$.
Since $d_n \to \infty$, we apply \eqref{eq:Green_bound} with $\tau = 3/2$ to estimate 
\begin{equation*}
I_n = 2\pi \iint G(X, Y) \zeta_{1,n}(X) \zeta_{2,n}(Y) \, X_2 Y_2 dX dY \lesssim \frac{1}{d_n^3} \iint X_2^3 \zeta_{1,n}(X) Y_2^3 \zeta_{2,n}(Y) \, dX dY \lesssim \frac{1}{d_n^3}.
\end{equation*}
Applying \eqref{eq:en_cross} along with the uniform $L^p_r$-bounds yields $II_n \lesssim P[\tilde{\zeta}_n]^{\frac{7p-6}{10p - 6}}$.
Since $E[\zeta_{i,n}] \le I_{p,\mu_{i,n}, \nu_{i,n}}$, we arrive at:
\begin{equation*}
E[\zeta_n] \le I_{p,\mu_{1,n}, \nu_{1,n}} + I_{p,\mu_{2,n}, \nu_{2,n}} + \frac{C}{d_n^3} + C' P[\tilde{\zeta}_n]^{\frac{7p-6}{10p-6}}.
\end{equation*}
Taking the limsup as $n \to \infty$ followed by $\varepsilon \to 0$, we extract subsequential limits $\mu_i$ and $\nu_i$ from $\mu_{i,n}$ and $\nu_{i,n}$. The dichotomy assumption ensures that $\mu_1 = \bar{\mu}$ and $\mu_2 = 1-\bar{\mu}$, yielding $\mu_1 + \mu_2 = 1$. This gives
\[
I_{p,1,1} \le I_{p,\mu_1, \nu_1} + I_{p,\mu_2, \nu_2},
\]
which contradicts strict superadditivity \eqref{eq:superadd_p} and \eqref{eq:superadd_inf}.
Therefore, dichotomy does not occur.\\

\noindent \textit{(i) Compactness}.
There exists a sequence $\{Y_n = (z_n, r_n)\} \subset \overline{\bbR^2_+}$ such that for any $\varepsilon > 0$, there exists $L > 0$ satisfying
\begin{equation}\label{eq:tight}
\liminf_{n \to \infty} \left(2\pi \int_{B_L(Y_n) \cap \bbR^2_+} X_2^3 \zeta_n \, dX\right) \ge 1 - \varepsilon, \quad \limsup_{n\to\infty} \left(2\pi \int_{\bbR^2_+ \setminus B_L(Y_n)} X_2^3 \zeta_n \, dX\right) \le \varepsilon.
\end{equation}
First, we rule out radial escape.
Suppose $\limsup_{n \to \infty} r_n = \infty$, so we may assume $r_n \ge 2L$ and $r_n \to \infty$.
The energy decomposes as
\begin{equation*}
E[\zeta_n] = \pi \int_{B_L(Y_n) \cap \bbR^2_+} \zeta_n \mathcal{G}\zeta_n  \, X_2 dX + \pi\int_{\bbR^2_+ \setminus B_L(Y_n)} \zeta_n\mathcal{G}\zeta_n  \, X_2 dX.
\end{equation*}
Using $X_2 \ge r_n - L$ on $B_L(Y_n)$ and \eqref{eq:str_2}, the first term is bounded by
\begin{equation*}
\int_{B_L(Y_n) \cap \bbR^2_+}  \zeta_n\mathcal{G}\zeta_n \, X_2 dX \le \|\mathcal{G}\zeta_n\|_{L^{\infty}(\mathbb{R}^{2}_{+})} \int_{B_L(Y_n) \cap \bbR^2_+} \frac{1}{X_2^2} (X_2^3 \zeta_n) \, dX \lesssim \frac{1}{(r_n - L)^2}.
\end{equation*}
For the second term, applying \eqref{eq:en_cross} yields
\begin{equation*}
\pi\int_{\bbR^2_+ \setminus B_L(Y_n)} \zeta_n\mathcal{G}\zeta_n \, X_2 dX = E_{\mathrm{int}}[\zeta_n, \zeta_n \mathbf{1}_{\bbR^2_+ \setminus B_L(Y_n)}] \lesssim P[\zeta_n \mathbf{1}_{\bbR^2_+ \setminus B_L(Y_n)}]^{\frac{7p-6}{10p-6}} \lesssim \varepsilon^{\frac{7p-6}{10p - 6}}.
\end{equation*}
Taking $n \to \infty$ followed by $\varepsilon \to 0$ yields $\lim_{n \to \infty} E[\zeta_n] = 0$, which contradicts $I_{p,1,1} > 0$.
Hence, $\{r_n\}$ is uniformly bounded.

To conclude compactness, we establish the convergence of a suitably translated sequence to a maximizer.
Since $\{r_n\}$ is bounded, by enlarging $L$ and translating the sequence along the $z$-axis, we may assume that the tightness condition \eqref{eq:tight} holds with $Y_n=(0,0)$.
Since $Q_p[\zeta_n] \to 1$, we extract a subsequence satisfying 
\begin{equation*}
\left\{
\begin{aligned}
&\zeta_n\rightharpoonup \zeta \quad \textrm{in}\ L^{p}(\mathbb{R}^{2}_{+};r)\quad \textrm{for}\ 1<p<\infty,\\
&\zeta_n\rightharpoonup^* \zeta \quad \textrm{in}\ L^{\infty}(\mathbb{R}^{2}_{+};r)\quad \textrm{for}\ p=\infty,
\end{aligned}
\right.
\end{equation*}
for some $\zeta$.
By the lower semi-continuity of the norms, we have $Q_p[\zeta] \le 1$ and $P[\zeta] \le  1$.
Thus $\zeta \in K_{p,1,1}$.

To show that $\zeta$ is a maximizer, we establish $\lim_{n \to \infty} E[\zeta_n] = E[\zeta]$. Using $B_L = B_L(0) \cap \bbR^2_+$, we decompose the energy into two regions as
\begin{equation*}
E[\zeta_n] =  \pi \int_{B_L}  \zeta_n\mathcal{G}\zeta_n \, X_2 dX + \pi \int_{\bbR^2_+ \setminus B_L} \zeta_n\mathcal{G}\zeta_n \, X_2 dX.
\end{equation*}
The first term can be written as 
\begin{align*}
\int_{B_L}  \zeta_n\mathcal{G}\zeta_n \, X_2 dX &= \int_{B_L} \int_{B_L} G(X,Y) \zeta_n(X) \zeta_n(Y) \, X_2 Y_2 dY dX \\
&\quad +  \int_{B_L} \int_{\bbR^2_+ \setminus B_L} G(X, Y) \zeta_n(X) \zeta_n(Y) \, X_2 Y_2 dY dX.
\end{align*}
By the symmetry $G(X, Y) = G(Y, X)$, the second term is bounded by
\begin{equation*}
 \int_{B_L} \int_{\bbR^2_+ \setminus B_L} G(X, Y) \zeta_n(X) \zeta_n(Y) \, X_2 Y_2 dX dY \le \int_{\bbR^2_+ \setminus B_L} \zeta_n \mathcal{G}\zeta_n \, X_2 dX.
\end{equation*}
Applying \eqref{eq:en_cross} along with the tightness \eqref{eq:tight} yields
\begin{align*}
0 &\le E[\zeta_n] - \pi \int_{B_L} \int_{B_L} G(X, Y) \zeta_n(X) \zeta_n(Y) \, X_2 Y_2 dY dX \\
&\le 2 \pi \int_{\bbR^2_+ \setminus B_L}  \zeta_n\mathcal{G}\zeta_n \, X_2 dX \lesssim P[\zeta_n \mathbf{1}_{\bbR^2_+ \setminus B_L}]^{\frac{7p-6}{10p - 6}} \lesssim \varepsilon^{\frac{7p-6}{10p-6}}.
\end{align*}
We apply an identical decomposition to $E[\zeta]$. Combining those estimates, we obtain
\begin{align*}
|E[\zeta_n] - E[\zeta]| \lesssim  \left| \iint_{B_L \times B_L}  G(X, Y) \big( X_2\zeta_n(X) Y_2\zeta_n(Y) - X_2\zeta(X) Y_2\zeta(Y) \big) \, dY dX \right| 
 + \varepsilon^{\frac{7p-6}{10p-6}}.
\end{align*}
For the conjugate exponent $p' \in [1, \infty)$, we have $G \in L^{p'}(B_L \times B_L)$ by \eqref{eq:Green_bound} with sufficiently small $\tau < 1/p'$.
The sequence $X_2\zeta_n(X) Y_2\zeta_n(Y)$ converges to $X_2 \zeta(X) Y_2\zeta(Y)$ weakly in $L^p(B_L \times B_L)$ for $1<p<\infty$ and weakly-$*$ in $L^\infty(B_L \times B_L)$ for $p = \infty$.
This ensures that the double integral in the inequality above vanishes as $n \to \infty$.
Letting $\varepsilon \to 0$ concludes $\lim_{n \to \infty} E[\zeta_n] = E[\zeta]$.
This guarantees $E[\zeta] = I_{p,1,1}$, confirming that $\zeta \in S_{p,1,1}$.\\

We show the stronger convergence \eqref{eq:conv2} for $1<p<\infty$.
Since $\zeta \in S_{p,1, 1}$, the norm constraint must be saturated by Remark \ref{sec:saturation} and hence $Q_p[\zeta] = 1$.
In particular, $\lim_{n \to \infty} Q_p[\zeta_n] = Q_p[\zeta]$. The weak convergence $\zeta_n \rightharpoonup \zeta$ combined with this norm convergence ensures the strong convergence \eqref{eq:conv2} by uniform convexity of $L^p(\bbR^2_+; r)$. Applying H\"older's inequality yields 
\begin{align*}
P[|\zeta_n-\zeta|]
=(X_2^{2}\mathbf{1}_{B_L(0)},|\zeta_n-\zeta|)+P[|\zeta_n-\zeta|\mathbf{1}_{B_L(0)^{c}}]
\le ||X_2^{2}\mathbf{1}_{B_L(0)}||_{L^{p'}_r}||\zeta_n-\zeta||_{L^{p}_r}+P[|\zeta_n-\zeta|\mathbf{1}_{B_L(0)^{c}}].
\end{align*}
The first term vanishes as $n\to\infty$. The limsup of the second term is at most $2\varepsilon$ by \eqref{eq:tight}. Letting $n \to \infty$ and then $\varepsilon \to 0$ yields the convergence \eqref{eq:conv1}.

It remains to prove the convergence \eqref{eq:conv1} for $p = \infty$. In this case, the maximizer is Hill's vortex $\zeta = \mathbf{1}_{\Omega}(\cdot+c,\cdot)$ for the set $\Omega=B_{\varrho_{H}}(0)_+$ and some constant $c$. By translating the sequence $\{\zeta_n\}$
along the $z$-direction, we may assume that $c=0$. We take an arbitrary $\delta>0$. Since $Q_\infty[\zeta_n]\to 1$, we may assume that $0 \le \zeta_n \le 1+\delta$ for all $n$. It follows that 
\begin{align*}
P[|\zeta_n-\zeta|]&=P[|\zeta_n-\zeta|\mathbf{1}_{\Omega}]+P[|\zeta_n-\zeta|\mathbf{1}_{\Omega^{c}}]\\
&=P[(\zeta_n-1)\mathbf{1}_{\{1<\zeta_n\le 1+\delta\}}\mathbf{1}_{\Omega}]
+P[(1-\zeta_n)\mathbf{1}_{\{\zeta_n\le 1\}}\mathbf{1}_{\Omega}]+P[\zeta_n\mathbf{1}_{\Omega^{c}}] \\
&\le \delta+P[(1-\zeta_n)\mathbf{1}_{\Omega}]+P[\zeta_n\mathbf{1}_{\Omega^{c}}]
=\delta+1+P[\zeta_n]-2P[\zeta_n\mathbf{1}_{\Omega}].
\end{align*}
Since $\Omega$ is bounded, the convergence \eqref{eq:conv2} ensures $P[\zeta_n\mathbf{1}_{\Omega}] \to P[\zeta\mathbf{1}_{\Omega}] = 1$. Since $P[\zeta_n] \to 1$, letting $n \to \infty$ followed by $\delta \to 0$ implies the desired convergence \eqref{eq:conv1}. The proof is complete.
\end{proof}

\begin{proof}[Proof of Theorem \ref{t:MT1}]
The existence of energy maximizers follows from Theorem \ref{thm:existence}. The compactness of their supports, Steiner symmetry, the Euler--Lagrange equation, and the identities for the Lagrange multipliers follow from Lemmas \ref{lem:compact}, \ref{lem:steiner}, \ref{lem:EL_p}, and \ref{lem:pohozaev}, respectively.
\end{proof}

\section{Orbital stability}\label{s:5}

We prove the orbital stability results stated in Theorems \ref{t:MT2} and \ref{t:MT3}(iv).
We first review Danchin's global well-posedness result for the Euler equations with axisymmetric initial data without swirl, and show that the resulting global solutions conserve the Casimir invariants, as well as energy and impulse.
Since the velocity field of Danchin's solutions is log-Lipschitz continuous, the associated flow map is uniquely determined, and the distribution function of the relative vorticity is preserved.
This in turn implies the conservation of the Casimir invariants.

\subsection{Conservation laws for global-in-time solutions}

We show conservation laws for global-in-time solutions to the axisymmetric Euler equations without swirl \eqref{eq:Euler1}.
We define the Zygmund space and the log-Lipschitz space by 
\begin{align}
C_{\star}^1(\mathbb R^3)
&=
\left\{
u\in L^\infty(\mathbb R^3):
\|u\|_{C_{\star}^1(\mathbb R^3)}<\infty
\right\}, \\
\|u\|_{C_{\star}^1(\mathbb R^3)}
&=
\|u\|_{L^\infty(\mathbb R^3)}
+
\sup_{\substack{x\in\mathbb R^3\\h\in\mathbb R^3\setminus\{0\}}}
\frac{|u(x+h)+u(x-h)-2u(x)|}{|h|},\\
LL(\mathbb R^3)
&=
\left\{
u\in L^\infty(\mathbb R^3):
\|u\|_{LL(\mathbb{R}^{3})}<\infty
\right\}, \\
\|u\|_{LL(\mathbb{R}^{3})}
&=
\|u\|_{L^\infty}
+
\sup_{\substack{x,y\in\mathbb R^3\\x\neq y}}
\frac{|u(x)-u(y)|}
{|x-y|\bigl(1+\log^+|x-y|^{-1}\bigr)},
\end{align}
where $\log^+s=(\log s)_+$. The space $C_{\star}^{1}(\mathbb{R}^{3})$ is continuously embedded into $LL(\mathbb{R}^{3})$; see \cite[Proposition 2.107]{BCD11}. For a Banach space $\mathcal{X}$, we denote by
$BC(\mathbb{R};\mathcal{X})$ the space of bounded continuous
$\mathcal{X}$-valued functions on $\mathbb{R}$,
and by $C_w(\mathbb{R};\mathcal{X})$ the space of
weakly continuous $\mathcal{X}$-valued functions on $\mathbb{R}$. The following global well-posedness result is due to Danchin \cite[Theorem 1]{Danchin}.

\begin{proposition}\label{p:Danchin}
For axisymmetric initial data $\zeta_0\in L^{3,1}(\mathbb{R}^{3})$ satisfying $r\zeta_0\in L^{\infty}\cap L^{3,1}(\mathbb{R}^{3})$, there exists a unique global-in-time axisymmetric solution $\zeta\in BC(\mathbb{R}; L^{3,1}(\mathbb{R}^{3}))$ of \eqref{eq:Euler1} satisfying 
\begin{align}
&u=\nabla \times (\mathcal{G}\zeta\nabla \theta)\in C_{w}(\mathbb{R}; C^{1}_{\star}(\mathbb{R}^{3})),  \label{eq:Danchinvelocity}\\
&r\zeta\in L^{\infty}_{\textrm{loc}}(\mathbb{R}; L^{\infty}\cap L^{3,1}(\mathbb{R}^{3})),\\
&\frac{1}{r^{2}}\partial_z \mathcal{G}\zeta\in BC(\mathbb{R}; L^{\infty}(\mathbb{R}^{3})). \label{eq:DanchinBS}
\end{align}
The solution $\zeta$ satisfies \eqref{eq:Euler1} in the distributional sense. Moreover, for the flow map $\chi(x,t)$ associated with the velocity field \eqref{eq:Danchinvelocity}, we have 
\begin{align}
\zeta(\chi(x,t),t)=\zeta_0(x)\quad \textrm{a.e. in}\ \mathbb{R}^{3}. \label{eq:Lagrangean}
\end{align}
In particular, for every $\lambda>0$, 
\begin{align}
d_{\zeta(t)}(\lambda)
=
d_{\zeta_0}(\lambda),  \label{eq:distribution}
\end{align}
where $d_f(\lambda)=|\{x\in\mathbb R^3:|f(x)|>\lambda\}|$ denotes the distribution function of $f$.
\end{proposition}

\begin{proof}
The existence of a unique global-in-time solution and all the stated properties, except for the continuity of $\zeta$ in $L^{3,1}(\mathbb R^3)$ and that of $r^{-2}\partial_z\mathcal{G}\zeta$ in $L^\infty(\mathbb R^3)$, are established in \cite{Danchin}.

By \eqref{eq:Danchinvelocity}, the velocity field takes values in the Zygmund space and is therefore log-Lipschitz continuous. Consequently, the associated flow map $\chi(x,t)$ is uniquely determined. For each $\lambda > 0$, let
\[
A_\lambda(t)
=
\{x\in\mathbb R^3:|\zeta(x,t)|>\lambda\}.
\]
Then, by \eqref{eq:Lagrangean},
\[
\begin{aligned}
\chi_t^{-1}(A_\lambda(t))
=\{x\in\mathbb R^3:\chi(x,t)\in A_\lambda(t)\}
&=\{x\in\mathbb R^3:|\zeta(\chi(x,t),t)|>\lambda\} \\
&=\{x\in\mathbb R^3:|\zeta_0(x)|>\lambda\}
=A_\lambda(0),
\end{aligned}
\]
where $\chi_t=\chi(x,t)$.
Since $\chi_t$ preserves Lebesgue measure, we have $|A_\lambda(t)| = |A_\lambda(0)|$ and \eqref{eq:distribution} holds.
The conservation of the $L^{3,1}$-norm then follows from the layer-cake representation, obtained by Tonelli's theorem,
\[
\|f\|_{L^{3,1}(\mathbb{R}^{3})}
=
3\int_0^\infty d_f(\lambda)^{1/3}\,d\lambda.
\]
In particular, for any $f\in L^{3,1}(\mathbb R^3)$,
\[
\|f\circ \chi(\cdot,t)\|_{L^{3,1}(\mathbb{R}^{3})}
=
\|f\|_{L^{3,1}(\mathbb{R}^{3})}.
\]
The continuity of $\zeta$ in $L^{3,1}(\mathbb R^3)$ follows from this invariance of the $L^{3,1}$-norm under composition with the flow map, together with the density of smooth compactly supported functions in $L^{3,1}(\mathbb R^3)$.
The Green's function estimate \eqref{eq:DanchinGreen} implies that the operator $r^{-2}\partial_z\mathcal G$ is bounded from $L^{3,1}(\mathbb R^3)$ to $L^\infty(\mathbb R^3)$.
This yields the continuity in \eqref{eq:DanchinBS}.
\end{proof}

\begin{proposition}[Casimir invariants]
Let $f:[0,\infty)\to \mathbb{R}$ be a Borel function such that $f(0)=0$ and let 
\begin{align}
C_{f}[\zeta]=\int_{\mathbb{R}^{3}}f(\zeta)dx.  \label{eq:fcasimir}
\end{align}
Under the assumptions of Proposition \ref{p:Danchin}, assume in addition that $\zeta_0$ is nonnegative and $C_f[\zeta_0]$ is finite.
Then we have 
\begin{align}
C_{f}[\zeta(t)]=C_{f}[\zeta_0],\quad t\in \mathbb{R}.  \label{eq:casimirconservation}
\end{align}
\end{proposition}

\begin{proof}
For the pushforward measure $\mu_t=\mathcal{D}(\zeta(t))$ and $\lambda>0$, we have
\begin{align*}
\mu_t(\lambda,\infty)=|\{x\in \mathbb{R}^{3}: \zeta(x,t)>\lambda\}|
=d_{\zeta(t)}(\lambda).
\end{align*}
By the conservation of the distribution function \eqref{eq:distribution}, $\mu_t=\mu_0$ and 
\begin{align*}
C_{f}[\zeta(t)]=\int_{\mathbb{R}^{3}}f(\zeta(t))dx=\int_{0}^{\infty}f(\zeta)d\mu_t(\zeta)=\int_{0}^{\infty}f(\zeta)d\mu_0(\zeta)=C_{f}[\zeta_0].
\end{align*}
\end{proof}

\begin{proposition}\label{p:conservation}
Under the assumptions of Proposition \ref{p:Danchin}, assume in addition that $\zeta_0$ is nonnegative and $\zeta_0\in L^{p}\cap L^{1}_{*} (\mathbb{R}^{3})$ for $1< p\le \infty$. Then we have  
\begin{align}
E[\zeta(t)]=E[\zeta_0],\quad 
P[\zeta(t)]=P[\zeta_0],\quad 
Q_p[\zeta(t)]=Q_p[\zeta_0],\quad t\in \mathbb{R}. \label{eq:conservation}
\end{align}
\end{proposition}

\begin{proof}
The conservation of the $L^p$-Casimir follows from \eqref{eq:casimirconservation}.
We write $x\in\mathbb R^3$, the velocity $u(x,t)$, and the flow map $\chi(x,t)$ in cylindrical coordinates as 
\begin{align*}
x&=re_r(\theta)+ze_z,\\
u(x,t)&=u^r(x,t)e_r(\theta)+u^z(x,t)e_z,\\
\chi(x,t)&=R(x,t)e_r(\theta)+Z(x,t)e_z,
\end{align*}
using the basis $e_r(\theta) = (\cos\theta,\sin\theta,0)$ and $e_z=(0,0,1)$.
By the flow equation, $R(x,t)$ satisfies
\[
\partial_t\log R(x,t)
=\left(\frac{u^r}{r}\right)(\chi(x,t),t).
\]
Since $u^r/r$ is bounded by \eqref{eq:DanchinBS}, there exists a constant $C>0$ such that $R(x,t)\le re^{C|t|}$.
By a change of variables, we have 
\[
\begin{aligned}
P[\zeta(t)]
=\int_{\mathbb R^3}r^2\zeta(x,t)\,dx
=\int_{\mathbb R^3}R(x,t)^2\zeta_0(x)\,dx.
\end{aligned}
\]
Hence $P[\zeta(t)]$ is finite for every $t$.
Moreover, the energy inequality \eqref{eq:variational-inequality2} implies that $E[\zeta(t)]$ is also finite.
In particular, $r^{-1}\nabla\psi\in L^\infty_{\mathrm{loc}}(\mathbb R;L^2(\mathbb R^3))$ for $\psi=\mathcal{G}\zeta$.

By differentiating $P[\zeta(t)]$, we have  
\[
\begin{aligned}
\frac{d}{dt}P[\zeta(t)]
=2\int_{\mathbb R^3}
R(x,t)\partial_tR(x,t)\zeta_0(x)\,dx
=2\int_{\mathbb R^3}
ru^r(x,t)\zeta(x,t)\,dx
=-2\int_{\mathbb R^3}
\partial_z\psi(x,t)\zeta(x,t)\,dx.
\end{aligned}
\]
This integral converges by the condition \eqref{eq:DanchinBS} and $r^2\zeta\in L^\infty_{\textrm{loc}}(\mathbb{R}; L^1(\mathbb{R}^{3}))$. Using the five-dimensional transformation \eqref{eq:5dtrans}, we have 
\[
\partial_{y_1}\varphi\in L^\infty_{\textrm{loc}}(\mathbb{R};L^\infty(\mathbb R^5)),
\qquad
\Delta\varphi\in L^\infty_{\textrm{loc}}(\mathbb{R};L^1(\mathbb R^5)),
\]
and the last integral coincides with 
\begin{align*}
\frac{2}{\pi}\int_{\mathbb R^5}
\partial_{y_1}\varphi\,\Delta_y\varphi\,dy.
\end{align*}
For $h\neq 0$ and the symmetric difference quotient 
\[
D_h\varphi(y)
=
\frac{\varphi(y_1+h,y')-\varphi(y_1-h,y')}{2h},
\]
by the symmetry of the Laplacian under translations, 
we have 
\[
\int_{\mathbb R^5}D_h\varphi\,\Delta_y\varphi\,dy=0.
\]
Thus, letting $h\to 0$ implies $\frac{d}{dt}P[\zeta(t)]=0$ and conservation of the impulse.

It remains to show the conservation of kinetic energy. For $X=(X_1,X_2)$ and $Y=(Y_1,Y_2)$, the Green's function estimate \eqref{eq:DanchinGreen} yields
\begin{align*}
|\nabla_XG(X,Y)|\lesssim Y_2^{2}\int_{0}^{2\pi}\frac{d\theta}{|x-y|^{2}},
\end{align*}
for $x=(X_2\cos\theta,X_2\sin\theta,X_1)$ and $y=(Y_2,0,Y_1)$. It follows that 
\begin{align*}
\int_{ \mathbb{R}^{3}}|\nabla_XG(X,Y)|\zeta(x,t)dx
&=2\pi \int_{ \mathbb{R}^{2}_+}|\nabla_XG(X,Y)|\zeta(X,t)X_2dX \\
&\lesssim  Y_2^{2}
\int_{ \mathbb{R}^{2}_+}\int_{0}^{2\pi}\frac{d\theta}{|x-y|^{2}}
\zeta(X,t)X_2dX \\
&\lesssim Y_2^{2} \int_{ \mathbb{R}^{3}}\frac{1}{|x-y|^{2}}\zeta(x,t)dx=Y_2^{2}\frac{1}{|y|^{2}}*\zeta(\cdot,t).
\end{align*}
Here, $\zeta(x,t)$ is also written as $\zeta(X,t)$ by axisymmetry. The resulting convolution is axisymmetric.
By applying Hölder's inequality in Lorentz spaces, we estimate 
\begin{align*}
\left\|\frac{1}{|y|^{2}}*\zeta(\cdot,t)\right\|_{L^{\infty}(\mathbb{R}^{3})}\lesssim  \left\|\frac{1}{|y|^{2}}\right\|_{L^{3/2,\infty}(\mathbb{R}^{3})}\left\|\zeta(t)\right\|_{L^{3,1}(\mathbb{R}^{3})}.
\end{align*}
It follows that 
\begin{align*}
\iint_{\mathbb{R}^{3}\times \mathbb{R}^{3}}|u(x,t)||\nabla_XG(X,Y)|\zeta(x,t)\zeta(y,t)dxdy
\lesssim ||u(t)||_{L^{\infty}(\mathbb{R}^{3})} ||\zeta(t)||_{L^{3,1}(\mathbb{R}^{3})}P[\zeta(t)].
\end{align*}
By a change of variables, we have  
\begin{equation*}
\begin{aligned}
E[\zeta(t)]
=\frac{1}{4\pi}\iint_{\mathbb{R}^{3}\times \mathbb{R}^{3}}G(Z(x,t),R(x,t),Z(y,t),R(y,t))\zeta_0(x)\zeta_0(y)dxdy.
\end{aligned}
\end{equation*}
The preceding estimate, together with the local-in-time bounds for $u$, $\zeta$, and $P[\zeta]$, justifies differentiating the Lagrangian representation of the energy. More precisely, the chain rule along particle trajectories and Fubini's theorem imply that $E[\zeta(t)]$ is locally absolutely continuous. Consequently, for almost every $t$, we have  
\begin{align*}
\frac{d}{dt}E[\zeta(t)]
&=\frac{1}{4\pi}\iint_{\mathbb{R}^{3}\times \mathbb{R}^{3}}u(\chi(x,t),t)\cdot \nabla_{X}G(Z(x,t),R(x,t),Z(y,t),R(y,t))\zeta_0(x)\zeta_0(y)dxdy\\
&+\frac{1}{4\pi}\iint_{\mathbb{R}^{3}\times \mathbb{R}^{3}}u(\chi(y,t),t)\cdot \nabla_{Y}G(Z(x,t),R(x,t),Z(y,t),R(y,t))
\zeta_0(x)\zeta_0(y)dxdy.
\end{align*}
Using the symmetry of the Green's function $G(X,Y)=G(Y,X)$, we have $\nabla_{Y}G(X,Y)=(\nabla_{X}G)(Y,X)$. Exchanging the integration variables $x$ and $y$, we find that the two integrals are equal. It follows that 
\begin{align*}
&\iint_{\mathbb{R}^{3}\times \mathbb{R}^{3}}u(\chi(x,t),t)\cdot \nabla_{X}G(Z(x,t),R(x,t),Z(y,t),R(y,t))\zeta_0(x)\zeta_0(y)dxdy \\
&=\iint_{\mathbb{R}^{3}\times \mathbb{R}^{3}}u(x,t)\cdot \nabla_{X}G(X,Y)\zeta(x,t)\zeta(y,t)dxdy 
=2\pi \int_{\mathbb{R}^{3}}(u \cdot \nabla\mathcal{G}\zeta) \zeta dx=0.
\end{align*}
Consequently, $\frac{d}{dt}E[\zeta(t)]=0$ for almost every $t$. This shows the conservation of the kinetic energy.
\end{proof}

\begin{remark}[Weak solutions]
Nobili and Seis \cite[Theorems 1-3]{NS22} constructed global weak solutions to the Euler equations \eqref{eq:Euler1} for axisymmetric initial data without swirl $u_0\in L^{2}_{\textrm{loc}}(\mathbb{R}^{3})$ satisfying $\zeta_0 \in L^p\cap L^{1}(\bbR^3)$ for $1<p<\infty$ by a vanishing viscosity limit of the Navier--Stokes flows. The constructed weak solution satisfies the renormalization property of DiPerna and Lions and conserves the $L^{p}$-norm of relative vorticity. Moreover, if $u_0\in L^{2}(\mathbb{R}^{3})$ and the nonnegative initial relative vorticity $\zeta_0\in L^{p}\cap L^{1}\cap L^{1}_{*}(\mathbb{R}^{3})$ for $3/2<p<\infty$ has finite impulse, the constructed weak solution conserves the kinetic energy. See also \cite{Abe20}.
\end{remark}

\subsection{Stability for \texorpdfstring{$1<p<\infty$}{1<p<infinity}}

The compactness theorem (Theorem \ref{thm:existence}) and the conservation laws (Proposition \ref{p:conservation}) imply the orbital stability of the set of energy maximizers $S_p$ for $1<p<\infty$.

\begin{proof}[Proof of Theorem \ref{t:MT2}]
Suppose, for a contradiction, that the assertion is false.
Then there exists $\varepsilon_0 > 0$, a sequence of nonnegative axisymmetric initial data $\{\zeta_{0,n}\} \subset L^{p}\cap L^{1}_{*}(\bbR^3)$ satisfying \eqref{eq:Danchin}, and a sequence of times $\{t_n\} \subset \bbR $ such that  
\begin{equation*}
\inf_{\eta\in S_p} \|\zeta_{0,n}-\eta\|_{L^p(\bbR^3)} + |P[\zeta_{0,n}]-1|  \le \frac{1}{n},
\end{equation*}
while the corresponding global solution $\zeta_n (t)$ satisfies
\begin{equation*}
\inf_{\eta\in S_p} \|\zeta_n(t_n)-\eta\|_{L^p\cap L^{1}_{*}(\bbR^3)}  \ge \varepsilon_0.
\end{equation*}
For $\eta\in S_p$, we have $E[\eta]=I_p$ and $P[\eta] = 1$. Applying \eqref{eq:en_diff} yields
\begin{align*}
\lim_{n\to \infty}|E[\zeta_{0,n}]-I_p|= 0.
\end{align*}
We set $\zeta_n=\zeta_n(t_n)$. By the conservation laws \eqref{eq:conservation}, we have
\[
\lim_{n\to\infty}E[\zeta_n] = I_p, \quad \lim_{n\to\infty}P[\zeta_n] = 1, \quad \limsup_{n\to\infty}Q_p[\zeta_n] \le 1.
\]
Applying Theorem~\ref{thm:existence}, there exists a subsequence and some $\zeta\in S_p$ such that, up to translation, we have
\[
\zeta_n \to \zeta \quad \text{in} \,\, L^p\cap L^1_*(\bbR^3).
\]
Therefore,
\begin{align*}
0=\lim_{n\to \infty}\|\zeta_n-\zeta\|_{L^p\cap L^1_{*}(\bbR^3)} \ge \liminf_{n\to \infty}\inf_{\eta\in S_p}  \|\zeta_n-\eta\|_{L^p\cap L^1_{*}(\bbR^3)}  \ge \varepsilon_0.
\end{align*}
This is a contradiction. Thus the theorem is proved.
\end{proof}

\begin{remark}[Comparison with previous stability results]\label{r:stabilitycomparison}
Previous orbital stability results for individual steadily translating vortex rings \cite{Choi2024,CQZZ23,CQYZZ22,CLQZZ26}
assume that the nonnegative initial relative vorticity
$\zeta_0$ satisfies
\[
\zeta_0\in L^\infty \cap L^1\cap L^{1}_{*}(\mathbb R^3),
\qquad
r\zeta_0\in L^\infty(\mathbb R^3),
\]
and establish stability by measuring perturbations in the norm of   
\[
L^2\cap L^1\cap L^{1}_{*} (\mathbb R^3).
\]
This regularity assumption is subcritical for the global well-posedness of the Euler equations for axisymmetric data without swirl. Under these assumptions, all $L^q$ Casimirs,
$1\le q\le\infty$, are finite and conserved.

In contrast, Theorem \ref{t:MT2} establishes orbital stability
of the entire set of energy maximizers,
without assuming uniqueness of the maximizer.
Our stability theorem does not require
bounded relative vorticity or finite circulation. Instead, for the Casimir exponent $1<p<\infty$
used in the energy maximization problem \eqref{eq:variational},
we assume
\[
\zeta_0\in L^{3,1}\cap L^p\cap L^{1}_{*}(\mathbb R^3),
\qquad
r\zeta_0\in L^\infty
\cap L^{3,1}(\mathbb R^3),
\]
and establish
orbital stability by measuring perturbations in the norm of 
\[
L^p\cap L^{1}_{*}(\mathbb R^3).
\]
The replacement of the $L^\infty$ assumption on
$\zeta_0$ by an $L^{3,1}$ condition is made possible
by Danchin's global well-posedness result. See also \cite{KJ22,BC26} for recent ill-posedness results for $\zeta_0\in L^{3,q}(\mathbb{R}^{3})$ for $q>1$ and \cite{DJ26} for the global well-posedness result under the condition $r\zeta_0\in W^{3,1}(\mathbb{R}^{3})$. On the other hand, the removal of the finite
circulation assumption follows from the energy maximization principle \eqref{eq:variational} based only on the impulse and a single $L^p$ Casimir. The initial condition \eqref{eq:stability}
requires only that the impulse of $\zeta_0$ be close
to that of an energy maximizer, rather than requiring
the stronger condition that
$P[|\zeta_0-\eta|]$ be small, cf. \cite[Remark 1.2]{AC2022}.
\end{remark}

\subsection{Stability for \texorpdfstring{$p=\infty$}{p=infinity}}

By Theorem \ref{t:MT3}(iii), the energy maximizers are precisely the translations of Hill's vortex.
Thus, the orbital stability of the set of energy maximizers reduces to that of Hill's vortex.
As in the case of finite $p$ with $1<p<\infty$, orbital stability can be established using the compactness of energy-maximizing sequences.
However, orbital stability does not hold in the $L^\infty$ norm. We therefore establish orbital stability in $L^q\cap L^1_*(\mathbb R^3)$ for every finite $1<q<\infty$.

A direct application of the compactness theorem requires an upper bound on the essential supremum of the initial data $Q_{\infty}[\zeta_{0}]\le 1$.
We first establish stability under this additional constraint and then show that it can be removed by exploiting the conservation of the Casimir invariants.

\begin{proposition} \label{p:conditional}
The assertion of Theorem~\ref{t:MT3}(iv) holds
under the additional assumption
$Q_\infty[\zeta_0]\le 1$.
\end{proposition}

\begin{proof}
By Theorem \ref{t:MT3}(iii), we have
\[
S_{\infty}
=
\{\zeta_H(\cdot+c,\cdot):c\in\mathbb R\}.
\]
We thus prove the orbital stability of the set $S_{\infty}$. Arguing by contradiction as in the proof of Theorem \ref{t:MT2},
we obtain a sequence of initial data $\{\zeta_{0,n}\}
\subset
L^\infty\cap L^q\cap L^1_*(\mathbb R^3)$ satisfying  \eqref{eq:Danchin} and a sequence of times $\{t_n\}$ such that 
\begin{equation}
\inf_{c\in \mathbb{R}}
\|\zeta_n(t_n)-\zeta_H(\cdot+c,\cdot)\|_{L^{q}\cap L^{1}_{*} (\mathbb{R}^3)}
>\varepsilon_0.
\end{equation}
Since $\zeta_{0,n}$ converges to Hill's vortex
$\zeta_H$ in $L^q(\mathbb R^3)$, the energy inequality \eqref{eq:en_diff} implies $E[\zeta_{0,n}]\to I_\infty$. Moreover, the initial data satisfy $P[\zeta_{0,n}]\to 1$ and $Q_\infty[\zeta_{0,n}]\le 1$. By the conservation laws \eqref{eq:conservation}, the sequence
$\zeta_n=\zeta_n(t_n)$ satisfies
\[
\lim_{n\to\infty}E[\zeta_n]=I_\infty,
\qquad
\lim_{n\to\infty}P[\zeta_n]=1,
\qquad
Q_\infty[\zeta_n]\le 1.
\]
Applying Theorem \ref{thm:existence}, we may pass to a subsequence
and apply suitable translations so that $\zeta_n\to\zeta_H$ in $L^1_*(\mathbb R^3)$. Together with the conservation of the Casimir invariants,
this implies, after passing to a further subsequence if necessary,
that
\begin{align*}
&\zeta_n\to\zeta_H
\quad\text{a.e. in }\mathbb R^3, \\
&Q_\infty[\zeta_n]\le 1, \\
\qquad
&Q_q[\zeta_n]
=
Q_q[\zeta_{0,n}]
\to Q_q[\zeta_H].
\end{align*}
Since $\zeta_H=\mathbf 1_{B_{\varrho_H}}$ for the open ball $B_{\varrho_H}=B_{\varrho_H}(0)$ of radius $\varrho_{H}$,
we obtain
\begin{align*}
\int_{\mathbb R^3}|\zeta_n-\zeta_H|^q\,dx
&=
\int_{B_{\varrho_H}}|\zeta_n-1|^q\,dx
+
\int_{\mathbb R^3\setminus B_{\varrho_H}}
|\zeta_n|^q\,dx\\
&=
\int_{B_{\varrho_H}}|\zeta_n-1|^q\,dx
+
Q_q^{q}[\zeta_{n}]
-
\int_{B_{\varrho_H}}|\zeta_n|^q\,dx
\to 0.
\end{align*}
Therefore, $\zeta_n\to\zeta_H$ in $L^q\cap L^1_*(\mathbb R^3)$, which yields the desired contradiction.
\end{proof}

\begin{proposition}
Let $\zeta(t)$ be a global-in-time solution to \eqref{eq:Euler1} for nonnegative axisymmetric initial data
$\zeta_0\in L^{q}\cap L^{1}_{*} (\mathbb{R}^{3})$ for $1<q<\infty$ satisfying \eqref{eq:Danchin}. We define $\xi^{i}(t)$ for $i=1,2$ by 
\begin{align}
\zeta(t)=\min\{1,\zeta(t)\}+(\zeta(t)-1)_{+}=\xi^{1}(t)+\xi^{2}(t).  \label{eq:dec}
\end{align}
Then, for $i=1,2$, 
\begin{align}
Q_{q}[\xi^{i}(t)]=Q_{q}[\xi^{i}_{0}],\quad t\in \mathbb{R}. \label{eq:casimirconserve2}
\end{align}
\end{proposition}

\begin{proof}
Since $\xi^{i}=g_i(\zeta)$, where $g_1(s)=\min\{1,s\}$ and $g_2(s)=(s-1)_{+}$, the Casimir conservation \eqref{eq:casimirconservation} yields 
\begin{align*}
Q_{q}[\xi^{i}(t)]=Q_{q}[g_i(\zeta(t))]=Q_{q}[g_i(\zeta_0)]=Q_{q}[\xi_0^{i}].
\end{align*}
This proves \eqref{eq:casimirconserve2}.
\end{proof}

\begin{proof}[Proof of Theorem \ref{t:MT3}(iv)]
It suffices to remove the conditions $\zeta_0\in L^{\infty}(\mathbb{R}^{3})$ and $Q_{\infty}[\zeta_0]\le 1$ in Proposition \ref{p:conditional}. Arguing by contradiction as in the proof of Proposition \ref{p:conditional},
we obtain a sequence of initial data $\{\zeta_{0,n}\}\subset L^q\cap L^1_*(\mathbb R^3)$ satisfying \eqref{eq:Danchin} and a sequence of times $\{t_n\}$ such that
\[
\zeta_{0,n}\to\zeta_H
\quad\text{in }L^q(\mathbb R^3),
\qquad
P[\zeta_{0,n}]\to 1,
\]
while the corresponding solutions $\zeta_n(t)$ satisfy
\[
\inf_{c\in\mathbb R}
\|\zeta_n(t_n)-\zeta_H(\cdot+c,\cdot)\|
_{L^q\cap L^1_*(\mathbb R^3)}
> \varepsilon_0
\]
for some $\varepsilon_0>0$. Let $\zeta_n(t)=\xi_n^1(t)+\xi_n^2(t)$ be the decomposition defined in \eqref{eq:dec}.
Then
\begin{align*}
&P[\xi_{0,n}^1]+P[\xi_{0,n}^2]
=P[\zeta_{0,n}]\to 1,\\
&Q_q^q[\zeta_{0,n}-\xi_{0,n}^1]
=Q_q^q[\xi_{0,n}^2]
=\int_{\mathbb R^3}(\zeta_{0,n}-1)_+^q\,dx
\le
\int_{\mathbb R^3}(\zeta_{0,n}-\zeta_H)_+^q\,dx
\le
\|\zeta_{0,n}-\zeta_H\|_{L^q(\mathbb R^3)}^q
\to 0.
\end{align*}
Moreover, the energy estimate \eqref{eq:en_diff} implies that $E[\zeta_{0,n}]\to I_\infty$.

We set $\zeta_n=\zeta_n(t_n)$ and
$\xi_n^i=\xi_n^i(t_n)$ for $i=1,2$.
By conservation of energy, $E[\zeta_n]\to I_\infty$. Furthermore, Casimir conservation \eqref{eq:casimirconserve2} yields
\begin{align*}
\lim_{n\to\infty}Q_q[\xi_n^2]
&=
\lim_{n\to\infty}Q_q[\xi_{0,n}^2]
=0,\\
\lim_{n\to\infty}Q_q[\xi_n^1]
&=
\lim_{n\to\infty}Q_q[\xi_{0,n}^1]
= Q_q[\zeta_H].
\end{align*}
Applying the energy estimates \eqref{eq:en_diff}-\eqref{eq:en_cross} for 
\[
E[\zeta_n]
=
E[\xi_n^1]
+2E_{\mathrm{int}}[\xi_n^1,\xi_n^2]
+E[\xi_n^2],
\]
we obtain
\[
\lim_{n\to\infty}E[\xi_n^1]=I_\infty,
\qquad
\limsup_{n\to\infty}P[\xi_n^1]\le 1,
\qquad
Q_\infty[\xi_n^1]\le 1.
\]
By Theorem \ref{thm:existence}, after passing to a subsequence and applying
suitable translations, we have $\xi_n^1\to\zeta_H$
in $L^1_*(\mathbb R^3)$. In particular, after passing to a further subsequence
if necessary,
\begin{align*}
\xi_n^1&\to\zeta_H
\quad\text{a.e. in }\mathbb R^3,\\
Q_\infty[\xi_n^1]&\le 1,\\
Q_q[\xi_n^1]
&=
Q_q[\xi_{0,n}^1]
\to Q_q[\zeta_H].
\end{align*}
Arguing as in the proof of Proposition \ref{p:conditional}, we obtain the convergence $\xi_n^1\to\zeta_H$ in $L^q\cap L^1_*(\mathbb R^3)$.
Consequently,
\begin{align*}
P[\xi_n^2]
&=
P[\zeta_n]-P[\xi_n^1]
\to 0,\\
Q_q[\xi_n^2]
&=
Q_q[\xi_{0,n}^2]
\to 0.
\end{align*}
Hence $\xi_n^2\to 0$ in $L^q\cap L^1_*(\mathbb R^3)$.
We conclude that $\zeta_n\to\zeta_H$ in $L^q\cap L^1_*(\mathbb R^3)$, which contradicts the choice of $t_n$ and $\zeta_{0,n}$.
This completes the proof.
\end{proof}

\begin{remark}
The stability results in Theorems \ref{t:MT2} and \ref{t:MT3}(iv) remain valid for weak solutions with initial data of lower regularity than $L^{3,1}(\mathbb{R}^{3})$,
provided that the conservation laws hold.
In this case, the results guarantee the existence of at least one weak solution satisfying the stability estimates.
See \cite{BNL13,AC2022,ACJ2025}.
\end{remark}

\section{Spherical vortex rings and the limit \texorpdfstring{$p\to\infty$}{p → ∞}}\label{s:6}

We complete the proof of Theorem \ref{t:MT3}. The quantities $I_p$, $W_p$, and $\kappa_p$ remain uniformly bounded as $p\to\infty$ by uniform estimates established in \S\ref{s:3}.
Together with the uniform $L^\infty$ bounds for the energy maximizers of $I_p$, this implies that these maximizers form a maximizing sequence for $I_\infty$ as $p\to\infty$.
Consequently, the energy maximizers and the associated constants converge to their counterparts for Hill's vortex.

The corresponding stream functions converge to the stream function of Hill's vortex in $C^1(\mathbb R^5)$ under the five-dimensional transformation.
This implies that the energy maximizers converge locally uniformly to Hill's vortex away from its boundary and the symmetry axis inside the vortex core. Furthermore, using the strict radial monotonicity of the stream function of Hill's vortex in five dimensions, we obtain a graphical representation of the vortex cores by the implicit function theorem.

\subsection{Convergence to Hill's vortex}

We first use the compactness of energy-maximizing sequences to show that the energy maximizers of $I_p$ converge to Hill's vortex in a weak topology.
We then exploit the Euler--Lagrange equation to upgrade this to local uniform convergence away from the vortex boundary and the symmetry axis inside the vortex.
We say that a maximizer is centered if it is Steiner symmetric with respect to $z=0$.

\begin{lemma}[Convergence to Hill's vortex]\label{lem:conv_hill}
Let $\zeta_H$ denote Hill's vortex \eqref{eq:normalizedHill0} and let $W_\infty$ be the constant defined in \eqref{eq:speedHill}.
Let $p_k\to \infty$.
Let $I_{p_k}$, $W_{p_k}$, and $\kappa_{p_k}$ be the constants defined in \eqref{eq:variational-inequality} and \eqref{eq:speed}.
Let $\zeta_k\in S_{p_k}$ be a centered maximizer.
Then, up to a subsequence, the following convergences hold: 
\begin{align}
I_{p_k} & \to  I_\infty,  \label{eq:lar_p_en_limit} \\
W_{p_k} &\to W_\infty,  
\label{eq:Wconv}\\
\kappa_{p_k}&\to \kappa_{\infty},\label{eq:kappaconv}\\
\zeta_k &\to \zeta_H \quad \text{in}\,\, L^1_*(\bbR^3), \quad \zeta_k \rightharpoonup^* \zeta_H \quad \text{in}\,\, L^\infty(\bbR^3). \label{eq:lar_p_vor_lim} 
\end{align} 
\end{lemma}

\begin{proof}
The constants $I_{p_k}, W_{p_k}, \kappa_{p_k}$ are uniformly bounded from above and below by \eqref{eq:unifrombound2}.
Moreover, by \eqref{eq:vortexsupbound},  
\[
\liminf_{n\to\infty}E[\zeta_k]\ge I_\infty,
\quad
P[\zeta_k]=1,\quad 
\limsup_{k\to\infty}Q_{\infty}[\zeta_k] \le 1.
\]
Applying Theorem~\ref{thm:existence} yields, up to a subsequence, the convergence \eqref{eq:lar_p_vor_lim}.
Applying the energy inequality \eqref{eq:en_diff} to $E[\zeta_k]-E[\zeta_H]=I_{p_k}-I_\infty$, we obtain the convergence \eqref{eq:lar_p_en_limit}.
The convergence \eqref{eq:Wconv} and \eqref{eq:kappaconv} follow from \eqref{eq:speed}.
\end{proof}

\begin{proposition}
Under the assumptions of Lemma \ref{lem:conv_hill}, we have 
\begin{align}
\frac{1}{r^{2}}\mathcal{G}\zeta_k &\to \frac{1}{r^{2}}\mathcal{G}\zeta_H\quad \textrm{in}\ C^1(\bbR^5),\label{eq:varphiconv} \\
\zeta_k&\to \zeta_H\quad \quad \textrm{locally uniformly in}\ \mathbb{R}^{5}\setminus (\partial B_{\varrho_H}(0)\cup (B_{\varrho_H}(0)\cap \{y'=0\})).  \label{eq:strongerconvzeta}
\end{align}
Both convergences are understood under the five-dimensional transformation \eqref{eq:5dtrans}.
\end{proposition}

\begin{proof}
Setting $\varphi_k=\mathcal{G}\zeta_k/r^{2}$ and $\varphi_H=\mathcal{G}\zeta_H/r^{2}$, we have $\varphi_k-\varphi_H=N*(\zeta_k-\zeta_H)$ for $N(y)=(8\pi^{2})^{-1}|y|^{-3}$.
By \eqref{eq:lar_p_vor_lim}, we have $\zeta_k\to\zeta_H$ in $L^{1}(\bbR^5)$.
Using $\sup_{k}Q_{\infty}[\zeta_k]<\infty$, for every $\delta > 0$, we estimate 
\begin{align*}
\|\varphi_k-\varphi_H\|_{L^\infty(\mathbb{R}^{5})} &\lesssim  \sup_{y \in \bbR^5} \int_{|w-y|<\delta} \frac{|\zeta_k(w)-\zeta_H(w)|}{|w-y|^3}\,dw + \delta^{-3}\|\zeta_k-\zeta_H\|_{L^1(\bbR^5)} \\
&\lesssim \delta^2 + \delta^{-3} \|\zeta_k-\zeta_H\|_{L^1(\bbR^5)}.
\end{align*}
Using $|\nabla N(y)| \lesssim |y|^{-4}$, the same decomposition gives
\[
\|\nabla\varphi_k - \nabla\varphi_H\|_{L^\infty(\mathbb{R}^{5})} \lesssim  \delta +  \delta^{-4}\|\zeta_k-\zeta_H\|_{L^1(\bbR^5)}.
\]
Letting $k\to\infty$ and then $\delta\to0$ yields the convergence \eqref{eq:varphiconv}.

By \eqref{eq:Wconv}, we have 
\[
\Phi_k = \frac{1}{r^{2}}\mathcal{G}\zeta_k - \frac{W_{p_k}}{2}\to    \frac{1}{r^{2}}\mathcal{G}\zeta_H - \frac{W_\infty}{2}=\Phi_H\quad \textrm{in}\ C^1(\bbR^5).
\]
Since $\zeta_k=\kappa_{p_k}^{-\frac{1}{p_k-1}}(r^{2}\Phi_k)_{+}^{\frac{1}{p_k-1}}$ and $t_+^{\frac{1}{p_k-1}}$ converges locally uniformly to $\mathbf{1}_{(0,\infty)}(t)$ on $\mathbb{R}\setminus \{0\}$, the convergence \eqref{eq:strongerconvzeta} follows.
\end{proof}

\subsection{Spherical vortex rings for large \texorpdfstring{$p$}{p}}

Since the stream function of Hill's vortex is strictly decreasing in the radial direction near the vortex boundary, the strong convergence of the stream functions allows us to represent the boundaries of the vortex cores of the energy maximizers as graphs by the implicit function theorem.

\begin{proof}[Proof of Theorem \ref{t:MT3}]
It remains to show (i) and (ii).
The convergence \eqref{eq:Iinftyconv} follows from \eqref{eq:lar_p_en_limit} and \eqref{eq:Wconv}.
Let $p_k$, $\zeta_k\in S_{p_k}$, and $W_{p_k}$ be as in Lemma \ref{lem:conv_hill}.
The maximizer $\zeta_k$ converges to $\zeta_H$ locally uniformly in $\mathbb{R}^{3}\setminus (\partial B_{\varrho_H}(0)\cup (\{r=0\}\cap B_{\varrho_H}(0)))$ by \eqref{eq:strongerconvzeta}. By \eqref{eq:5Dstreamhill}, the function $\Phi_H=\mathcal{G}\zeta_H/r^{2}-W_{\infty}/2$ satisfies 
\[
\Phi_H=0, \qquad \partial_{\varrho}\Phi_H = -\frac{\varrho_H}{5} \quad\text{on } \{\varrho=\varrho_H\}.
\]
Since $\Phi_k=\mathcal{G}\zeta_k/r^{2}-W_{p_k}/2$ converges to $\Phi_H$ in $C^{1}(\mathbb{R}^{5})$, for every sufficiently small $\delta> 0$, there exists $N_{\delta}$ such that for $k\ge N_{\delta}$, we have 
\begin{align*}
\Phi_k > 0 \quad\text{on }\,\, \{\varrho\le \varrho_H - \delta\}, \quad \Phi_k < 0 \quad\text{on }\,\,  \{\varrho \ge \varrho_H + \delta\},\quad 
\partial_{\varrho} \Phi_k \le -\frac{\varrho_H}{10}\quad\text{on }\,\, \{|\varrho - \varrho_H| \le \delta\}.
\end{align*}
This means that every ray from the origin intersects $\{\Phi_k=0\}$ exactly once. By the implicit function theorem, there exists a function $\varrho_{k} \in C^1[0,\pi]$ such that 
\[
\{(z,r)\in \overline{\mathbb{R}^{2}_{+}} : \Phi_k\ge 0\} = \left\{ (\varrho\cos\vartheta,\varrho\sin\vartheta): 0 \le \varrho \le \varrho_{k}(\vartheta),\ 0 \le \vartheta \le \pi \right\}.
\]
Since $\{r^{2}\Phi_k>0\} = \{\zeta_k > 0 \}$, taking the closure yields the representation \eqref{eq:suppgraph}. The $C^1$ convergence of $\Phi_k$ to $\Phi_H$, together with the identity $\partial_{\varrho}\Phi_H = -\varrho_H/5$ on $\{\varrho = \varrho_H\}$, yields
\[
\lim_{k\to\infty}\|\varrho_{k} - \varrho_H\|_{C^1[0,\pi]}=0.
\]
Since the sequence $p_k\to\infty$ is arbitrary, this proves \eqref{eq:lar_p_C1}. The proof is complete.
\end{proof}

\section{Toroidal vortex rings and the limit \texorpdfstring{$p\to1^{+}$}{p → 1}}\label{s:7}

We prove Theorem \ref{t:MT4}. As outlined in \S\ref{s:1.6.3}, we proceed as follows.
We first establish the asymptotic formula for $I_p$ and show that the measures associated with the impulse and Casimir densities converge to the Dirac mass $\delta_{(0,e^{1/3})}$.
These results yield assertions (ii) and (iii) of Theorem \ref{t:MT4}.

Using the Euler--Lagrange equation and energy inequalities, we further show that, for $p$ sufficiently close to $1$, the support of each maximizer stays away from the symmetry axis $r=0$.
Consequently, the maximizers form toroidal vortex rings whose circulation converges to $(2\pi e^{2/3})^{-1}$ as $p\to1^+$.
This establishes assertions (i) and (iv) of Theorem \ref{t:MT4}.

\subsection{Sharp lower bound}

We begin with a sharp lower bound for the maximal energy.

\begin{lemma}[Sharp lower bound]\label{lem:p1_low}
For $R>1$, we have
\begin{align}
 \liminf_{p\to1^+}(p-1)I_p
 \ge \frac{\log R}{8\pi^2 R^3}.  \label{eq:SLBR1}
\end{align}
In particular, 
\begin{align}
\liminf_{p \to 1^{+}}(p-1)I_p\ge \frac{1}{24\pi^{2}e}.  \label{eq:p1_en_subseq1}
\end{align}
\end{lemma}

\begin{proof}
For $0<\varrho<R$, we set $\tilde{\zeta}=R^{2/(p-1)}  \mathbf{1}_{B_{\varrho}}$ and $B_{\varrho}=B_{\varrho}(0,R)$. By elementary computation,
\begin{align*} 
P[\tilde{\zeta}]=2\pi^{2}R^{\frac{2}{p-1}+3}\varrho^{2}\left(1+\frac{3}{4}\left(\frac{\varrho}{R}\right)^{2}\right),\quad Q_{p}^{p}[\tilde{\zeta}]=2\pi^{2}R^{\frac{2p}{p-1} +1 }\varrho^{2}.
\end{align*}
As a function of $0<\varrho<R$,
$P[\tilde{\zeta}]$ is strictly increasing from $0$
and exceeds $1$ at $\varrho=R$. We choose the unique $\varrho_p>0$ such that $P[\tilde{\zeta}]=1$. The constant $\varrho_p$ satisfies 
\begin{align*}
0 < \varrho^{2}_{p}\le \varrho^{2}_{p}\left(1+\frac{3}{4}\left(\frac{\varrho_{p}}{R}\right)^{2}\right) = \frac{1}{2\pi^2R^{\frac{2}{p-1}+3}} \to 0\quad \textrm{as}\ p\to 1^{+}.
\end{align*}
It follows that 
\begin{align*}
Q_p^{p}[\tilde{\zeta}]
=2\pi^{2}R^{\frac{2p}{p-1} + 1}\varrho_p^{2}
=\left(1+\frac{3}{4}\left(\frac{\varrho_p}{R}\right)^2\right)^{-1}\le 1.
\end{align*}
Using $P[\tilde{\zeta}]=1$, we have 
\begin{equation*}
    \begin{aligned}
\log{\frac{1}{\varrho_p}}&=\frac{1}{p-1}\log{R}+O(1), \\
R^{\frac{2}{p-1}}|B_{\varrho_p}|&=\pi R^{\frac{2}{p-1}}\varrho_p^{2}=\frac{1}{2\pi R^{3}}+o(1),\quad \textrm{as}\ p\to 1^{+}.        
\end{aligned}
\end{equation*}
For $X,Y\in B_{\varrho_p}$, by the asymptotics of the Green function \eqref{eq:p1_loc_ker}, 
\begin{align*}
\pi G(X,Y)X_2Y_2 = \frac{ {(X_2Y_2)^{3/2}}}{2}\log\frac{1}{|X-Y|} + O(1)\quad \textrm{as}\ p\to 1.
\end{align*}
It follows that 
\begin{align*}
E[\widetilde{\zeta}] =\frac{R^{\frac{4}{p-1}+3}}{2}\iint_{B_{\varrho_p}^{2}}\log\frac{1}{|X-Y|}\,dXdY + O(1).
\end{align*}
After the change of variables, we find
\[
\begin{aligned}
\iint_{B_{\varrho_p}^{2}}\log\frac{1}{|X-Y|} \,dXdY = \varrho_p^4\iint_{B_1(0)^{2}}\left(\log\frac{1}{\varrho_p}+\log\frac{1}{|X-Y|}\right)\,dXdY
=  |B_{\varrho_p}|^2 \log\frac{1}{\varrho_p} + O\left(\varrho_p^4\right).
\end{aligned}
\]
It follows that
\[
E[\widetilde{\zeta}] =\frac{1}{2}R^{\frac{4}{p-1}+3}|B_{\varrho_p}|^2 \log\frac{1}{\varrho_p} +O(1).
\]
Using the asymptotics for $\varrho_p$ and $|B_{\varrho_p}|$, we obtain 
\[
(p-1)I_p\ge (p-1)E[\widetilde{\zeta}] = \frac{\log{R}}{8\pi^{2}R^{3}}+o(1)\quad \textrm{as}\ p\to 1^{+}.
\]
Thus \eqref{eq:SLBR1} holds. Maximizing the right-hand side at  $R=e^{1/3}$ yields \eqref{eq:p1_en_subseq1}.
\end{proof}

\subsection{Compactness of limiting measures}

We next establish the convergence of the probability measures associated with the impulse and the Casimir.
In what follows, we choose a sequence $p_k \to 1^+$ such that $(p_k-1)I_{p_k}$ converges, and denote its limit by
\[
\mathcal{I}_*=\lim_{k\to\infty}(p_k-1)I_{p_k}.
\]
For $\mu,\nu\ge0$, we define
\[
\mathcal{I}_{*,\mu,\nu}
=\mathcal{I}_*\mu^{1/2}\nu^{3/2}.
\]
By Remark \ref{r:strictp_1}, this function satisfies the strict superadditivity
\begin{align}
\mathcal{I}_{*,\mu_1,\nu_1}+\mathcal{I}_{*,\mu_2,\nu_2}<\mathcal{I}_{*,\mu_1+\mu_2,\nu_1+\nu_2},  \label{eq:limitss}
\end{align}
for nonnegative $\mu_i,\nu_i\ge 0$ satisfying $\mu_i+\nu_i>0$ for $i=1,2$ and $(\mu_1+\mu_2)(\nu_1+\nu_2)>0$.
We first apply the concentration-compactness principle to establish the tightness of the impulse and Casimir measures on $\overline{\mathbb{R}^2_+}$, using the positivity of $\mathcal I_{*,\mu,\nu}$ for $\mu,\nu>0$ and its strict superadditivity.
We then use an energy inequality \eqref{eq:variational-inequality2} to show that the masses of these measures near the boundary are uniformly small, thereby establishing their tightness on $\mathbb{R}^2_+$.

\begin{proposition}[Compactness up to translations]\label{lem:p1_concent}
Let $p_k \to 1^+$ be a sequence along which $(p_k - 1) I_{p_k}$ converges.
Then, we have 
\begin{align}
 \mathcal{I}_*=\lim_{k\to\infty}(p_k-1)I_{p_k} 
 \ge
 \frac{1}{24\pi^{2}e}.  \label{eq:p1_en_subseq}
\end{align}
For $\zeta_k\in S_{p_k}$, define
\begin{align}
 d\mathfrak{p}_k = 2\pi r^3\zeta_k\,dz\,dr,
 \qquad
 d\mathfrak{q}_k = 2\pi r\zeta_k^{p_k}\,dz\,dr. \label{eq:meas_cons}
\end{align}
Then both $\mathfrak{p}_k$ and $\mathfrak{q}_k$ are probability measures.
For every $\varepsilon>0$, there exists $\varrho>0$ such that, up to translations in the $z$-direction, we have 
\begin{align}
\liminf_{k\to\infty}\mathfrak{p}_k \left(B_{\varrho}(0)\cap\bbR^2_+\right) \ge 1 - \varepsilon,
\quad
\liminf_{k\to\infty}\mathfrak{q}_k \left( B_{\varrho}(0)\cap\bbR^2_+\right) \ge 1 - \varepsilon.  \label{eq:Compactness up to translations}
\end{align}
\end{proposition}

\begin{proof}
The measures $\mathfrak{p}_k$ and $\mathfrak{q}_k$ are probability measures by $P[\zeta_k]=Q_{p_k}[\zeta_k]=1$. We define the averaged probability measure by 
\begin{align}
\mathfrak{r}_k=\frac{1}{2}(\mathfrak{p}_k+\mathfrak{q}_k).\label{eq:averaged}
\end{align}
We apply the concentration-compactness principle (Lemma~\ref{lem:lions}) to exclude the possibility of (ii) vanishing and (iii) dichotomy of the sequence, and deduce (i) compactness, which yields \eqref{eq:Compactness up to translations}. \\

\noindent \textit{(ii) Vanishing.} Suppose that vanishing occurs. Namely, 
\begin{align*}
\delta_k = \sup_{X \in \bbR^2_+} \mathfrak{r}_{k}(B_2(X) \cap \bbR^2_+)=\frac{1}{2}\sup_{X\in \mathbb{R}^{2}_{+}}\left(P[\zeta_k\mathbf{1}_{B_2(X)}]+Q_{p_{k}}^{p_k}[\zeta_k \mathbf{1}_{B_2(X)}] \right)\to 0.
\end{align*}
For $X_0 \in \bbR^2_+$, we apply \eqref{eq:variational-inequality2} to estimate 
\begin{align*}
E[\zeta_k \mathbf{1}_{B_2(X_0)}] \le I_{p_k} P[\zeta_k \mathbf{1}_{B_2(X_0)}]^{\alpha_{k}} Q_{p_k}[\zeta_k \mathbf{1}_{B_2(X_0)}]^{\beta_{k}}, 
\end{align*}
with the constants $\alpha_k$ and $\beta_k$ defined by \eqref{eq:exponent}. Using $\alpha_{k} + \frac{\beta_{k}}{p_k} - 1 = \frac{2p_k}{5p_k - 3}$, we have
\[
E[\zeta_k \mathbf{1}_{B_2(X_0)}] \lesssim I_{p_k} \delta_k^{\frac{2p_k}{5p_k - 3}} \left(P[\zeta_k \mathbf{1}_{B_2(X_0)}] + Q_{p_k}^{p_k}[\zeta_k \mathbf{1}_{B_2(X_0)} ] \right).
\]
For $X,Y\in \mathbb{R}^{2}_{+}$ satisfying $|X-Y|<1$,
\begin{align*}
\int_{\mathbb{R}^{2}_{+}}\mathbf{1}_{\{|X - X_0| < 2\}} \mathbf{1}_{\{|Y - X_0| < 2\}} dX_0=|B_2(X)\cap B_2(Y) \cap \mathbb{R}^{2}_{+}|\gtrsim \mathbf{1}_{\{|X- Y| < 1\}}.
\end{align*}
Thus, the near-field energy satisfies
\[
 \iint_{ \{ |X - Y| < 1 \}} G(X, Y) \zeta_k(X) \zeta_k(Y) X_2 Y_2\,dXdY \lesssim \int_{\bbR^2_+} E[\zeta_k \mathbf{1}_{|X - X_0| < 2}] dX_0  \lesssim I_{p_k} \delta_k^{\frac{2p_k}{5p_k - 3}}.
\]
On the other hand, the far-field energy contribution is uniformly bounded by \eqref{eq:Green_bound} for $\tau =  3/2$, and hence 
\begin{align*}
0<\mathcal{I}_*=\lim_{k\to \infty}(p_k - 1)I_{p_k}=\lim_{k\to \infty}(p_k - 1)E[\zeta_k]= 0.
\end{align*}
This is a contradiction. Thus vanishing does not occur.\\

\noindent \textit{(iii) Dichotomy.} 
Suppose that dichotomy occurs.
Then, after passing to a subsequence, there exist
$\tau\in(0,1)$ and two measurable sets $A_k,B_k\subset\bbR^2_+$
such that
\begin{align*}
 \dist(A_k,B_k)\to\infty,&
 \qquad
 \mathfrak{r}_{k}(\bbR^2_+\setminus(A_k\cup B_k))\to0,\\
 \mathfrak{r}_{k}(A_k)\to\tau,&
 \qquad
 \mathfrak{r}_{k}(B_k)\to1-\tau.
\end{align*}
We set 
\[
 f_k=\zeta_k\mathbf 1_{A_k},
 \qquad
 g_k=\zeta_k\mathbf 1_{B_k}.
\]
Using $0\le P[f_k], P[g_k]\le 1$ and $0\le Q_{p_k}^{p_k}[f_k],Q_{p_k}^{p_k}[g_k]\le 1$, and passing to a further subsequence, we may assume that
\begin{align*}
 P[f_k]\to\mu_1&,\qquad
 P[g_k]\to\mu_2, \\
Q_{p_k}^{p_k}[f_k]\to\nu_1&,\qquad
 Q_{p_k}^{p_k}[g_k]\to\nu_2,
\end{align*}
for some constants $0\le \mu_i,\nu_i\le 1$. Since
\begin{align*}
 2\mathfrak{r}_{k}(A_k)
 =
 P[f_k]+Q_{p_k}^{p_k}[f_k],\quad 
 2\mathfrak{r}_{k}(B_k)
 =
 P[g_k]+Q_{p_k}^{p_k}[g_k],
\end{align*}
by the definition of the averaged measure, we have
\[
 \mu_1+\nu_1=2\tau,
 \qquad
 \mu_2+\nu_2=2(1-\tau).
\]
Moreover, using $P[\zeta_k]=1$ and $Q_{p_k}^{p_k}[\zeta_k]=1$, we have $\mu_1+\mu_2\le 1$ and $\nu_1+\nu_2\le 1$. Since $\mu_1+\mu_2+\nu_1+\nu_2=2$, it follows that 
\begin{align*}
\mu_1+\mu_2= 1,\quad \nu_1+\nu_2= 1.
\end{align*}
We set 
\[
\zeta_k=f_k+g_k+h_k.
\]
The energy decomposition yields 
\[
\mathcal{I}_*
=
\lim_{k\to \infty}(p_k-1)I_{p_k}
=
\lim_{k\to \infty}(p_k-1)E[\zeta_k]
\le
\limsup_{k\to \infty}(p_k-1)(E[f_k]+E[g_k]).
\]
Applying \eqref{eq:variational-inequality2}, we have 
\begin{align*}
\limsup_{k\to \infty}(p_k-1)E[f_k]\le \mathcal{I}_{*,\mu_1,\nu_1},\quad
\limsup_{k\to \infty}(p_k-1)E[g_k]\le \mathcal{I}_{*,\mu_2,\nu_2}.
\end{align*}
The strict superadditivity \eqref{eq:limitss} yields 
\begin{align*}
\mathcal{I}_* \le \mathcal{I}_{*,\mu_1,\nu_1}+\mathcal{I}_{*,\mu_2,\nu_2}<\mathcal{I}_{*,\mu_1+\mu_2,\nu_1+\nu_2}= \mathcal{I}_*.
\end{align*}
We obtain a contradiction. Thus, dichotomy does not occur.\\

\noindent \textit{(i) Compactness.} Therefore, there exist $X_k=(z_k,r_k)\in\overline{\bbR^2_+}$ such that for every $\varepsilon>0$,
there exists $\varrho>0$ satisfying
\[
 \liminf_{k\to\infty}
 \mathfrak{r}_k\bigl(B_\varrho(X_k)\cap\bbR^2_+\bigr)
 \ge 1-\varepsilon.
\]
By \eqref{eq:p1_en_subseq} and \eqref{eq:speed}, we have  $\liminf_{k\to\infty}(p_k-1)W_{p_k}>0$. Using \eqref{eq:str_2} and \eqref{eq:EL_intro0}, we find that for $X=(z,r)\in \supp \zeta_k$,  
\begin{equation*}
    r^2 \le \frac{2}{W_{p_k}}\|\mathcal{G}\zeta_k\|_{L^{\infty}(\mathbb{R}^{2}_{+})} \le C.
\end{equation*}
Hence $\{r_k\}$ is bounded.
By enlarging the radius and translating the measure, we may assume that $X_k=0$. Since $\mathfrak{p}_k$ and $\mathfrak{q}_k$ are nonnegative measures, we have
\[
 \mathfrak{p}_k\le 2\mathfrak{r}_k,
 \qquad
 \mathfrak{q}_k\le 2\mathfrak{r}_k.
\]
Hence
\[
 \mathfrak{p}_k\bigl(B_R(0)^c\cap\bbR^2_+\bigr)
 \le
 2\mathfrak{r}_k\bigl(B_R(0)^c\cap\bbR^2_+\bigr),
\]
and the same estimate holds for $\mathfrak{q}_k$.  Since both $\mathfrak{p}_k$
and $\mathfrak{q}_k$ are probability measures, by replacing
$\varepsilon$ by $\varepsilon/2$, \eqref{eq:Compactness up to translations} follows.
\end{proof}

\begin{proposition}\label{lem:rad_en}
There exists a constant $C>0$ such that, for every $1<p\le\infty$, every $Z,R>0$, and every $\zeta$ satisfying $Q_p[\zeta] \le 1$, we have
\begin{align}
E[\zeta\mathbf{1}_{\{|X_1|\le Z,\ 0<X_2<R\}}]\le \frac{p}{p-1}C|Z|^{\alpha_p\left(1-\frac{1}{p}\right)} R^{2\alpha_p\left(2-\frac{1}{p}\right)},  \label{eq:nearboundary}
\end{align}
where $\alpha_p$ is given by \eqref{eq:exponent}.
\end{proposition}

\begin{proof}
By \eqref{eq:variational-inequality2} and \eqref{eq:unifrombound2}, we estimate
\begin{align*}
E[\zeta\mathbf{1}_{\{|X_1|\le Z,\ 0<X_2<R\}}]\le I_pP[\zeta\mathbf{1}_{\{|X_1|\le Z,\ 0<X_2<R\}}]^{\alpha_p}Q_p[\zeta]^{\beta_p}
\le C\frac{p}{p-1}P[\zeta\mathbf{1}_{\{|X_1|\le Z,\ 0<X_2<R\}}]^{\alpha_p}.
\end{align*}
By H\"older's inequality, we estimate 
\[
P[\zeta\mathbf{1}_{\{|X_1|\le Z,\ 0<X_2<R\}}] 
\le R^2 Q_p[\zeta]Q_{p'}[\mathbf{1}_{\{|X_1|\le Z,\ 0<X_2<R\}}]
\lesssim |Z|^{1-\frac{1}{p}} R^{4-\frac2p}.
\]
Substituting this estimate, we obtain \eqref{eq:nearboundary}.
\end{proof}

We denote by $C_b(\mathbb{R}^{2}_{+})$ the space of bounded continuous functions on $\mathbb{R}^{2}_{+}$. We say that $\mathfrak{m}_k\in\mathcal{P}(\mathbb{R}^{2}_{+})$
converges weakly to $\mathfrak{m}\in\mathcal{P}(\mathbb{R}^{2}_{+})$
if
\[
\int_{\mathbb{R}^{2}_{+}} \phi\,d\mathfrak{m}_k
\to
\int_{\mathbb{R}^{2}_{+}} \phi\,d\mathfrak{m}
\quad \text{as } k\to\infty
\]
for every $\phi\in C_b(\mathbb{R}^{2}_{+})$.
We write
$\mathfrak{m}_k\rightharpoonup \mathfrak{m}$ in $\mathcal{P}(\mathbb{R}^{2}_{+})$.

\begin{lemma}[Compactness of limiting measures]\label{lem:meas_tight}
There exist probability measures $\mathfrak{p}$ and $\mathfrak{q}$ such that, after passing to a subsequence, we have 
\begin{align}
 \mathfrak{p}_k \rightharpoonup \mathfrak{p},
 \qquad
 \mathfrak{q}_k \rightharpoonup \mathfrak{q}
 \quad\text{in }\mathcal P(\mathbb{R}^2_+).  \label{eq:meas_lim}
\end{align}
\end{lemma}

\begin{proof}
By Proposition \ref{lem:p1_concent}, for arbitrary $\varepsilon>0$, there exist $Z>0$ and $R>0$ such that 
\begin{align*}
\limsup_{k\to\infty}Q_{p_k}[\zeta_k\mathbf{1}_{\{|X_1|<Z,X_2<R\}^{c}}]\le \varepsilon.
\end{align*}
For $0<\delta<R$, using 
\begin{align*}
\{|X_1|<Z,X_2<R\}^{c}\cup \{|X_1|<Z,0<X_2<\delta\}=\{|X_1|<Z,\delta<X_2<R\}^{c}, 
\end{align*}
we decompose the $L^{p}$-Casimir as 
\begin{align*}
Q_{p_k}^{p_k}[\zeta_k\mathbf{1}_{\{|X_1|<Z,X_2<R\}^{c}}]+Q_{p_k}^{p_k}[\zeta_k\mathbf{1}_{\{|X_1|<Z,0<X_2<\delta\}}]
=Q_{p_k}^{p_k}[\zeta_k\mathbf{1}_{\{|X_1|<Z,\delta<X_2<R\}^{c}}]\,.
\end{align*}
It follows that 
\begin{align*}
\limsup_{k\to\infty}Q_{p_k}^{p_k}[\zeta_k\mathbf{1}_{\{|X_1|<Z,\delta<X_2<R\}^{c}}]\le \varepsilon+\limsup_{k\to\infty}Q_{p_k}^{p_k}[\zeta_k\mathbf{1}_{\{|X_1|<Z,0<X_2<\delta\}}]\,.
\end{align*}
We show that the last term is bounded by $\varepsilon$ for sufficiently small $\delta>0$. We set 
\begin{align*}
\zeta_k=\zeta_k\mathbf{1}_{\{|X_1|\le Z,\ 0<X_2<\delta\}}+\zeta_k(1-\mathbf{1}_{\{|X_1|\le Z,\ 0<X_2<\delta\}}).
\end{align*}
We apply \eqref{eq:nearboundary} to estimate 
\begin{align*}
(p_k-1) E[\zeta_k\mathbf{1}_{\{|X_1|\le Z,\ 0<X_2<\delta\}}]\le p_kC|Z|^{\alpha_{p_k}\left(1-\frac{1}{p_k}\right)} \delta^{2\alpha_{p_k}\left(2-\frac{1}{p_k}\right)}\,. 
\end{align*}
Using 
\begin{align*}
P[\zeta_k(1-\mathbf{1}_{\{|X_1|\le Z,\ 0<X_2<\delta \}})]&\le 1,\\
Q_{p_k}^{p_k}[\zeta_k(1-\mathbf{1}_{\{|X_1|\le Z,\ 0<X_2<\delta\}})]&=1-Q_{p_k}^{p_k}[\zeta_k\mathbf{1}_{\{|X_1|\le Z,\ 0<X_2<\delta\}}],
\end{align*}
we apply \eqref{eq:variational-inequality2} to estimate 
\begin{align*}
(p_k-1)E[\zeta_k(1-\mathbf{1}_{\{|X_1|\le Z,\ 0<X_2<\delta\}})]\le (p_k-1)I_{p_k}(1-Q_{p_k}^{p_k}[\zeta_k\mathbf{1}_{\{|X_1|\le Z,\ 0<X_2<\delta\}}])^{\beta_{p_k}}.
\end{align*}
We apply the triangle inequality \eqref{eq:p1_en_tri} for $E[\zeta_k]$ and multiply it by $\sqrt{p_k-1}$.
It follows that 
\begin{align*}
\sqrt{(p_k-1)I_{p_k}}&=\sqrt{(p_k-1)E[\zeta_k]} \\
&\le \sqrt{p_kC|Z|^{\alpha_{p_k}\left(1-\frac{1}{p_k}\right)} \delta^{2\alpha_{p_k}\left(2-\frac{1}{p_k}\right)}}+\sqrt{(p_k-1)I_{p_k}}(1-Q_{p_k}^{p_k}[\zeta_k\mathbf{1}_{\{|X_1|\le Z,\ 0<X_2<\delta\}}])^{\beta_{p_k}/(2p_k)}.
\end{align*}
Taking the liminf as $k\to\infty$ and then taking the liminf as $\delta\to 0$ yields
\[
\sqrt{\mathcal{I}_*}\le\sqrt{\mathcal{I}_*}(1-\limsup_{\delta \to 0}\limsup_{k\to\infty}Q_{p_k}[\zeta_k\mathbf{1}_{\{|X_1|\le Z,\ 0<X_2<\delta\}}])^{3/4}\le \sqrt{\mathcal{I}_*}.
\]
Thus, for sufficiently small $\delta>0$, we have  
\begin{align*}
\limsup_{k\to\infty}Q_{p_k}[\zeta_k\mathbf{1}_{\{|X_1|\le Z,\ 0<X_2<\delta\}}]\le \varepsilon.
\end{align*}
We obtain 
\begin{align*}
\limsup_{k\to\infty}Q_{p_k}^{p_k}[\zeta_k\mathbf{1}_{\{|X_1|<Z,\ \delta<X_2<R\}^{c}}]\le 2\varepsilon.
\end{align*}
Therefore, $\{\mathfrak{q}_k\}$ is tight, and its weak convergence along a subsequence follows from Prokhorov's theorem.
The same argument applies to $\{\mathfrak{p}_k\}$. 
\end{proof}

\subsection{Atomic estimate of the limiting energy}

We call $X_0=(z_0,r_0)$ an atom of the limit measure $\mathfrak{q}$ if $\mathfrak{q}(\{X_0\}) > 0$.
For the atom $X_0$, we refer to  
\begin{align*}
 \mathsf p_0=\mathfrak{p}(\{X_0\}),
 \qquad
 \mathsf q_0=\mathfrak{q}(\{X_0\})
\end{align*}
as its impulse weight and Casimir weight, respectively.
We denote the atoms of the limiting measure $\mathfrak q$ by $\{X_j=(z_j,r_j)\}_{j\ge1}$.
We show that only these atoms contribute to the limit of the rescaled energy $(p_k-1)E[\zeta_k]$.

\begin{figure}[htbp]
\centering
\includegraphics[width=0.8\textwidth]{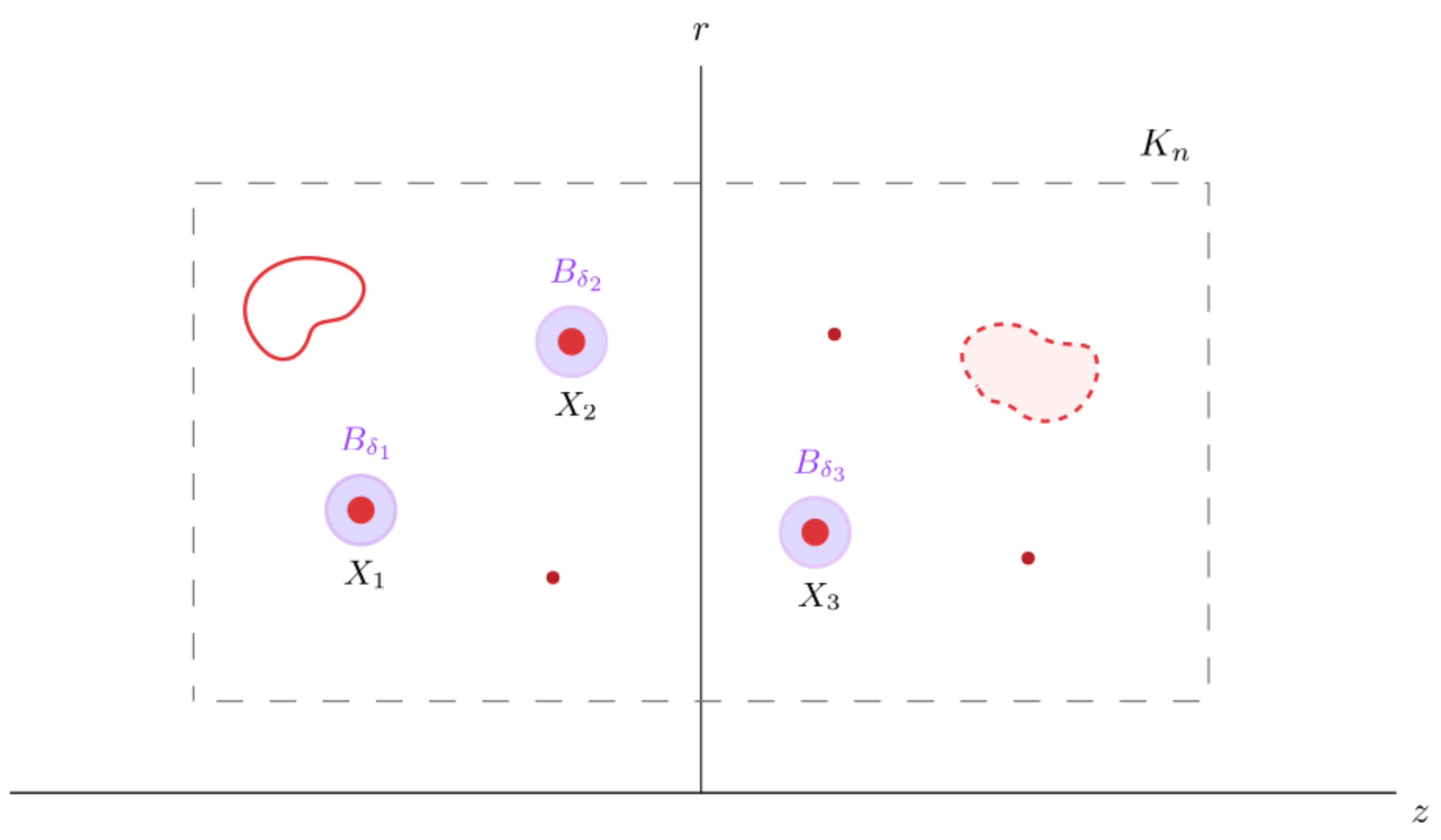}
\caption{The limiting measure $\mathfrak{q}$ obtained from the weak convergence of the probability measures $\mathfrak{q}_k$ may, in general, have atomic ({\textcolor{red}{$\bullet$}}), absolutely continuous (\raisebox{-0.3ex}{\redregion}), and singular non-atomic components ({\textcolor{red}{$\sim$}}).
We approximate $\mathbb{R}^2_+$ by bounded domains $K_n$ and cover each atom $X_j$ in $K_n$ with mass greater than $\tau$ by a disk $B_{\delta_j}(X_j)$.
Based on this decomposition, we split the sequence $\zeta_k$ into components supported on the corresponding regions.}
\label{fig:atomic}
\end{figure}

To this end, we consider an increasing sequence of rectangular domains $K_n\Subset\mathbb R^2_+$ such that
\[
\varepsilon_n
=\mathfrak p(\mathbb R^2_+\setminus K_n)
+\mathfrak q(\mathbb R^2_+\setminus K_n)
\longrightarrow0,
\]
and the atoms in $K_n$ whose Casimir weights
$q_j$ exceed $\tau>0$. Namely, we set 
\begin{align*}
J_{n,\tau}
=
\{j\in \mathbb{N}\ :\ X_j\in K_n,\ \mathsf{q}_j>\tau\}.  
\end{align*}
We isolate the portions of $\zeta_k$ associated with these large atoms by choosing disks $B_{\delta_j}(X_j)$ around each $X_j$ and defining
\[
\zeta_k^{\mathrm{atom}}
=\sum_{j\in J_{n,\tau}}
\zeta_k\mathbf 1_{B_{\delta_j}(X_j)}.
\]
We denote the remaining part of $\zeta_k$ inside $K_n$ by $\zeta_k^{\mathrm{rem}}$, and the part outside $K_n$ by $\zeta_k^{\mathrm{far}}$.
This yields the decomposition
\[
\zeta_k
=\zeta_k^{\mathrm{atom}}
+\zeta_k^{\mathrm{rem}}
+\zeta_k^{\mathrm{far}}.
\]
In general, the limiting measure $\mathfrak{q}$ may contain multiple atoms, an absolutely continuous part with respect to the Lebesgue measure, and a singular continuous part.
These possible components of $\mathfrak{q}$ are illustrated schematically in Figure \ref{fig:atomic}.

For the atomic part $\zeta_k^{\mathrm{atom}}$, we derive a sharp upper bound for its energy using the following logarithmic-type HLS inequality, exploiting the logarithmic singularity of the leading term of the Green function.
On the other hand, the energy contributions from the remainder and far-field parts are controlled in terms of $\tau$ and $\varepsilon_n$.

\begin{proposition}[Logarithmic-type HLS inequality]\label{lem:p1_log_HLS}
The inequality
\begin{align}
\varepsilon \iint_{\bbR^2 \times \bbR^2 } \log\frac{D}{|X-Y|} f(X)f(Y)\,dXdY \le \frac{M^2}{2} \log\left( D^{2\varepsilon}C_\varepsilon \frac{N}{M^{1+\varepsilon}} \right),  \label{eq:logHLS}
\end{align}
holds for $D>0$, $0<\varepsilon<1$, and $f\ge0$ with the constants 
\begin{align}
M = \int_{\bbR^2}f\,dX, \quad N = \int_{\bbR^2}f^{1+\varepsilon}\,dX,\quad C_\varepsilon = \frac{\pi^\varepsilon}{1-\varepsilon}.
\end{align}
\end{proposition}

\begin{proof}
We apply Jensen's inequality 
\begin{align*}
\exp\left(\iint_{\mathbb{R}^{2} \times \bbR^2}g(X,Y)d\mathfrak {m}(X)d\mathfrak {m}(Y)\right)\le\iint_{\mathbb{R}^{2} \times \bbR^2}e^{g(X,Y)}d\mathfrak {m}(X)d\mathfrak {m}(Y)
\end{align*}
for the function and the probability measure
\begin{align*}
g(X,Y)=2\varepsilon \log{\frac{D}{|X-Y|}},\quad d\mathfrak {m}(X)=\frac{f(X)}{M}dX.
\end{align*}
We then obtain 
\begin{align*}
\exp\left(\frac{2\varepsilon}{M^{2}}\iint_{\mathbb{R}^{2}\times \bbR^2}\log{\frac{D}{|X-Y|}}f(X)f(Y)dXdY\right)\le \frac{D^{2\varepsilon}}{M^{2}}\iint_{\mathbb{R}^{2} \times \bbR^2}\frac{1}{|X-Y|^{2\varepsilon}}f(X)f(Y)dXdY.
\end{align*}  
By the HLS inequality \cite[Theorem 4.3]{LL},
\begin{align*}
\iint_{\bbR^2 \times \bbR^2} \frac{f(X)f(Y)}{|X-Y|^{2\varepsilon}}\,dXdY \le C_\varepsilon \|f\|_{L^{\frac{2}{2-\varepsilon}}(\bbR^2) }^2.
\end{align*}
Interpolation yields $\|f\|_{L^{\frac{2}{2-\varepsilon}}(\bbR^2)}^2 \le M^{1-\varepsilon}N$.
This proves \eqref{eq:logHLS}.
\end{proof}

\begin{proposition}[Limiting energy estimate near a single atom]\label{lem:atom_loc}
\noindent
(i) For every atom $X_0=(z_0,r_0) \in \bbR^2_+$, there exist radii $\delta_n$ satisfying 
\begin{align}
\delta_n \to0,\quad \mathfrak{p}(\partial B_{\delta_n}(X_0))&=\mathfrak{q}(\partial B_{\delta_n}(X_0))=0.  \label{eq:conditiondelta}
\end{align}

\noindent
(ii) For every atom $X_0=(z_0,r_0)\in \bbR^2_+$, the constants $\mathsf p_0=\mathfrak{p}(\{X_0\})$ and $\mathsf q_0=\mathfrak{q}(\{X_0\})$ satisfy 
\begin{align}
\mathsf p_0\le r^2_0\mathsf{q}_0.  \label{eq:mbR}
\end{align}

\noindent 
(iii) We define the atomic energy by 
\begin{equation} \label{eq:En_def}
e_{r_0}(\mathsf{p}_0,\mathsf{q}_0)=
\begin{cases}
\ \displaystyle \frac{\mathsf{p}^2_0}{16\pi^2r^3_{0}}\log\frac{r^2_{0}\mathsf{q}_0}{\mathsf{p}_0}, & \mathsf{p}_0>0,\\[8pt]
\ 0, &\mathsf{p}_0=0.
\end{cases}
\end{equation}
For every atom $X_0=(z_0,r_0)\in \bbR^2_+$ and radii $\delta=\delta_n$ satisfying \eqref{eq:conditiondelta}, we have 
\begin{align}
\lim_{\delta \to 0} \limsup_{k\to\infty}(p_k-1)E[\zeta_k\mathbf1_{B_\delta(X_0)}]\le e_{r_0}(\mathsf{p}_0,\mathsf{q}_0).  \label{eq:atombound}
\end{align}
\end{proposition}

\begin{proof}
The set of $\delta>0$ such that $\mathfrak{p}(\partial B_\delta)>0$ for $B_{\delta}=B_\delta(X_0)$ is at most countable, since the circles $\partial B_\delta$ are pairwise disjoint and $\mathfrak{p}$ is finite.
The same holds for $\mathfrak{q}$. Hence, there exist radii satisfying \eqref{eq:conditiondelta}. 

Since $r\le r_0+\delta$ on $B_\delta$, H\"older's inequality yields
\begin{align*}
\mathfrak{p}_k(B_\delta)
=P[\zeta_k\mathbf{1}_{B_\delta}]
\le 
(r_0+\delta)^2Q_{1}[\zeta_k\mathbf{1}_{B_\delta}]\le (r_0+\delta)^2Q_{p_k}[\zeta_k\mathbf{1}_{B_\delta}]Q_{p_k'}[\mathbf{1}_{B_\delta}] 
=(r_0+\delta)^2Q_{p_k'}[\mathbf{1}_{B_\delta}]\mathfrak{q}_k(B_\delta)^{1/p_k}.
\end{align*}
We take $\delta=\delta_n$ satisfying \eqref{eq:conditiondelta}. By Portmanteau's theorem, the weak convergence \eqref{eq:meas_lim} and the boundary condition \eqref{eq:conditiondelta} imply that 
\begin{align*}
\lim_{k\to\infty}\mathfrak{p}_k(B_\delta)=\mathfrak{p}(B_\delta),\quad 
\lim_{k\to\infty}\mathfrak{q}_k(B_\delta)=\mathfrak{q}(B_\delta).
\end{align*}
Letting $k \to \infty$ and then $\delta_n\to 0$ yields \eqref{eq:mbR}.

Let $0<\delta_0<r_0$. By \eqref{eq:p1_loc_ker}, there exists a constant $D_{r_0}>0$ such that for $X,Y\in B_\delta=B_\delta(X_0)$ and $0<\delta<\delta_0$, 
\begin{equation*}
\pi G(X,Y)X_2Y_2\le A_{r_0,\delta} \log\frac{D_{r_0}}{|X-Y|}, \quad A_{r_0,\delta}=\frac{(r_0+\delta)^3}{2}.
\end{equation*}
We set $f_k=\zeta_k\mathbf1_{B_\delta}$ and observe that 
\begin{align*}
E[f_k]\le A_{r_0,\delta}\iint_{\mathbb{R}^{2} \times \bbR^2 }\log\frac{D_{r_0}}{|X-Y|}f_k(X)f_k(Y)dXdY.
\end{align*}
We apply the logarithmic-type HLS inequality \eqref{eq:logHLS} for $\varepsilon=\varepsilon_k=p_{k}-1$, $D=D_{r_0}$, and $f=f_{k}$ to estimate 
\begin{align*}
\varepsilon_k E[f_k]\le\frac{A_{r_0,\delta} M^2_k}{2}\log\left(D_{r_0}^{2\varepsilon_k}C_{\varepsilon_k}\frac{N_k}{M^{1+\varepsilon_k}_k}\right),\quad M_k=\int_{\mathbb{R}^{2}} f_k\,dX, \quad N_k=\int_{\mathbb{R}^{2}} f_k^{p_k}\,dX.
\end{align*}
The constants $M_k$ and $N_k$ satisfy 
\begin{equation*}
\begin{aligned}
\frac{\mathfrak{p}_k(B_\delta)}{2\pi(r_0+\delta)^3}\le M_k \le\frac{\mathfrak{p}_k(B_\delta)}{2\pi(r_0-\delta)^3}, \qquad
\frac{\mathfrak{q}_k(B_\delta)}{2\pi(r_0+\delta)}\le N_k\le\frac{\mathfrak{q}_k(B_\delta)}{2\pi(r_0-\delta)}.
\end{aligned}
\end{equation*}
Passing first to the limit as $k\to\infty$ and then as $\delta=\delta_n\to 0$, the preceding bounds yield the limiting values
\begin{equation*}
M=\frac{\mathsf{p}_0}{2\pi r_0^3}, \quad N=\frac{\mathsf{q}_0}{2\pi r_0}.
\end{equation*}
For $\mathsf{p}_0>0$, we have $\mathfrak{p}(B_\delta)>0$ for sufficiently small $\delta$. We obtain  
\[
\lim_{\delta\to0}\limsup_{k\to\infty}\varepsilon_k E[f_k]\le\frac{r_0^{3} M^2}{ 4}\log\frac{N}{M}=\frac{\mathsf{p}^2_0}{16\pi^2r_0^3}\log\frac{r_0^2\mathsf{q}_0}{\mathsf{p}_0}=e_{r_0}(\mathsf{p}_0,\mathsf{q}_0).
\]
For $\mathsf{p}_0=0$, applying the energy estimate \eqref{eq:variational-inequality2} implies
\[
\varepsilon_k E[f_k]\le\varepsilon_k I_{p_k}\mathfrak{p}_k(B_\delta)^{\alpha_{k}}.
\]
It follows from \eqref{eq:p1_en_subseq} that 
\[
\limsup_{k\to\infty}\varepsilon_k E[f_k]\le \mathcal{I}_* \mathfrak{p}(B_\delta)^{1/2}.
\]
The right-hand side vanishes as $\delta\to0$. We obtain \eqref{eq:atombound}.
\end{proof}

We denote by $\mathcal{M}_{+}(K)$ the space of finite nonnegative Borel measures on $K$.

\begin{proposition}[Limiting energy estimate for small atoms]\label{lem:small_atom}
Let $\tau>0$ and $K\Subset\bbR^2_+$. Let $g_k\ge0$ and $\mathfrak{m}_k=2\pi r g_k^{p_k}\,dzdr$. Assume that
\begin{align}
\sup_{k\ge 1}Q_{p_k}[g_k]&<\infty,  \label{eq:SAEE1}\\
\supp g_k&\subset K,\quad k\ge 1,\label{eq:SAEE2} \\
\mathfrak{m}_k&\rightharpoonup\mathfrak{m}\quad \textrm{in}\ \mathcal{M}_{+}(K),  \label{eq:SAEE3}\\
\sup_{X\in K}\mathfrak{m}(\{X\})&\le \tau. \label{eq:SAEE4}
\end{align}
Then, we have
\begin{align}
\limsup_{k\to\infty}(p_k-1)E[g_k]\le C_K\tau,  \label{eq:SAEE5}
\end{align}
where $C_K$ is independent of $\tau$.
\end{proposition}

\begin{proof}
By the finiteness of the measure $\mathfrak{m}$ and the compactness of $K$, we have 
\[
\lim_{\delta \to 0}\sup_{X\in K} \mathfrak{m}(B_\delta(X)) = \sup_{X\in K}\mathfrak{m}(\{X\})\le \tau.
\]
For arbitrary $\varepsilon>0$, we take $\varrho>0$ such that 
\begin{align*}
\sup_{X\in K}\mathfrak{m}(B_\varrho(X))\le \tau+\varepsilon.
\end{align*}
We cover $K$ by finitely many squares $\{K_j\}_{1\le j\le N_{\delta}}$ of side length $4\delta$, with bounded overlap and $\mathfrak{m}(\partial K_j)=0$, such that 
\begin{align*}
\{(X,Y)\in K\times K: |X-Y|<\delta \}\subset \bigcup_{j=1}^{N_{\delta}}K_j\times K_j.
\end{align*}
This condition can be written as 
\begin{align*}
\mathbf{1}_{\{|X-Y|<\delta\}}\le \sum_{j=1}^{N_{\delta}}\mathbf{1}_{K_j}(X)\mathbf{1}_{K_j}(Y),\quad X,Y\in K.
\end{align*}
We choose $\delta>0$ such that $4\sqrt{2}\delta<\varrho$.
Then, for $X\in K_j\cap K$ and for $1\le j\le N_{\delta}$, it follows that  $K_j\cap K\subset B_{\varrho}(X)$.
This yields  
\[
\mathfrak{m}(K_j)=\mathfrak{m}(K_j\cap K) \le \tau + \varepsilon.
\]
The Portmanteau theorem yields
\begin{align*}
\lim_{k\to\infty}\mathfrak{m}_k(K_j)=\mathfrak{m}(K_j)\le \tau+\varepsilon.
\end{align*}
Since there are only finitely many $j$, it follows that 
\begin{align*}
\limsup_{k\to\infty}\left(\sup_{1\le j\le N_{\delta}}\mathfrak{m}_k(K_j) \right)\le \tau+\varepsilon.
\end{align*}
We multiply
\begin{align*}
\mathbf{1}_{\{|X-Y|<\delta\}}g_k(X)g_k(Y)\le \sum_{j=1}^{N_{\delta}}g_{k,j}(X)g_{k,j}(Y),\quad g_{k,j}(X)=g_k(X)\mathbf1_{K_j}(X)
\end{align*}
by $\pi G(X,Y)X_2Y_2$ and integrate it over $\mathbb{R}^{2}_+$.
It follows from \eqref{eq:variational-inequality2} that  
\[
\pi\iint_{\{|X-Y|<\delta\}}G(X,Y)g_k(X)g_k(Y)X_2Y_2\,dXdY\le\sum_{1\le j \le N_{\delta}}E[g_{k,j}]\le I_{p_k}\sum_{1\le j \le N_{\delta}}P[g_{k, j}]^{\alpha_{k}}  Q_{p_k}[ g_{k, j} ]^{\beta_{k}}.
\]
By H\"older's inequality, we have $P[g_{k,j}]\le C_K Q_{p_k}[g_{k,j}]=C_K\mathfrak{m}_k(K_j)^{1/p_k}$.
Using $\alpha_k+\beta_k=2$, we estimate 
\begin{align*}
\sum_{1\le j \le N_{\delta}}P[g_{k, j}]^{\alpha_{k}}  Q_{p_k}[ g_{k, j} ]^{\beta_{k}}
\le C_{K}^{\alpha_k}\sum_{1\le j \le N_{\delta}}  \mathfrak{m}_k(K_j)^{2/p_k}
\lesssim C_{K}^{\alpha_k} \sup_{1\le j\le N_{\delta}}\mathfrak{m}_k(K_j)^{2/p_k-1}
  \mathfrak{m}_k(K).
\end{align*}
We thus obtain 
\[
\limsup_{k\to\infty}\ (p_k-1)\pi\iint_{\{|X-Y|<\delta\}}G(X,Y)g_k(X)g_k(Y)X_2Y_2\,dXdY\le C_K'(\tau+\varepsilon).
\]
Since $\pi G(X,Y)X_2Y_2$ is bounded on $\{(X,Y)\in K\times K:|X-Y| \ge \delta\}$, we apply \eqref{eq:SAEE1} to obtain 
\[
\lim_{k\to\infty}\ (p_k-1)\pi \iint_{\{|X-Y|\ge \delta\}}G(X,Y)g_k(X)g_k(Y)X_2Y_2\,dXdY = 0.
\]
Thus,
\[
\limsup_{k\to\infty}(p_k-1)E[g_k]\le C_K'(\tau+\varepsilon).
\]
Letting $\varepsilon\to0$ yields \eqref{eq:SAEE5}.
\end{proof}

\begin{lemma}[Atomic estimate of the limiting energy]\label{prop:atom_limsup}
For atoms $X_j=(z_j,r_j)$ of the limiting measure $\mathfrak{q}$, the constants  
\begin{align}
\mathsf{p}_j=\mathfrak{p}(\{X_j\}),
\qquad
\mathsf{q}_j=\mathfrak{q}(\{X_j\}), \label{eq:mjbj}
\end{align}
satisfy 
\begin{align}
0\le \sum_{j=1}^{\infty}\mathsf{p}_j\le1, \quad 0\le  \sum_{j=1}^{\infty}\mathsf{q}_j\le1.  \label{eq:mjbj2}
\end{align}
Moreover, we have 
\begin{align}
\frac{1}{24\pi^{2}e} \le  \mathcal{I}_*
 \le
 \sum_{j=1}^{\infty} e_{r_j}(\mathsf{p}_j,\mathsf{q}_j).
 \label{eq:Lambdabound}
\end{align}
\end{lemma}

\begin{proof}
We choose an increasing sequence of rectangles $K_n=[-Z_n,Z_n]\times[1/R_n,R_n]\Subset\bbR^2_+$ exhausting $\bbR^2_+$ such that $\mathfrak{p}(\partial K_n)=\mathfrak{q}(\partial K_n)=0$ and
\[
\varepsilon_n=\mathfrak{p}(\bbR^2_+\setminus K_n)+\mathfrak{q}(\bbR^2_+\setminus K_n) \to 0 \quad \textrm{as} \,\, n \to \infty.
\]
For a fixed $\tau>0$, we consider the atoms of $\mathfrak{q}$ in $K_n$
with atomic weight greater than $\tau$, and denote their index set by
\[
J_{n,\tau}=\{j\in \mathbb{N}\ :\ X_j\in K_n,\ \mathsf{q}_j>\tau\}.
\]
The set $J_{n,\tau}$ is finite since
\begin{align*}
\sum_{j=1}^{\infty}\mathsf{q}_j= \sum_{j=1}^{\infty}\mathfrak{q}(\{X_j\})=\mathfrak{q}\left(\bigcup_{j=1}^{\infty} \{X_j\}\right)\le 1.
\end{align*}
For a large atom at $X_j$ with $j\in J_{n,\tau}$, we choose a ball $B_{\delta_j}=B_{\delta_j}(X_j) \Subset K_n^\circ$ with pairwise disjoint closures and satisfying $\mathfrak{p}(\partial B_{\delta_j})=\mathfrak{q}(\partial B_{\delta_j})=0$. For fixed $n$ and $\tau$, we extract large atoms of $\zeta_k$ in $K_n$ by setting  
\begin{align*}
\zeta_k^{\textrm{atom}}=\sum_{j\in J_{n,\tau}}\zeta_{k,j}^{\textrm{atom}},\quad \zeta_{k,j}^{\textrm{atom}}=\zeta_k\mathbf1_{B_{\delta_j}}.
\end{align*}
We define the remainder and far-field parts by 
\[
\zeta_k^{\textrm{rem}}=\zeta_k\mathbf1_{K_n\setminus\bigcup_{j\in J_{n,\tau}}B_{\delta_j}}, \quad \zeta_k^{\textrm{far}}=\zeta_k\mathbf1_{\bbR^2_+\setminus K_n}.
\]
Using the decomposition $\zeta_k=\zeta_k^{\textrm{atom}}+\zeta_k^{\textrm{rem}}+\zeta_k^{\textrm{far}}$, we estimate 
\begin{align*}
E[\zeta_k]=E[\zeta_k^{\textrm{atom}}]+E[\zeta_k^{\textrm{rem}}]+E[\zeta_k^{\textrm{far}}]+2E_{\rm int}[\zeta_k^{\textrm{atom}},\zeta_k^{\textrm{rem}}]+2E_{\rm int}[\zeta_k^{\textrm{far}},\zeta_k^{\textrm{atom}}+\zeta_k^{\textrm{rem}}].
\end{align*}
The first term is further decomposed into a sum of atomic energies by using 
\[
E[\zeta_k^{\textrm{atom}}]=\sum_{j\in J_{n,\tau}}E[\zeta_{k,j}^{\textrm{atom}}]+ 2 \sum_{\substack{i,j\in J_{n,\tau}\\ i<j}}
E_{\mathrm{int}}[\zeta_{k,i}^{\textrm{atom}},\zeta_{k,j}^{\textrm{atom}}].
\]
Since the balls have pairwise disjoint closures, the Green function estimate \eqref{eq:Green_bound} gives $E_{\mathrm{int}}[\zeta_{k,i}^{\textrm{atom}},\zeta_{k,j}^{\textrm{atom}}]=O(1)$ as $k\to\infty$. Thus, we have 
\[
\lim_{k\to\infty}(p_k-1)\left( E[\zeta_k^{\textrm{atom}}]-\sum_{j\in J_{n,\tau}}E[\zeta_{k,j}^{\textrm{atom}}] \right) = 0.
\]
The same holds for smaller concentric balls.
For the radii $\delta_j\to 0$ whose boundaries have zero $\mathfrak{p}$- and $\mathfrak{q}$-measure, applying the single atomic bound \eqref{eq:atombound} yields
\[
\lim_{\delta_{j}\to 0}\limsup_{k\to\infty}(p_k-1) E[\zeta_k^{\textrm{atom}}]  \le \sum_{j\in J_{n,\tau}}e_{r_j}(\mathsf{p}_j,\mathsf{q}_j).
\]
To estimate $E[\zeta_k^{\textrm{rem}}]$, we define the remainder region $A_n$ by 
\begin{align*}
A_n=K_n\setminus \bigcup_{j\in J_{n,\tau}}B_{\delta_j}.
\end{align*}
Using $\mathfrak{q}(\partial K_n)=0$ and $\mathfrak{q}(\partial B_{\delta_j})=0$, we have 
\begin{align*}
0\le \mathfrak{q}(\partial A_n)\le \mathfrak{q}(\partial K_n)+\sum_{j\in J_{n,\tau}} \mathfrak{q}(\partial B_{\delta_j})=0.
\end{align*}
Since $\mathfrak{q}_k$ weakly converges to $\mathfrak{q}$ by \eqref{eq:meas_lim} and $\mathfrak{q}(\partial A_n)=0$, Portmanteau's theorem yields the convergence of the restricted measure 
\begin{align*}
\mathfrak{q}_k\!\restriction_{A_n}\rightharpoonup  \mathfrak{q}\!\restriction_{A_n}\quad \textrm{in}\ \mathcal{M}_+(K_n),
\end{align*}
where $\mathfrak{q}\!\restriction_{A_n}$ denotes the restriction of $\mathfrak{q}$ to $A_n$. Then, $g_k=\zeta_k^{\textrm{rem}}$ and $\mathfrak{m}_k=2\pi r g_k^{p_k}\,dz\,dr=\mathfrak{q}_k\!\restriction_{A_n}$ satisfy the conditions \eqref{eq:SAEE1}-\eqref{eq:SAEE4} for $K=K_n$ and $\mathfrak{m}=\mathfrak{q}\!\restriction_{A_n}$.
It follows from \eqref{eq:SAEE5} that 
\begin{equation*}
\limsup_{k\to\infty}(p_k-1)E[\zeta_k^{\textrm{rem}}]\le C_{K_n}\tau.
\end{equation*}
Here, we emphasize that $Q_{p_k}[g_k] \le 1$ and the constant $C_{K_n}$ is independent of $\tau$ and $\delta_j$.
To estimate $E[\zeta_k^{\textrm{far}}]$, we invoke 
\[
P[\zeta_k^{\textrm{far}}]\to\mathfrak{p}(\bbR^2_+\setminus K_n), \quad Q_{p_k}^{p_k}[\zeta_k^{\textrm{far}}]\to\mathfrak{q}(\bbR^2_+\setminus K_n),
\]
by the choice of $K_n$. Thus, the energy estimate \eqref{eq:en_bound} and \eqref{eq:unifrombound2} imply 
\begin{equation*}
\limsup_{k\to\infty}(p_k-1)E[\zeta_k^{\textrm{far}}]\le C\varepsilon_n^{2}.
\end{equation*}
Using the pointwise estimates $\zeta_k^{\textrm{atom}},\zeta_k^{\textrm{rem}},\zeta_k^{\textrm{far}}\le \zeta_k$ and the Cauchy--Schwarz inequality \eqref{eq:p1_en_CS}, we estimate other remainder terms by 
\begin{align*}
\limsup_{k\to\infty}(p_k-1)|E_{\rm int}[\zeta_k^{\textrm{atom}},\zeta_k^{\textrm{rem}}]|&\le 2\sqrt{\limsup_{k\to\infty}(p_k-1)|E[\zeta_k^{\textrm{atom}}]}\sqrt{C_{K_n}\tau},\\
\limsup_{k\to\infty}(p_k-1)|E_{\rm int}[\zeta_k^{\textrm{far}},\zeta_k^{\textrm{atom}}+\zeta_k^{\textrm{rem}}]|&\le 2 \sqrt{C}\varepsilon_n\left(\sqrt{\limsup_{k\to\infty}(p_k-1)|E[\zeta_k^{\textrm{atom}}]}+\sqrt{C_{K_n}\tau}\right).
\end{align*}
Letting $\delta_j\to 0$ yields 
\[
\mathcal{I}_* \le\sum_{j\in J_{n,\tau}}e_{r_j}(\mathsf{p}_j,\mathsf{q}_j)
+C_{K_n}\tau
+C\varepsilon_n^{2}
+2\sqrt{\sum_{j\in J_{n,\tau}}e_{r_j}(\mathsf{p}_j,\mathsf{q}_j)}\sqrt{C_{K_n}\tau}
+2 \sqrt{C}\varepsilon_n\left(\sqrt{\sum_{j\in J_{n,\tau}}e_{r_j}(\mathsf{p}_j,\mathsf{q}_j)}+\sqrt{C_{K_n}\tau}\right).
\]
Letting $\tau\to0$ and then $n\to\infty$ yields \eqref{eq:Lambdabound}.
\end{proof}

\subsection{Sharp optimization and rigidity}\label{sec:rigidity}

\begin{proof}[Proof of Theorem \ref{t:MT4}(ii) and (iii)]
By definition of the atomic energy \eqref{eq:En_def} and \eqref{eq:mbR}, $e_{r_j}(\mathsf{p}_j,\mathsf{q}_j)$ for $\mathsf{p}_j>0$ is expressed as
\begin{align*}
 e_{r_j}(\mathsf{p}_j,\mathsf{q}_j)
 =
\frac{\mathsf{p}_j^{1/2}\mathsf{q}_j^{3/2}}{16\pi^2}
 \frac{\log t_j}{t_j^{3/2}},
 \qquad
 t_j=\frac{r_j^2\mathsf{q}_j}{\mathsf{p}_j}\ge1.
\end{align*}
Since the function $t^{-3/2}\log t$ attains its unique maximum
$2/(3e)$ at $t=e^{2/3}$, we have 
\begin{align*}
 e_{r_j}(\mathsf{p}_j,\mathsf{q}_j)
 \le
 \frac{1}{24\pi^2e}\mathsf{p}_j^{1/2}\mathsf{q}_j^{3/2}.
\end{align*}
This inequality also holds when $\mathsf{p}_j=0$. Using the conditions \eqref{eq:mjbj2}, we obtain
\begin{equation*}
\frac{1}{24\pi^2e}
 \le \mathcal{I}_*
 \le
 \frac{1}{24\pi^2e}
 \sum_{j=1}^{\infty}\mathsf{p}_j^{1/2}\mathsf{q}_j^{3/2}
 \le
 \frac{1}{24\pi^2e}
 \left(\sum_{j=1}^{\infty}\mathsf{p}_j\right)^{1/2}
\left(\sum_{j=1}^{\infty}\mathsf{q}_j^3\right)^{1/2}
 \le
 \frac{1}{24\pi^2e}.
\end{equation*}
Hence $\mathcal{I}_* = (24\pi^2e)^{-1}$. Since every subsequential limit equals $(24\pi^2e)^{-1}$, the convergence \eqref{eq:p1_exact_en} follows.
Equality in this estimate implies
\begin{align*}
\sum_{j=1}^{\infty}\mathsf{p}_j^{1/2}\mathsf{q}_j^{3/2}= \left(\sum_{j=1}^{\infty}\mathsf{p}_j\right)^{1/2}
\left(\sum_{j=1}^{\infty}\mathsf{q}_j^3\right)^{1/2} = 1.
\end{align*}
By the conditions \eqref{eq:mjbj2}, we have 
\begin{align*}
\sum_{j=1}^{\infty}\mathsf{q}_j^3= 1.
\end{align*}
Thus, there exists a unique $j_*$ such that $\mathsf{q}_{j}=\delta_{j_*,j}$.
Moreover,  
\[
 1 = \sum_{j=1}^{\infty}\mathsf{p}_j^{1/2}\mathsf{q}_j^{3/2}=\mathsf{p}_{j_*}^{1/2}.
\]
Hence $\mathsf{p}_{j}=\delta_{j_*,j}$. Since $e_{r_{j_*}}(\mathsf{p}_{j_{*}},\mathsf{q}_{j_{*}})$ attains the maximum at $t_{j_*}=r_{j_*}^{2}\mathsf{q}_{j_{*}}/\mathsf{p}_{j_{*}}=e^{2/3}$ for the atom $X_{j_*}=(z_{j_*},r_{j_*})$, we have $r_{j_*}=e^{1/3}$.
After a translation in the $z$-direction, we may assume $z_{j_*}=0$.
Then $\mathfrak{p}=\mathfrak{q}=\delta_{(0,e^{1/3})}$, which proves \eqref{eq:p1_mu_concen}.
\end{proof}

\subsection{Toroidal vortex rings for small \texorpdfstring{$p$}{p}}\label{sec:no_touch_fin}

It remains to show Theorem \ref{t:MT4}(i) and (iv).
We show that maximizers admit a vorticity-free strip near the axis for sufficiently small $p>1$ using the Euler--Lagrange equation \eqref{eq:EL_intro0} and the energy inequality \eqref{eq:variational-inequality2}.

\begin{proposition}[Vorticity-free strip]\label{lem:vacuum}
Let $p_k\to 1^{+}$ and let $\zeta_k\in S_{p_k}$ satisfy \eqref{eq:p1_mu_concen}.
For every $0<\delta_1<\delta_2<e^{1/3}$, the following holds for all sufficiently large $k$: 
\begin{align}
\zeta_k(z,r)=0 \quad \textrm{on} \quad \bbR\times[\delta_1,\delta_2].  \label{eq:VFS}
\end{align}
\end{proposition}

\begin{proof}
For $X=(z,r)\in A=\bbR\times[\delta_1,\delta_2]$ and $\delta>0$, we set 
\begin{align*}
\varphi_k =\frac{\mathcal{G}\zeta_k}{r^{2}}= 
\frac{1}{r^2}\int_{\mathbb{R}^{2}_{+}}G(X,Y)\zeta_k(Y)Y_2\,dY 
=\frac{1}{X_2^2}  \left(\int_{\{|X-Y|<\delta\}}+ \int_{\{|X-Y|\ge\delta\}} \right)= \varphi_k^1 + \varphi_k^2.
\end{align*}
We take $0< \delta_0<\delta_1$ and $\delta_2 < \delta_3 < e^{1/3}$, and choose $\delta>0$ so that $2\delta < \min\{\delta_1-\delta_0, \delta_3-\delta_2 \}$. Then, $Y\in A^+=\bbR\times[\delta_0,\delta_3]$ if $|X-Y|<\delta$.
By \eqref{eq:p1_loc_ker}, for $X$ and $Y$ satisfying $|X-Y|<\delta$, we have 
\[
\frac{1}{X_2^2}G(X,Y)\le C\left(1+\log_+\frac{1}{|X-Y|}\right).
\]
Making the change of variables $V  = Y-X$ and using $\delta_0 \le X_2+V_2 \le \delta_3$, we estimate 
\[
\begin{aligned}
    \left(\int_{B_\delta(X)} \left(1+\log_+\frac{1}{|X-Y|}\right)^q Y_2\,dY\right)^{1/q} &= \left(\int_{B_\delta(0)} \left(1+\log_+\frac{1}{|V|}\right)^q (X_2+V_2)\,dV\right)^{1/q}\\
&\le \delta_3^{1/q} \left\|1+\log_+\frac{1}{|\cdot|}\right\|_{L^q(B_\delta(0))} \le Cq.
\end{aligned}
\]
Applying H\"older's inequality with $q=p_k/(p_k-1)$, we obtain
\[
(p_k-1)\sup_A\varphi_k^1
\le Cp_kQ_{p_k}[\zeta_k \mathbf{1}_{A^+}]=Cp_k \mathfrak{q}_k(A^{+})^{1/p_k}.
\]
The right-hand side vanishes as $k\to\infty$ by $\delta_3 < e^{1/3}$ and the weak convergence of the measure \eqref{eq:p1_mu_concen}.

If $|X-Y| \ge \delta$, $G(X,Y)/X_2^2 \le C_\delta Y_2^2$ by the pointwise estimate \eqref{eq:Green_bound} for $\tau=3/2$.
Therefore,
\[
\sup_A\varphi_k^{2} \le C_\delta\int_{\bbR^2_+}Y_2^3\zeta_k(Y)\,dY=\frac{C_\delta}{2\pi}P[\zeta_k]=\frac{C_\delta}{2\pi}.
\]
We thus obtain 
\begin{align*}
\lim_{k\to\infty}(p_k-1)\sup_{ A^{+}}\varphi_k=0.
\end{align*}
By the Euler--Lagrange equation \eqref{eq:EL_intro0}, 
\[
\kappa_{p_k}\zeta_k^{p_k-1} = r^2\left(\varphi_k-\frac{W_{p_k}}{2}\right)_+ = \frac{r^2}{p_k - 1}\left( \left( p_k - 1 \right) \left( \varphi_k - \frac{W_{p_k}}{2} \right)\right)_+. 
\]
Since $(p_k-1)W_{p_k}$ converges to a positive constant by \eqref{eq:p1_exact_en}, on the strip $A$, the right-hand side is zero for sufficiently large $k$. We obtain \eqref{eq:VFS}.
\end{proof}

\begin{lemma}\label{lem:noin}
The property \eqref{eq:VFS} is strengthened to 
\begin{align}
\zeta_k=0\quad \textrm{on}\ \mathbb{R}\times (0,\delta_2].  \label{eq:VFS2}
\end{align}
\end{lemma}

\begin{proof}
By \eqref{eq:VFS}, we set
\begin{align*}
\zeta_k =\zeta_k\mathbf1_{\{r<\delta_1\}}
+\zeta_k\mathbf1_{\{r>\delta_2\}}
=\zeta_k^{1} + \zeta_k^{2}.
\end{align*}
Accordingly, we define the stream functions by 
\begin{align*}
\varphi_{k}= \frac{1}{r^{2}}\mathcal{G}\zeta^{1}_k+\frac{1}{r^{2}}\mathcal{G}\zeta^{2}_k=
\varphi_{k}^{1}+\varphi_{k}^{2}.
\end{align*}
Using $P[\zeta_k]=Q_{p_k}^{p_k}[\zeta_k]=1$, we define the constants $\mu_k$ and $\nu_k$ by 
\begin{align*}
P[\zeta_k^{1}]&=\mu_k,\quad Q_{p_k}^{p_k}[\zeta_k^{1}]=\nu_k,\\
P[\zeta_k^{2}] &= 1 - \mu_k,\quad Q_{p_k}^{p_k}[\zeta_k^{2}]=1-\nu_k.
\end{align*}
The constants $\mu_k$ and $\nu_k$ vanish by the weak convergence of the measure \eqref{eq:p1_mu_concen}. 

It suffices to show that $\nu_k =0$ for all large $k$. We apply the Cauchy--Schwarz inequality \eqref{eq:p1_en_CS} to estimate  
\begin{align*}
I_{p_k}=E[\zeta_k]=E[\zeta_k^{1}]+E[\zeta_k^{2}]+2E_{\textrm{int}}[\zeta_{k}^{1},\zeta_{k}^{2}]\le \left(\sqrt{E[\zeta^{1}_{k}]}+\sqrt{E[\zeta^{2}_{k}]}\right)^{2}.
\end{align*}
Using the energy inequality \eqref{eq:variational-inequality2}, we obtain 
\begin{align*}
1\le \mu_k^{\alpha_k/2}\nu_k^{\beta_k/(2p_k)}+(1-\mu_k)^{\alpha_k/2}(1-\nu_k)^{\beta_k/(2p_k)}
\end{align*}
Since $\beta_k/(2p_k)<1$ for large $k$, by the pointwise estimate 
\begin{align*}
(1-\nu_k)^{\beta_k/(2p_k)}\le 1-\frac{\beta_k}{2p_k}\nu_k,
\end{align*}
it follows that 
\begin{align}
\frac{\beta_k}{2p_k}\nu_{k}^{1-\beta_k/(2p_k)}\le \mu_k^{\alpha_k/2}.  \label{eq:upperbound}
\end{align}
On the other hand, using $\dist (\supp \zeta_k^{1}, \supp \zeta_k^{2}) \ge \delta_2 - \delta_1 > 0$ and the pointwise estimate \eqref{eq:Green_bound} for $\tau=3/2$, we have 
\begin{equation*}
\sup_{X_2<\delta_1}\varphi_{k}^{2}\le \frac{C}{(\delta_2-\delta_1)^{3}}P[\zeta_k^{2}]\le \frac{C}{(\delta_2-\delta_1)^{3}}.
\end{equation*}
Since $W_{p_k}$ diverges, the right-hand side is smaller than $W_{p_k}/4$ for sufficiently large $k$. By the Euler--Lagrange equation \eqref{eq:EL_intro0}, we have $\varphi_k>W_{p_k}/2$ on $\{\zeta_k^{1} >0\} \subset \supp \zeta_k$.
Therefore $\varphi_{k}^{1}\ge W_{p_k}/4$ on $\{\zeta_k^{1} >0\}$.
Thus,
\[
E[\zeta_k^{1}]= \pi \int_{\bbR^2_+}\varphi_{k}^{1}\, \zeta_k^{1}X_2^3\,dX\ge\frac{W_{p_k}}{8}\mu_k.
\]
It follows from \eqref{eq:variational-inequality2} that 
\begin{align*}
\frac{W_{p_k}}{8}\mu_k\le E[\zeta_k^{1}] \le I_{p_k}\mu_k^{\alpha_{k}}\nu_k^{\beta_{k}/p_k}.
\end{align*}
Since $W_k/2=\alpha_{k}I_{p_k}$ by \eqref{eq:speed}, we have 
\begin{align}
\mu_k^{1-\alpha_{k}} \le \frac{4}{\alpha_{k}} \nu_k^{\beta_{k}/p_k}. \label{eq:lowerbound}
\end{align}
Suppose that $\nu_k>0$ along a subsequence.
Since $\alpha_k<1$ for large $k$, by taking \eqref{eq:lowerbound} to the power $\alpha_k/[2(1-\alpha_k)]$ and substituting into \eqref{eq:upperbound}, we get
\begin{align*}
\frac{\beta_k}{2p_k}
\left(\frac{\alpha_k}{4}\right)^{\alpha_k/[2(1-\alpha_k)]}
\le \nu_k^{\beta_k/[2p_k(1-\alpha_k)]-1}.
\end{align*}
Since $\alpha_k\to 1/2$ and $\beta_k\to 3/2$, the left-hand side tends to a positive constant while the exponent on the right-hand side tends to $1/2$.
The right-hand side therefore tends to zero as $\nu_k \to 0$, which yields a contradiction.
Thus $\nu_k=0$ for all large $k$, and \eqref{eq:VFS2} follows.
\end{proof}

\begin{proof}[Proof of Theorem \ref{t:MT4}]
We first show (i). By Lemma \ref{lem:p2_simply}, the support of $\zeta\in S_p$ is simply connected for $1<p<2$. Suppose, for a contradiction, that for a sequence $p_k\to 1^{+}$, there exists $\zeta_k\in S_{p_k}$ such that $\dist(\supp\zeta_k,\{r=0\})=0$.
After translation and passage to a subsequence, by \eqref{eq:VFS2}, for any $0<\delta_2<e^{1/3}$, we have 
\[
\supp\zeta_k\subset\{r\ge \delta_2\}
\]
for all sufficiently large $k$. This gives a contradiction. Thus there exists $1<p_S<2$ such that for every $1<p<p_S$ and every $\zeta\in S_p$, the vorticity support of $\zeta$ forms a solid torus in $\mathbb{R}^{3}$.

For all sufficiently large $k$, we have 
\[
2\pi \Gamma[\zeta_k]= 2\pi\int_{\mathbb R^2_+}\zeta_kX_2\,dX =\int_{\mathbb R^2_+}\frac{1}{X_2^{2}} \left(2\pi X_2^3\zeta_k\,dX\right)=\mathfrak{p}_k[X_2^{-2}\mathbf{1}_{\{X_2>\delta_2\}}] \to  e^{-2/3}.
\]
This proves the convergence \eqref{eq:L1masslimit} of (iv), and the proof is complete.
\end{proof}

\section{Core profile and Kelvin--Hicks formula}\label{s:8}

We now complete the proof of Theorem~\ref{t:MT5}.
Our first goal is to derive a two-term asymptotic expansion of the maximal energy $I_p$ as $p \to 1^+$.
We begin by constructing a test function whose core profile is given by the log-HLS optimizer, obtaining a lower bound for 
\[
\liminf_{p\to1^+} \left( I_p-\frac{1}{24\pi^2e(p-1)} \right).
\]
Next, using a lower bound for the logarithmic potential, we show that the diameter $d_p$ of the vortex core decays exponentially:
\[
d_p\le C\exp\left(-\frac{c}{p-1}\right).
\]
Combining this estimate with a log-type HLS inequality, we derive a sharp upper bound for the logarithmic interaction energy.
We then use the exponential decay estimate and the asymptotic expansion of $F(s)$ to express the energy associated with the Green kernel in terms of the logarithmic interaction energy, with an error of $o(1)$.
This yields a matching sharp upper bound for $I_p$ and establishes its two-term asymptotic expansion.

\subsection{Refined sharp lower bound}

We establish a refined lower bound for $I_p$.
To this end, we first construct a probability density by suitably modifying the log-HLS optimizer $\omega_*$ and then rescale it.
We then extract the leading contribution to the energy by using an asymptotic expansion of the Green's kernel.

\begin{proposition}
Let $\omega_*$ be the optimizer \eqref{eq:optimizer}.
For $L \ge 1$, we set 
\begin{align}
\omega_{*,L}(X)=\omega_{*}(X)\left(1+\frac{1}{L^{2}}\right)\mathbf1_{\{|X|<L\}}. \label{eq:etaT}
\end{align}
Then the following identities hold:
\begin{align}
\int_{\mathbb{R}^{2}} X_2^{2m+1}\omega_{*,L}^{q}dX &= 0,\quad m\in \mathbb{N}_0,\,\, q\ge 1,   \label{eq:etamoment}\\
\int_{\mathbb{R}^{2}} \omega_{*,L} dX &= 1, \label{eq:etamass}\\
H[\omega_{*,L}]&=\frac{1}{L^{2}}-\frac{1}{L^{4}}\log(1+L^{2}) \label{eq:etaH},\\
\log \int_{\bbR^2} \omega^p_{*,L}\,dX &=  -(p-1) S[\omega_{*,L}]+ o(p-1) \quad \textrm{as }\,\, p\to 1^+.  \label{eq:logetaasym}
\end{align}
\end{proposition}

\begin{proof}
The property \eqref{eq:etamoment} follows from the radial symmetry of $\omega_{*,L}$. With the normalization factor $1+1/L^{2}$, $\omega_{*,L}$ satisfies \eqref{eq:etamass}. By direct computation, we obtain 
\begin{align*}
    -S[\omega_{*,L}] & = -\log\pi-2
+\log\left(1+\frac{1}{L^2}\right)
+\frac{2}{L^2}\log(1+L^2),\\
U[\omega_{*,L}]-1-\log\pi & = 1+\frac1{L^2}
-\log\left(1+\frac1{L^2}\right)
-\left(\frac{2}{L^{2}}+\frac{1}{L^{4}} \right)\log(1+L^2).
\end{align*}
Summing up the two equalities yields \eqref{eq:etaH}. Since $\omega_{*,L} $ is bounded, we expand $\omega^{p}_{*,L}$ as 
\begin{align*}
\omega^p_{*,L} = \omega_{*,L} e^{(p-1)\log \omega_{*,L}} =  \omega_{*,L} + (p-1)\omega_{*,L} \log \omega_{*,L} + O((p-1)^2).
\end{align*}
Integrating over $\mathbb{R}^{2}$, we obtain 
\begin{align*}
\int_{\bbR^2}\omega^p_{*,L}\,dX = 1 - (p-1)S[\omega_{*,L}] + O((p-1)^2).
\end{align*}
Since $\log (1 + t) = t + O(t^2)$, the desired result follows.
\end{proof}

\begin{lemma}[Refined sharp lower bound]\label{lem:p1_W_low}
We have
\begin{equation}\label{eq:p1_W_low}
\liminf_{p\to1^{+}}\left(I_p-\frac1{24\pi^2e(p-1)}\right) \ge \frac{1 }{ 8 \pi^2 e}\left(\log8\sqrt2 \pi - \frac{2}{3}\right).
\end{equation}
\end{lemma}

\begin{proof}
We refine the thin-ring construction in Lemma~\ref{lem:p1_low} using the log-HLS optimizer $\omega_*$ and with the ring radius fixed at $R = e^{1/3}$.
Let $\omega_{*,L}$ be given by \eqref{eq:etaT}.
For $\varrho > 0$, we rescale it around the point $X_0=(0, e^{1/3})$ as
\[
\omega(X) = \frac{A}{\varrho^2} \omega_{*,L} \left(\frac{X - X_0}{\varrho} \right),\quad A=\frac{1}{2\pi e^{2/3}}.
\]
For $0<\varrho < L^{-1}e^{1/3}$, we have $\supp \omega \subset B_{\varrho L}(X_0)$.
We set $\zeta=\omega/e^{1/3}$. By a change of variables and \eqref{eq:etamoment}, we obtain   
\begin{align*}
    Q_p^p[\zeta]  =  \frac{2\pi A^{p}}{ e^{p/3}\varrho^{2(p-1)}}\int_{\bbR^2}\omega^p_{*,L}(X) (e^{1/3} + \varrho X_2)\,dX 
=  \frac{2\pi A^{p}}{e^{(p-1)/3}\varrho^{2(p-1)}}\int_{\bbR^2}\omega^p_{*,L}(X)\,dX.
\end{align*}    
We choose $\varrho=\varrho_p$ so that $Q_{p}^{p}[\zeta]= 1$.
Namely, 
\begin{equation*}
\varrho_p =\frac{(2\pi A^{p})^{\frac{1}{2(p-1)}}}{e^{1/6}} \left(\int_{\bbR^2} \omega^p_{*,L}(X) \,dX\right)^{\frac{1}{2(p-1)}}.
\end{equation*}
By \eqref{eq:logetaasym} and $2\pi A<1$, the constant $\varrho_p$ tends to zero.
In particular, $\varrho_p  <L^{-1} e^{1/3}$ for $p>1$ sufficiently close to $1$. By taking the logarithm,
\[
(p - 1)\log\frac1{\varrho_p} = -\frac{1}{2} \log 2\pi -\frac{p}{2}\log A+ \frac{1}{6}\left(p-1\right) - \frac{1}{2}\log\int_{\bbR^2}\omega^p_{*,L}\,dX.
\]
Applying \eqref{eq:logetaasym} for $\omega_{*,L}$, we obtain
\[
\log\frac1{\varrho_p} = \frac{1}{3(p-1)} + \frac{1}{2}\log(2\pi e)+\frac{1}{2}S[\omega_{*,L}] + o(1).
\]
It follows from \eqref{eq:streamfunction} that 
\begin{align*}
E[\zeta]=\frac{1}{2e^{2/3}} \iint_{\mathbb{R}^{2}_{+}\times \mathbb{R}^{2}_{+}}(X_2Y_2)^{3/2}F\left(\frac{|X-Y|^{2}}{X_2Y_2}\right)\omega(X)\omega(Y)dXdY.
\end{align*}
We set 
\begin{align*}
s = \frac{\varrho_p^{2} |X-Y|^{2}}{(e^{1/3}+\varrho_p X_2)(e^{1/3}+\varrho_p Y_2)}.
\end{align*}
By a change of variables,
\begin{align*}
E[\zeta]
=\frac{A^{2}}{2e^{2/3}}\iint_{B_{L}(0)\times B_{L}(0)}(e^{1/3}+\varrho_p X_2)^{3/2}(e^{1/3}+\varrho_p Y_2)^{3/2} F\left(s\right)\omega_{*,L}(X)\omega_{*,L}(Y)dXdY.
\end{align*}
Since $X,Y\in B_{L}(0)$, we have 
\begin{align*}
(e^{1/3}+\varrho_p X_2)(e^{1/3}+\varrho_p Y_2)&=e^{2/3}+O(\varrho_p),\\
\log\frac{(e^{1/3}+\varrho_p X_2)(e^{1/3}+\varrho_p Y_2)}{e^{2/3}}&=O(\varrho_p),\\
s&=O(\varrho_{p}^{2}),\\
\log\frac{8}{\sqrt{s}}&=\log{\frac{8e^{1/3}}{\varrho_p}}+\log\frac{1}{|X-Y|}+O(\varrho_p).
\end{align*}
For fixed $L >  0$ and $X,Y\in B_L(0)$ with $X\neq Y$, the expansion \eqref{eq:p1_F_log} yields
\begin{align*}
(e^{1/3}+\varrho_p X_2)^{3/2}(e^{1/3}+\varrho_p Y_2)^{3/2}F(s)
\quad 
&=e \cdot \left( \log\frac{8e^{1/3}}{\varrho_p} +\log\frac{1}{|X-Y|}-2 \right)
\qquad \\
& + O_L\left( \varrho_p\left( \log\frac{1}{\varrho_p} +|\log|X-Y|| + 1 \right)
\right).
\end{align*}
Since $\omega_{*,L}$ is bounded and compactly supported, the error term contributes $O_L(\varrho_p|\log \varrho_p|)$ after integration over $B_L(0)\times B_L(0)$.
Using this expansion, we deduce that 
\begin{equation*}
E[\zeta] = \frac{A^2e^{1/3}}{2}\left( \log\frac{8e^{1/3}}{\varrho_p} - \frac{3}{2} +\frac{1}{2}\log\pi -\frac{1}{2}U[\omega_{*,L}] \right) + O\left(\varrho_p\log\frac1{\varrho_p}\right).
\end{equation*}
Substituting the expansion of $-\log \varrho_p$, we obtain 
\begin{equation*}
\log\frac{8e^{1/3}}{\varrho_p} - \frac{3}{2} +\frac{1}{2}\log\pi-\frac{1}{2}U[\omega_{*,L}]
=\log{8\sqrt{2}\pi}-\frac{2}{3}-\frac{1}{2}H[\omega_{*,L}]+\frac{1}{3(p-1)}+o(1).
\end{equation*}
It follows that
\begin{equation*}
\begin{aligned}
E[\zeta] - \frac{1}{24\pi^2e(p-1)}
& = \frac{1}{8\pi^2e}\bigg( \log{8\sqrt{2}\pi}-\frac{2}{3}-\frac{1}{2}H[\omega_{*,L}] \bigg) + o(1).
\end{aligned}
\end{equation*}
Since 
\begin{align*}
P[\zeta]  =  \frac{1}{e}\int_{\bbR^2}(e^{1/3} + \varrho_pX_2)^3 \omega_{*,L}(X)\,dX = 1 + \frac{3\varrho_p^2}{e^{2/3}}\int_{\bbR^2}X_2^2\omega_{*,L}(X)\,dX \to   1,
\end{align*}
we normalize $\zeta$ by 
\[
\widetilde\zeta = \frac{1}{P[\zeta]}\zeta \in K_p.
\]
By the maximality of $I_p$, 
\begin{align*}
I_p- \frac{1}{24\pi^2e(p-1)}
&\ge E[\tilde{\zeta}]- \frac{1}{24\pi^2e(p-1)}=\frac{1}{P[\zeta]^{2}}E[\zeta]- \frac{1}{24\pi^2e(p-1)}\\
&=\frac{1-P[\zeta]^{2}}{P[\zeta]^{2}}E[\zeta]+\frac{1}{8\pi^{2}e}\left(\log{8\sqrt{2}\pi}-\frac{2}{3}-\frac{1}{2}H[\omega_{*,L}] \right)+o(1).
\end{align*}
By \eqref{eq:etaH}, taking the liminf as $p\to 1^{+}$ and then letting $L \to \infty$ yields \eqref{eq:p1_W_low}.
\end{proof}

\subsection{Exponential decay of the vortex-core diameter}

We next derive an exponential decay estimate for the diameter of the vortex core.
We first show that, for any fixed radius, the support of a maximizer is contained in a ball centered at $(0,e^{1/3})$ for $p >  1$ sufficiently close to $1$.
This yields uniform bounds on the $X_2^2$-moment and the $L^p$ norm of the probability density.
Using these bounds, we establish a lower bound of order $(p-1)^{-1}$ for the logarithmic potential, which leads to the exponential decay estimate for the core diameter.

\begin{proposition}
Let $1<p\le \infty$ and let $\zeta_p\in S_p$ be a centered maximizer.
Let $\hat{\omega}_p$ and $R_p$ be as defined in \eqref{eq:hatomega} of Theorem~\ref{t:MT5}.
For every $\delta>0$, there exists $p_{\delta}>1$ such that 
\begin{equation}\label{eq:supportzetaball}
\supp\hat{\omega}_p \subset B_{\delta}(0,R_p),
\end{equation}
for every $1<p\le p_{\delta}$.
Moreover, $\hat{\omega}_p$ satisfies  
\begin{align}\label{eq:normalizedconserve}
\int_{\bbR^2} X_2^2\hat{\omega}_p(X)\,dX = R_p^2,\quad \int_{\bbR^2} \hat{\omega}_p(X)^{p}\,dX
\lesssim 1. 
\end{align}
\end{proposition}

\begin{proof}
For some constants $0 < a_1 < a_2$, we have $\supp\zeta_p \subset \bbR \times [a_1,a_2]$ by \eqref{eq:VFS2} and \eqref{eq:rbounded}.
Taking the constants $a_0<a_1<a_2<a_3$ and $\sigma_0<\sigma_1$, we set $A=\bbR \times [a_1,a_2]\setminus B_{\sigma_1}(0,e^{1/3})$ and $A^{+}=\bbR \times [a_0,a_3]\setminus B_{\sigma_0}(0,e^{1/3})$ so that $A\subset A^{+}$.
We take $\delta>0$ so that the conditions $X\in A$ and $|X-Y|<\delta$ imply $Y\in A^{+}$.
As in the proof of Proposition \ref{lem:vacuum}, for $p_k\to 1^{+}$ and centered $\zeta_k\in S_{p_k}$, we set  
\begin{align*}
\varphi_k =\frac{\mathcal{G}\zeta_k}{r^{2}}
= \frac{1}{r^2}\int_{\mathbb{R}^{2}_{+}}G(X,Y)\zeta_k(Y)Y_2\,dY 
=\frac{1}{X_2^2}  \left(\int_{\{|X-Y|<\delta\}}+ \int_{\{|X-Y|\ge\delta\}} \right)= \varphi_k^1 + \varphi_k^2,
\end{align*}
and estimate 
\[
(p_k-1)\sup_A\varphi_k^1 \le Cp_kQ_{p_k}[ \mathbf{1}_{A^+}]=Cp_k\mathfrak{q}_k(A^{+})^{1/p_k}.
\]
The right-hand side tends to zero as $p_k\to 1^{+}$ by the convergence of measures \eqref{eq:p1_mu_concen}.
Since $\varphi^{2}_{k}$ is bounded on $A$, we have 
\[
\lim_{k\to \infty}(p_k-1)\sup_{A} \varphi_k = 0.
\]
By the Euler--Lagrange equation \eqref{eq:EL_intro0}, 
\[
\kappa_{p_k}\zeta_k^{p_k-1} = r^2\left(\varphi_k-\frac{W_k}{2}\right)_+ = \frac{r^2}{p_k - 1}\left( \left( p_k - 1 \right) \left( \varphi_k - \frac{W_k}{2} \right)\right)_+. 
\]
Since $(p_k-1)W_{k}$ converges to a positive constant, we obtain $\zeta_k=0$ on $A=\mathbb{R}\times [a_1,a_2]\setminus B_{\sigma_1}(0,e^{1/3})$ for sufficiently large $k$. Since $R_{p_k}\to e^{1/3}$ by \eqref{eq:L1masslimit}, we obtain \eqref{eq:supportzetaball}. The properties \eqref{eq:normalizedconserve} follow from the properties of the maximizer $P[\zeta_p]=Q_p^{p}[\zeta_p]=1$ and $\supp\hat{\omega}_p\subset \mathbb{R}\times [a_1,a_2]$.
\end{proof}

\begin{lemma}\label{lem:p1_ex_supp}
Let $d_p=\operatorname{diam}(\supp\hat{\omega}_p)$. Then
\begin{align}
\supp\hat{\omega}_p \subset \overline{B_{d_p}(0,R_p)}.\label{eq:p1_W_supp_1}
\end{align}
Moreover, there exist constants $C, c > 0$ and $p_0>1$ such that, for every $1<p\le p_0$,
\begin{equation}\label{eq:diambound}
d_p\le C\exp\left(-\frac{c}{p-1}\right).
\end{equation}  
\end{lemma}

\begin{proof}
We write $\Omega_p=\supp{\hat{\omega}_p}$. Since ${\hat{\omega}_p}$ is a probability density, it follows from \eqref{eq:normalizedconserve} that 
\begin{align*}
\min_{X\in \Omega_p}X_2\le R_p\le \max_{X\in \Omega_p}X_2.
\end{align*}
Let $\operatorname{conv}\Omega_p$ denote the convex hull of $\Omega_p$, defined by
\begin{align*}
\operatorname{conv}\Omega_p
=
\left\{
\sum_{j=1}^{N}\lambda_j X_j :
N\in\mathbb{N},\
X_j\in \Omega_p,\
\lambda_j\ge 0,\
\sum_{j=1}^{N}\lambda_j=1
\right\}.
\end{align*}
By the symmetry of $\hat\omega_p$ with respect to $z=0$ and the preceding bounds, we have $(0,R_p)\in \operatorname{conv}\Omega_p$. For $X\in \Omega_p$ and $Y=\sum_{j=1}^{N}\lambda_j X_j\in \operatorname{conv}\Omega_p$, it follows that 
\begin{align*}
|X-Y|=\left|\sum_{j=1}^{N}\lambda_j (X-X_j)\right|\le d_p.
\end{align*}
In particular, this proves \eqref{eq:p1_W_supp_1}.

For $X,Y\in \Omega_p\subset B_{\delta}(0,e^{1/3})$, it follows from \eqref{eq:p1_loc_ker} that  
\begin{align*}
G(X,Y)\le C_1\log\frac{1}{|X-Y|}+C_2.
\end{align*}
Multiplying this inequality by $\hat{\omega}_p(Y)$ and integrating over $\mathbb{R}^{2}$, we obtain   
\begin{align*}
\frac{1}{\Gamma_p}\mathcal{G}\zeta_p(X)=\int_{\mathbb{R}^{2}}G(X,Y)\hat{\omega}_p(Y)dY\le C_1\int_{\mathbb{R}^{2}}\log\frac{1}{|X-Y|}\hat{\omega}_p(Y)dY+C_2.
\end{align*}
Since $\mathcal{G}\zeta_p \ge W_pr^2/2$ on $\supp \zeta_p$ and $(p-1)W_p$ converges to a positive constant, we obtain
\[
\frac{C_3}{p - 1}\le \int_{\bbR^2} \log\frac{1}{|X - Y|} \hat{\omega}_p(Y)\,dY.
\]
This estimate holds for every $X \in \Omega_p$. 

We may assume that $d_p < 1$ by \eqref{eq:supportzetaball}. Let $N \ge 2$ be an integer. Choose $X_1,X_N\in \Omega_p$ such that $|X_1 - X_N| = d_p$. Since $\Omega_p$ is connected by Lemma \ref{lem:p2_simply}, the continuous function $f(X)=|X-X_1|$ maps $\Omega_p$ onto $[0,d_p]$. Hence there exist points $\{X_i\}_{i=1}^{N}\subset \Omega_p$ such that
\[
f(X_i)= \frac{(i - 1)d_p}{N - 1}.
\]
For $i\neq j$, we have 
\begin{align*}
|X_i-X_j|\ge \left||X_i-X_1|-|X_j-X_1|  \right|=\frac{|i-j|d_p}{N-1} \ge \frac{d_p}{N-1}.
\end{align*}
For any point $Y \in \Omega_p$ and any $i\neq j$, we have 
\begin{align*}
\frac{d_p}{N-1}\le |X_i-X_j|\le |X_i-Y|+|X_j-Y|.
\end{align*}
Thus, there exists at most one $i=i_Y$ such that 
\begin{align*}
|X_{i_Y}-Y|< \frac{d_p}{2(N-1)}.
\end{align*}
If such $i_Y$ exists,
\begin{align*}
\sum_{1\le i\le N}\log\frac{1}{|X_i-Y|}
= \sum_{\substack{1\le i\le N\\ i\neq i_Y}}\log\frac{1}{|X_i-Y|}+\log\frac{1}{|X_{i_Y}-Y|}
\le 
(N-1)\log\frac{2(N-1)}{d_p}+\max_{1\le i\le N}\log^{+}\frac{1}{|X_i-Y|}.
\end{align*}
Otherwise, choosing any $k\in\{1,\ldots,N\}$, we have
\begin{align*}
\sum_{1\le i\le N}\log\frac{1}{|X_i-Y|}
= \sum_{\substack{1\le i\le N\\ i\neq k}}\log\frac{1}{|X_i-Y|}+\log\frac{1}{|X_{k}-Y|}
\le 
(N-1)\log\frac{2(N-1)}{d_p}+\max_{1\le i\le N}\log^{+}\frac{1}{|X_i-Y|}.
\end{align*}
We thus obtain
\begin{align*}
\sum_{i=1}^{N}\log\frac{1}{|X_i-Y|}
\le (N-1)\log\frac{2(N-1)}{d_p}+\max_{1\le i\le N}\log^{+}\frac{1}{|X_i-Y|}.
\end{align*}
Multiplying this inequality by $\hat{\omega}_p$ and integrating over $\mathbb R^2$, we obtain 
\begin{align*}
\frac{C_3N}{p-1}\le \sum_{i=1}^{N}\int_{\mathbb{R}^{2}}\log\frac{1}{|X_i-Y|}\hat{\omega}_p(Y)dY
\le (N-1)\log\frac{2(N-1)}{d_p}+\int_{\mathbb{R}^{2}}\max_{1\le i\le N}\log^{+}\frac{1}{|X_i-Y|}\hat{\omega}_p(Y)dY.
\end{align*}
We estimate the last term using $d_p<1$ and $\|\log(1/|\cdot|)\|_{L^q(B_1)} \lesssim q$ for $q = p/(p - 1)$. It follows from \eqref{eq:normalizedconserve} that 
\begin{align*}
\int_{\mathbb{R}^{2}}\max_{1\le i\le N}\log^{+}\frac{1}{|X_i-Y|}\hat{\omega}_p(Y)dY
&\le \left(\int_{\mathbb{R}^{2}}\max_{1\le i\le N}\left(\log^{+}\frac{1}{|X_i-Y|}\right)^{q}dY\right)^{1/q}||\hat{\omega}_p||_{L^{p}(\mathbb{R}^{2})} \\
&\le \left(\sum_{i=1}^{N}\int_{\mathbb{R}^{2}}\left(\log^{+}\frac{1}{|X_i-Y|}\right)^{q}dY\right)^{1/q}||\hat{\omega}_p||_{L^{p}(\mathbb{R}^{2})} \lesssim N^{1/q}q.
\end{align*}
We obtain 
\begin{align*}
\frac{N C_3}{p - 1} \le (N - 1)\log\frac{2(N - 1)}{d_p} + C_4 N^{1/q}q.
\end{align*}
We first choose $N$ sufficiently large so that $NC_3>2C_4$. Since $pN^{1/q}\to1$ as $p\to1^+$, we may then assume that $pN^{1/q}\le2$. Rearranging the preceding inequality, we obtain
\[
0<\frac{NC_3 - 2C_4}{N - 1}\le 
(p - 1)\log\frac{2(N - 1)}{d_p}.
\]
This means that 
\begin{equation*}
d_p \le 2(N-1)\exp\left(-\frac{1}{p - 1}\left(\frac{NC_3-2C_4}{N-1}\right)\right).
\end{equation*}
Thus \eqref{eq:diambound} holds. 
\end{proof}

\subsection{Refined asymptotics}

We next derive a sharp upper bound for $I_p$ by estimating the energy of a maximizer.
To this end, we estimate the logarithmic interaction energy 
\begin{align}
E_{\mathrm{log}}[\omega]
=\frac12(1+\log\pi-U[\omega]) 
=\iint_{\mathbb R^2\times\mathbb R^2}
\log\frac1{|X-Y|}\omega(X)\omega(Y)\,dX\,dY.
 \label{logenergy} 
\end{align}
We first establish a sharp upper bound for the logarithmic interaction energy using a logarithmic-type HLS inequality and the exponential decay estimate for the vortex core diameter.
We then extract the logarithmic interaction energy as the leading term of the Green kernel to obtain a sharp upper bound for $I_p$.
The resulting refined asymptotics of $I_p$ imply the convergence of the free energy, showing that the probability densities $\hat\omega_p$ form a minimizing sequence for the log-HLS inequality.

\begin{proposition}
As $p\to 1^{+}$, 
\begin{align}
&0 < E_{\textrm{log}}[\hat{\omega}_p]  \le \left( \frac{1}{p-1} + \frac{3}{2} \right) \log R_p + \frac{1}{2}\log(2 \pi^2) + \frac{1}{2} + o(1),\label{eq:p1_next_log} \\
&\frac{1}{p-1}\log{\left(\int_{\mathbb{R}^{2}}\hat{\omega}_p^{p}dX\right)}-\log{2\pi}-\left(3+\frac{2}{p-1}\right)\log{R_p}=o(1). \label{eq:logintomegap}
\end{align}    
\end{proposition}

\begin{proof}
The logarithmic interaction energy $E_{\textrm{log}}[\hat{\omega}_p]$ is positive since $d_p=\operatorname{diam}(\supp \hat{\omega}_p)< 1$ by \eqref{eq:diambound} and $|X-Y|<1$ for $X,Y\in \operatorname{supp}\hat{\omega}_{p}$. We apply the logarithmic-type HLS inequality \eqref{eq:logHLS} for $M = D = 1$ and $\varepsilon = p-1$ to estimate
\begin{equation*}
E_{\textrm{log}}[\hat{\omega}_p] \le \frac{1}{2(p-1)}\log{N} + \frac{\log \pi}{2} -   \frac{\log(2-p) }{2(p-1)},\quad N=\int_{\mathbb{R}^{2}}\hat{\omega}_p^{p}dX.
\end{equation*}
We use the Casimir constraint
\[
1 = Q_p^{p}[\zeta_p] 
=\frac{1}{(2\pi)^{p-1}R_p^{2p}}\int_{\bbR^2} \frac{1}{X_2^{p-1}} \hat{\omega}_p^p(X) \,dX.
\]
By \eqref{eq:p1_W_supp_1}, we have
\begin{align*}
|X_2-R_p|\le Cd_p,\quad X\in \operatorname{supp}\hat{\omega}_p. 
\end{align*}
This implies that 
\begin{align*}
\left|\frac{1}{X_2^{p-1}}-\frac{1}{R_p^{p-1}} \right|\le Cd_p,\quad X\in \operatorname{supp}\hat{\omega}_p.
\end{align*}
It follows that 
\begin{align*}
\left|(2\pi)^{p-1}R_p^{2p}-\frac{1}{R_p^{p-1}}N \right|
=\left|\int_{\mathbb{R}^{2}}\left(\frac{1}{X_2^{p-1}}-\frac{1}{R_p^{p-1}} \right)\hat{\omega}_p^{p}dX \right|
\le Cd_pN. 
\end{align*}
This yields
\begin{align*}
\frac{1}{p-1}\log(1-Cd_pR_p^{p-1})
\le \frac{1}{p-1}\log{N}-\log{2\pi}-\left(3+\frac{2}{p-1}\right)\log{R_p}
\le \frac{1}{p-1}\log(1+Cd_pR_p^{p-1}).
\end{align*}
Both bounds are $o(1)$ and hence \eqref{eq:logintomegap} follows. Using 
\[
-\frac{\log(2-p)}{(p-1)} =  1 + O(p-1),
\]
we obtain \eqref{eq:p1_next_log}.
\end{proof}

\begin{lemma}\label{thm:p1_O(1)}
As $p \to 1^{+}$,
\begin{align}
I_p &= \frac{1}{24\pi^2e(p-1)} + \frac{1 }{ 8 \pi^2 e}\left(\log8\sqrt2 \pi - \frac{2}{3}\right) + o(1),
\label{eq:p1_en_O(1)}\\[4pt]
W_p &= \frac{1}{24\pi^2e(p-1)} + \frac{1 }{ 8 \pi^2 e}\left(\log8\sqrt2 \pi + \frac{5}{6}\right) + o(1).
\label{eq:p1_speed_O(1)}
\end{align}
Moreover, we have
\begin{align}
&I_p = \frac{1}{8\pi^2R_p^3} \left( \left( \frac{1}{p-1} + \frac{5}{2} \right) \log R_p + \log(8\sqrt2 \pi) - \frac{3}{2} \right) + o(1), \label{eq:p1_upper_2}\\
&E_{\textrm{log}}[\hat{\omega}_p]-\left( \frac{1}{p-1} + \frac{3}{2} \right) \log R_p -  \log(\sqrt{2} \pi) - \frac{1}{2} = o(1), \label{eq:logEasym} \\
&U[\hat{\omega}_p]+\left(\frac{2}{p-1}+3\right)\log{R_p}+\log{2\pi}=o(1).  \label{eq:Uasymp}
\end{align}
\end{lemma}

\begin{proof}
The asymptotic expansion \eqref{eq:p1_speed_O(1)} follows from \eqref{eq:p1_en_O(1)} and \eqref{eq:speed} since 
\[
W_p = \frac{14p-12}{5p-3} I_p = \left( 1 + \frac{9}{2}(p-1) + O((p-1)^2) \right)I_p.
\]
We expand $I_p = E[\zeta_p]$ for $\zeta_p = \Gamma_p \hat{\omega}_p/r$ by  
\begin{equation*}
I_p = \frac{\Gamma_p^2}{2} \iint_{\bbR^2 \times \bbR^2} \sqrt{X_2Y_2} F(s) \hat{\omega}_p(X) \hat{\omega}_p(Y) \,dXdY, \quad s = \frac{|X-Y|^2}{X_2Y_2}.
\end{equation*}
Using \eqref{eq:p1_F_log} and
\[
\log\frac8{\sqrt s} = \log\frac1{|X-Y|} + \log(8R_p) + \log\frac{\sqrt{X_2Y_2}}{R_p},
\]
we expand the kernel as 
\begin{align*}
\sqrt{X_2Y_2}F(s) = & R_p \left(\log\frac1{|X-Y|}+\log(8R_p) - 2\right)\\
& + (\sqrt{X_2Y_2} - R_p) \left(\log\frac1{|X-Y|} + \log(8R_p) - 2\right)\\
& + \sqrt{X_2Y_2} \log\frac{\sqrt{X_2Y_2}}{R_p} + O\left(\sqrt{X_2Y_2} s|\log s|\right).
\end{align*}
Substituting this expansion into the above and using $\Gamma_p =(2\pi R_p^{2})^{-1}$, we obtain
\[
\begin{aligned}
I_p 
&= \frac{\Gamma_p^2R_p}{2} \iint_{\mathbb{R}^{2}\times \mathbb{R}^{2} } \left(\log\frac1{|X-Y|}+\log(8R_p)-2\right) \hat{\omega}_p(X)\hat{\omega}_p(Y)\,dXdY\\
& +  \frac{\Gamma_p^2}{2} \iint_{\mathbb{R}^{2}\times \mathbb{R}^{2} }(\sqrt{X_2Y_2}-R_p) \left(\log\frac1{|X-Y|}+\log(8R_p)-2\right) \hat{\omega}_p(X)\hat{\omega}_p(Y)\,dXdY \\
&  +  \frac{\Gamma_p^2}{2} \iint_{\mathbb{R}^{2}\times \mathbb{R}^{2} }\sqrt{X_2Y_2}\log\frac{\sqrt{X_2Y_2}}{R_p} \hat{\omega}_p(X)\hat{\omega}_p(Y)\,dXdY \\
& + \frac{\Gamma_p^2}{2} \iint_{\mathbb{R}^{2}\times \mathbb{R}^{2} }\sqrt{X_2Y_2} \left(F(s)-\log\frac8{\sqrt s}+2\right) \hat{\omega}_p(X)\hat{\omega}_p(Y)\,dXdY \\
& =: \frac{1}{8\pi^2R_p^3} \left(E_{\textrm{log}}[\hat{\omega}_p]+\log(8R_p)-2\right) +\mathcal E_{1,p}+\mathcal E_{2,p}+\mathcal E_{3,p}.
\end{aligned}
\]
By the support bound \eqref{eq:p1_W_supp_1}, for $X, Y \in \supp \hat{\omega}_p$, we have
\begin{equation*}
\left|\sqrt{X_2Y_2} - R_p\right| + \left|\log\frac{\sqrt{X_2Y_2}}{R_p}\right| \lesssim d_p, \quad s|\log s|\lesssim \frac{d_p^{2}}{p-1}.
\end{equation*}
Since $\Gamma_p$ and $\log(8R_p) $ are uniformly bounded, it follows that 
\begin{align*}
|\mathcal E_{1,p}| &\lesssim d_p \iint_{\mathbb{R}^{2}\times \mathbb{R}^{2} }\left(1 + \log\frac1{|X-Y|}\right) \hat{\omega}_p(X)\hat{\omega}_p(Y)\,dXdY=d_p(1+E_{\textrm{log}}[\hat{\omega}_p]) \lesssim  \frac{d_p}{p-1},\\
|\mathcal E_{2,p}| &\lesssim d_p \iint_{\mathbb{R}^{2}\times \mathbb{R}^{2} } \hat{\omega}_p(X)\hat{\omega}_p(Y)\,dXdY \lesssim d_p,\\
|\mathcal E_{3,p}| &\lesssim \frac{d_p^{2}}{p-1}.
\end{align*}
Combining these estimates,
we obtain
\begin{equation*}
\begin{aligned}
\left| I_p-\frac{1}{8\pi^2R_p^3} \left(E_{\textrm{log}}[\hat{\omega}_p]+\log(8R_p)-2\right) \right| \lesssim \frac {d_p}{p-1}= o(1).
\end{aligned}
\end{equation*}
Therefore, \eqref{eq:p1_next_log} yields
\begin{align*}
I_p &\le \frac{1}{8\pi^2R_p^3} \left( \left(\frac{1}{p-1}+\frac{5}{2}\right)\log R_p + \log(8 \sqrt{2}\pi) - \frac{3}{2} \right) + o(1) \\
&\le \sup_{R>0} \frac{1}{8\pi^2R^3} \left( \left(\frac{1}{p-1}+\frac{5}{2}\right) \log R + \log(8\sqrt2\pi) - \frac{3}{2} \right) + o(1).
\end{align*}
The supremum is attained at
\[
\log R = \frac{1}{3} - \frac{\log(8\sqrt2\pi)-3/2}{(p-1)^{-1}+5/2}.
\]
Evaluating this yields
\[
I_p \le\frac{1}{24\pi^2e(p-1)} + \frac{\log(8\sqrt2 \pi) - 2/3}{8\pi^2e}+o(1).
\]
By the lower bound \eqref{eq:p1_W_low}, we obtain \eqref{eq:p1_en_O(1)}. Consequently, the asymptotic identity \eqref{eq:p1_upper_2} follows. The asymptotics of the logarithmic interaction energy \eqref{eq:logEasym} follow from \eqref{eq:p1_upper_2} and the above difference estimate for $I_p$. The asymptotic formula \eqref{eq:Uasymp} is equivalent to \eqref{eq:logEasym} by \eqref{logenergy}.
\end{proof}

\begin{lemma}
As $p\to 1^{+}$, 
\begin{align}
&S[\hat{\omega}_p]+\log{2\pi}+\left(3+\frac{2}{p-1}\right)\log{R_p}=o(1),  \label{eq:entropyasym}\\
&H[\hat{\omega}_p]=o(1).\label{eq:Hasymp}
\end{align}
\end{lemma}

\begin{proof} 
By Jensen's inequality,
\begin{align*}
e^{-(p-1)S[\hat{\omega}_p]}=\exp\left(\int_{\mathbb{R}^{2}}(\log{\hat{\omega}_p^{p-1}})\hat{\omega}_pdX\right)\le \int_{\mathbb{R}^{2}}\exp(\log{\hat{\omega}_p^{p-1}})\hat{\omega}_pdX=\int_{\mathbb{R}^{2}}\hat{\omega}_p^{p}dX.
\end{align*}
Taking the logarithm, we have 
\begin{align*}
J_p[\hat{\omega}_p]:=\frac{1}{p-1}\log\left(\int_{\mathbb{R}^{2}}\hat{\omega}_p^{p}dX\right)+S[\hat{\omega}_p]\ge 0.
\end{align*}
We combine this with the log-HLS inequality \eqref{eq:HlogHLS}. By the definition of the logarithmic interaction energy \eqref{logenergy}, we have  
\begin{align*}
0\le H[\hat{\omega}_p]+J_p[\hat{\omega}_p]
=U[\hat{\omega}_p]+\frac{1}{p-1}\log\left(\int_{\mathbb{R}^{2}}\hat{\omega}_p^{p}dX\right)
=-2E_{\textrm{log}}[\hat{\omega}_p]+1+\log\pi+\frac{1}{p-1}\log\left(\int_{\mathbb{R}^{2}}\hat{\omega}_p^{p}dX\right).
\end{align*}
By the asymptotics \eqref{eq:logEasym} and \eqref{eq:logintomegap}, the right-hand side is $o(1)$ as $p\to 1^{+}$. We obtain \eqref{eq:entropyasym} and \eqref{eq:Hasymp}.
\end{proof}

\begin{proposition}\label{rmk:further}
Let $\epsilon_p=\epsilon[\hat{\omega}_p]$. As $p\to 1^{+}$, we have 
\begin{align}
R_p &= e^{1/3} \left[ 1 + \left( \frac{3}{2} - \log(8\sqrt{2} \pi) \right)(p-1) + o(p-1) \right], \label{eq:W_radius}\\[4pt]
\Gamma_p &= \frac{1}{2\pi e^{2/3}} \left[ 1 + \left(2 \log(8\sqrt{2}\pi) - 3\right)(p-1) + o(p-1) \right], \label{eq:W_circulation}\\[4pt]
\epsilon_p &= 8 e^{-3}  \exp\left( -\frac{1}{3(p-1)} \right)(1 + o(1)). \label{eq:W_core}
\end{align}
\end{proposition}

\begin{proof}
Multiplying the Euler--Lagrange equation $\kappa_p \zeta^{p-1} =(\mathcal{G}\zeta-W_pr^{2}/2)_+$ of $\zeta\in S_p$ by $\partial_r \zeta$ and integrating over $\mathbb{R}^{2}_{+}$ with respect to $2\pi rdzdr$, we obtain 
\begin{align*}
\kappa_p (\zeta^{p-1},\partial_r \zeta)=(\mathcal{G}\zeta,\partial_r \zeta)-\frac{W_p}{2}(r^{2},\partial_r \zeta),
\end{align*}
where $(\cdot,\cdot)$ is defined in \eqref{eq:product}.
By the conditions $P[\zeta]=Q_{p}^{p}[\zeta]=1$ and $r-R_p=O(d_p)$ for $(z,r)\in \supp\zeta$, 
\begin{align*}
(r^{2},\partial_r \zeta) = -\frac{3}{R_p} + O\left(d_p\right),\quad 
(\zeta^{p-1},\partial_r \zeta) = -\frac1{pR_p}+O\left(d_p\right).
\end{align*}
Using the asymptotics of $F'(s)$ in \eqref{eq:p1_F_log}, for $s=|X-Y|^2/(X_2Y_2)$, we have 
\begin{align*}
(\mathcal{G}\zeta,\partial_r \zeta) &= -\frac{\Gamma_p^2}{2} \iint_{\mathbb{R}^{2}\times \mathbb{R}^{2}}\frac{X_2+Y_2}{\sqrt{X_2Y_2}} \left(\frac32F(s)-sF'(s)\right) \hat{\omega}_p(X)\hat{\omega}_p(Y)\,dXdY
= -\frac{3I_p}{R_p} - \frac{\Gamma_p^2}{2} + O \left(d_p\right).
\end{align*}
Substituting this into the above yields
\begin{align*}
\frac{1}{8\pi^{2}R_p^{3}}=\frac{6(p-1)}{5p-3}I_p+o((p-1)).
\end{align*}
Combining this with \eqref{eq:p1_upper_2}, we obtain 
\begin{align*}
\log{R_p}=\frac{1}{3}+\left(\frac{3}{2}-\log{8\sqrt{2}\pi} \right)(p-1)+o((p-1)).
\end{align*}
This implies \eqref{eq:W_radius}. Then \eqref{eq:W_circulation} and \eqref{eq:W_core} follow from $\Gamma_p = (2\pi  R_p^2)^{-1}$ and \eqref{eq:entropyasym}.
\end{proof}

\begin{proof}[Proof of Theorem \ref{t:MT5}]
The asymptotics \eqref{eq:Rp_asymptotics} and \eqref{eq:GRe} follow from \eqref{eq:W_radius} and \eqref{eq:W_core}. The convergence \eqref{eq:stabilityapplogHLS} follows from the stability estimate \eqref{eq:logHLSstability} and the convergence of the free energy \eqref{eq:Hasymp}. Using $(2\pi R_p^2)^{-1} = \Gamma_p$, we express \eqref{eq:p1_upper_2} as 
\begin{align*}
I_p = \frac{\Gamma_p}{4\pi R_p} \left( \left( \frac{1}{p-1} + \frac{5}{2} \right) \log R_p + \log(8\sqrt2 \pi) - \frac{3}{2} \right) + o(1).
\end{align*}
By \eqref{eq:scalefactor} and \eqref{eq:entropyasym}, it follows that 
\begin{align*}
\left(\frac{3}{2}+\frac{1}{p-1}\right)\log{R_p}=-\log\epsilon_p-1-\frac{1}{2}\log{2\pi^{2}}+o(1).
\end{align*}
Substituting this into the above, we obtain
\[
I_p = \frac{\Gamma_p}{4\pi R_p} \left( \log\frac{8R_p}{\epsilon_p} - \frac{5}{2} + o(1) \right).
\]
Since $-(p-1)\log\epsilon_p=1/3+O((p-1))$, it follows that 
\begin{align*}
    W_p &= \frac{14p-12}{5p-3} I_p = \left( 1 + \frac{9}{2}(p-1) + O((p-1)^2)\right)I_p\\
    &= \frac{\Gamma_p}{4\pi R_p} \left( \log \frac{8R_p}{\epsilon_p} - \frac{5}{2} + \frac{9}{2}(p-1)\log\frac{8R_p}{\epsilon_p} + o(1) \right)
    = \frac{\Gamma_p}{4\pi R_p} \left( \log\frac{8R_p}{\epsilon_p} - 1 + o(1) \right).
\end{align*}
We obtain \eqref{eq:W_Kelvin0} and the proof is now complete.
\end{proof}

\noindent\textbf{Acknowledgments.} 
KW is deeply grateful to In-Jee Jeong for many illuminating discussions and valuable suggestions.
KA was supported by the JSPS through the Grant-in-Aid for Scientific Research (C) 24K06800, MEXT Promotion of Distinctive Joint Research Center Program JPMXP0723833165, and Osaka Metropolitan University Strategic Research Promotion Project (Development of International Research Hubs). 
KW was supported by SNF Project 212573 FLUTURA -- Fluids, Turbulence, Advection.

\medskip
	
\noindent\textbf{Conflict of interest statement.} The authors declare no conflicts of interest.
	
\medskip
	
\noindent\textbf{Data availability statement.} Data sharing is not applicable to this article.

\medskip

\noindent\textbf{AI declaration.}
All mathematical content of this paper is the work of the authors.
ChatGPT (GPT-6) was used solely for English editing and assistance in preparing Figures 1--4 using Python and Ipe.
All AI-assisted code and figures were reviewed and verified by the authors.

\bibliographystyle{alpha-article-year}
\bibliography{Casimir}

\end{document}